\def\be{\begin{eqnarray}}
\def\en{\end{eqnarray}}
\def\BE{\begin{equation}}
\def\EE#1{\label{#1}\end{equation}}
\def\T{{\mathbb T}}
\def\subs#1#2{\mbox{\small${#1\atop#2}$}}
\def\nL{\nabla^{\rm L}}
\def\rf#1{(\ref{#1})}
\def\i{{\rm i}}
\def\e{{\rm e}}

\def\bv{{\bm v}}
\def\bx{{\bm x}}

\def\bxi{{\bm\xi}}
\def\ba{{\bm a}}
\def\bk{{\bm k}}
\def\bom{{\bm\omega}}
\def\oz{{\bm\omega}^{\rm(init)}}
\def\vz{{\bm v}^{\rm(init)}}
\def\tf{t_{\rm fin}}
\def\Rlb{{R_{\rm bound}}}
\def\Rei{{R_{\rm inf}}}
\documentclass[12pt,5p]{Melsarticle}
\usepackage{amsfonts,amssymb,amsmath,epsfig,bm,url,wrapfig}
\begin{document}
\title{THE CAUCHY--LAGRANGIAN METHOD\\FOR NUMERICAL ANALYSIS OF EULER FLOW}
\author[mitpan]{O.~Podvigina}
\author[mitpan]{V.~Zheligovsky}
\author[oca]{U.~Frisch}
\address[mitpan]{Institute of Earthquake Prediction Theory and
Mathematical Geophysics, Russian Acad.~Sci.,\\84/32 Profsoyuznaya St, 117997 Moscow, Russian Federation}
\address[oca]{Lab.~Lagrange,~UCA,~OCA,~CNRS,~CS~34229,~06304~Nice~Cedex~4,~France}

\begin{abstract}
A novel semi-Lagrangian method is introduced to solve numerically the Euler
equation for ideal incompressible flow in arbitrary space dimension.
It exploits the time-analyticity of fluid particle trajectories and requires,
in principle, only limited spatial smoothness of the initial data. Efficient
generation of high-order time-Taylor coefficients is made possible
by a recurrence relation that follows from the Cauchy invariants formulation
of the Euler equation (Zheligovsky \& Frisch, J.~Fluid Mech.~2014, {\bf 749},
404--430). Truncated time-Taylor series of very high order allow the use
of time steps vastly exceeding the Courant--Friedrichs--Lewy limit, without
compromising the accuracy of the solution. Tests performed on the
two-dimensional Euler equation indicate that the Cauchy--Lagrangian method is
more --- and occasionally much more --- efficient and less prone to instability
than Eulerian Runge--Kutta methods, and less prone to rapid growth
of rounding errors than the high-order Eulerian time-Taylor algorithm. We also
develop tools of analysis adapted to the Cauchy--Lagrangian method, such as
the monitoring of the radius of convergence of the time-Taylor series. Certain
other fluid equations can be handled similarly.
\end{abstract}
\maketitle

\section{Introduction}\label{intro}

As is well known, fluid flow can be characterised in terms of the
current positions of fluid particles (Eulerian coordinates), or in terms
of their initial positions (Lagrangian coordinates). For ideal (inviscid)
incompressible fluid both formulations were introduced in the 18th century
\cite{eu,la}. From both a theoretical and a numerical point of view, the
Eulerian formulation seems to have a significant edge because it gives an
explicit quadratic expression for the time-derivative of the velocity;
the Lagrangian formulation has even been
qualified ``an agony'', because of its complexity \cite{pr}.

However, the Eulerian time-stepping methods have one serious well-known
drawback, the Courant--Friedrichs--Lewy (CFL) condition \cite{cfl}, which
constrains the time step to be less than a dimensionless constant multiplied by
the time needed to sweep across the spatial mesh at the maximum flow speed. As
a consequence, the complexity of computations with
$N$ collocation points in each spatial direction is roughly O$(N^4)$ in three
dimensions. Hence, progress in high-resolution numerical simulation is
creepingly slow, in spite of Moore's law. (As we will see in
Section~\ref{s:rounding}, this is not the only drawback of Eulerian schemes.)

Foremost because of the CFL constraint, which may conflict, for example,
with the desire to quickly produce well-resolved numerical weather
forecast, there has been a strong incentive to develop
{\it semi-Lagrangian} (SL) schemes, in which some form of Lagrangian
integration -- with time steps not constrained by CFL -- alternates with a
reversion to an Eulerian grid (see, e.g., \cite{sc} and references therein).
The SL algorithms used so far, e.g., in geophysical fluid
dynamics, engineering, mechanics and plasma physics, were designed for
situations where satisfactory results can be obtained with rather
low-order temporal schemes. They are thus not appropriate for
numerical experimentation on delicate questions, such as the issue of
blow-up in three-dimensional flow (see, e.g., \cite{gi}).

We propose a new SL algorithm, which we call the Cauchy--Lagrangian
algorithm. It relies on Cauchy's Lagrangian formulation of the
equations of ideal incompressible flow \cite{cau,fv} and on recent
results about the time-analyticity of Lagrangian
trajectories \cite{fz,zf}. The Cauchy--Lagrangian algorithm is
particularly well-suited for problems where high precision is a
prerequisite, and it is actually superior to Eulerian schemes.

Cauchy's 1815 Lagrangian equations are formulated in terms of the Lagrangian
map from Lagrangian to Eulerian coordinates, $\ba\to\bx=\bx(\ba,t)$, defined as
the solution of the ordinary differential equation for fluid particle trajectory,
$\dot{\bx}=\bv(\bx,t)$, with the initial condition $\bx(\ba,0)=\ba$. Here,
the dot denotes a (Lagrangian) time derivative. The three-dimensional Cauchy
invariants equations are (see Section~\ref{s:background} for a simplified
derivation):
\BE\sum_{k=1}^3\nL\dot{x}_k\times\nL x_k=\oz,\qquad\det(\nL\bx)=1,\EE{caueq}
where $\oz=\nL\times\vz$ denotes the initial vorticity and $\nL$ is
the Lagrangian gradient, i.e., the space derivatives in $\ba$.
Since the r.h.s.~of the first equation in \rf{caueq} does not depend on time,
its l.h.s.~is obviously a Lagrangian invariant,
whose scalar components are called the ``Cauchy invariants''. These
invariants were much later interpreted as a consequence of Noether's
theorem applied to a continuous symmetry of the Euler
equation, called the relabelling invariance (see \cite{fv}, Section 5.2).

It was shown in \cite{fz,zf} that the Cauchy invariants equations
\rf{caueq} imply analyticity of fluid trajectories for initial
conditions that have certain rather weak regularity. The proof of
analyticity is derived from an explicit recurrence relation for coefficients
of the time-Taylor series
for the Lagrangian displacement $\bxi=\bx-\ba$. It can also be used to
construct a numerical Lagrangian method of very high order in time in
both two dimensions (2D) and three dimensions (3D). Its time steps are
only constrained by the radius of convergence of the Taylor series
(typically, the inverse of the largest initial velocity gradient). For more
details, see Section~\ref{s:background} and \cite{fz,zf}.

The paper is organised as indicated hereafter. In
Section~\ref{s:background} we recall some of the known results about
Cauchy's Lagrangian formulation and its application to the
time-analyticity of the Lagrangian map. In Section~\ref{s:CLalgorithm}
we describe the Cauchy--Lagrangian (CL) algorithm in detail: we begin
with an overview and show that CL may be considered as a
semi-Lagrangian algorithm of arbitrary high order
(Section~\ref{ss:overview}), present the interpolation technique for
reverting to Eulerian coordinates at the end of each time step
(Section~\ref{ss:interpolation}) and show how to optimise the choices
of the time step and of the order of the Taylor expansion truncation
to minimise computational complexity (Section~\ref{ss:optimal}).
Section~\ref{s:testing} is devoted to testing the CL algorithm in the
2D case against various numerical methods, and to CPU
benchmarks. Section~\ref{s:rounding} is a comparison of high-order
time-Taylor expansions in Eulerian and Lagrangian coordinates: the
Lagrangian method suffers much less from the rounding errors. In
Section~\ref{s:depletion} we determine how quickly we have to decrease
the time step in the CL method because the radius of convergence
$R(t)$ of the time-Taylor series around time $t$ generally shrinks
as $t$ increases, as more and more small-scale eddies (high spatial
Fourier harmonics) are generated. Finally, in Section~\ref{s:conclusion},
we present concluding remarks and point out that the Cauchy--Lagrangian method
is well adapted for certain other problems concerning ideal fluid flow.

\section{Mathematical background}\label{s:background}

Here we recall some of the background material which is used to
develop the Cauchy--Lagrangian method: derivation of the Cauchy
invariants equations and the recurrence relation for time-Taylor
coefficients from which follows the analyticity of the time-Taylor
series (Section~\ref{ss:Lagfo}). We then present, in a
heuristic way, a result allowing the determination of the radius of
convergence of such Taylor series, an important new tool for analysing
the output of Cauchy--Lagrangian computations (Section~\ref{ss:lpth}).

\subsection{Cauchy's Lagrangian formalism for ideal incompressible fluid flow}
\label{ss:Lagfo}

The flow of an ideal incompressible fluid, described in Eulerian
coordinates, is governed by the Euler equation
\BE\partial_t\bv+(\bv\cdot\nabla)\bv=-\nabla p,\qquad\nabla\cdot\bv=0,\EE{euler}
where $\bv(\bx,t)$ is the velocity and $p(\bx,t)$ the pressure (divided
by the density, which in the incompressible case is constant). Here the spatial
differentiation $\nabla$ is performed in the Eulerian coordinates~$\bx$.

The Eulerian equation \rf{euler} is notorious for difficulties in its
investigation, both analytically and numerically. It turns out that
its Lagrangian analogue can be once integrated in time, which yields
the Cauchy invariants equations \rf{caueq}. Since these equations are
central to our numerical method, we briefly recall how they are derived.

In Lagrangian coordinates, fluid particles are characterised by their initial
positions $\ba$ and the time $t$. The subsequent positions of fluid
particles are then described by the Lagrangian map $\ba\mapsto\bx(\ba,t)$.
The spatial gradient in the Lagrangian coordinate $\ba$ is denoted $\nL$.

Cauchy's derivation of \rf{caueq} is now briefly recalled \cite{cau}.
Observe that the l.h.s.~of \rf{euler} is the acceleration $\ddot{\bx}$
of the fluid particle; thus, $\ddot{\bx}=-\nabla p$. Multiplying this
equation by the transpose of the Jacobian matrix $\nL\bx$,
we transform the Eulerian gradient in the r.h.s.~into the Lagrangian one and
obtain, in the component form:
\BE\sum_{k=1}^3\ddot{x}_k\nL_ix_k=-\nL_ip.\EE{laggrad}
Cauchy then applies a Lagrangian curl to \rf{laggrad},
and finds that the l.h.s. can be exactly integrated in
time to yield the first equation of the set \rf{caueq}. The second one,
the Jacobian equation, expresses the conservation
of volume, equivalent to incompressibility.

The Cauchy invariants equations are quadratically nonlinear equations,
not resolved in the time derivatives; the Jacobian equation
is also quadratic (in 2D) or cubic (in 3D). At first glance, it is unclear,
in what respect \rf{caueq} are advantageous compared to the Eulerian equations
\rf{euler}. To see what is gained, following \cite{fz,zf}, we consider
the equations for the displacement $\bxi=\bx-\ba$, that ensue from \rf{caueq}
and are in 3D as follows:
\begin{align}
\nL\times\dot{\bxi}&+\sum_{k=1}^3\nL\dot{\xi}_k\times\nL\xi_k=\oz,\label{drot}\\
\nL\cdot\bxi&+\sum_{1\le i<j\le3}(\nL_i\xi_i\,\nL_j\xi_j
-\nL_i\xi_j\,\nL_j\xi_i)+\det(\nL\bxi)=0,\label{dgra}
\end{align}
and expand the displacement in the time-Taylor series
\BE\bxi(\ba,t)=\sum_{s=1}^{\infty}\bxi^{(s)}(\ba)t^s.\EE{xiSeries}
The structure of equations \rf{drot}--\rf{dgra} enables us
to derive recurrence relations for the coefficients $\bxi^{(s)}(\ba)$,
which in 3D take the following form \cite{fz,zf}:
\begin{align}
\nL\!\times\bxi^{(s)}=&\,\delta^s_1\oz-\sum_{k=1}^3~\sum_{m=1}^{s-1}\,
{m\over s}\,\nL\xi^{(m)}_k\times\nL\xi^{(s-m)}_k,\label{rrot}\\
\nL\cdot\bxi^{(s)}=&\sum_{1\le i<j\le3}~\sum_{m=1}^{s-1}\,\left(
\nL_j\xi^{(m)}_i\nL_i\xi^{(s-m)}_j-\nL_i\xi^{(m)}_i\nL_j\xi^{(s-m)}_j\right)\nonumber\\
&-\sum_{i,j,k}~\sum_{l+m+n=s}\!\!\varepsilon_{ijk}
\nL_i\xi^{(l)}_1\nL_j\xi^{(m)}_2\nL_k\xi^{(n)}_3.\label{rgra}
\end{align}
Here $\delta^s_1$ is the Kronecker delta and $\varepsilon_{ijk}$
is the totally antisymmetric tensor.
In particular, for $s=1$, \rf{rrot}--\rf{rgra}
are equivalent to the relation $\bxi^{(1)}=\vz$.

To determine $\bxi^{(s)}$, knowing its (Lagrangian) curl and divergence,
is a Helmholtz--Hodge problem, with a unique solution
for suitable boundary conditions. For example, for space-periodic flow
we demand that the averages of all $\bxi^{(s)}$
over the periodic box vanish, and construct the solution of
\rf{rrot}--\rf{rgra} by taking the gradient of \rf{rgra},
subtracting from it the curl of \rf{rrot} and thus obtaining a Poisson equation
for $\bxi^{(s)}$. This yields the following single (tensorial)
recurrence relation for the (Lagrangian) gradients of the time-Taylor
coefficients ($s\ge1$)
\be&&\hspace*{-15mm}\nL_\mu\xi_\nu^{(s)}=\nL_\mu v_\nu^{\rm (init)}\delta^s_1+\sum_{\subs{1\le j\le3,}{j\ne\nu}}
{\cal C}_{\mu j}\left(\,\sum_{\subs{1\le k\le3,}{0<m<s}}{2m-s\over s}\,
(\nL_\nu\xi^{(m)}_k)\,\nL_j\xi^{(s-m)}_k\right)\label{bigXis}\\
&&\hspace*{-15mm}+\,{\cal C}_{\mu\nu}\left(\sum_{\subs{1\le i<j\le3}{0<m<s}}
\left(\rule{0mm}{1em}(\nL_j\xi^{(m)}_i)\nL_i\xi^{(s-m)}_j-(\nL_i\xi^{(m)}_i)
\nL_j\xi^{(s-m)}_j\right)-\!\!\sum_{\subs{i,j,k}{l+m+n=s}}
\varepsilon_{ijk}(\nL_i\xi^{(l)}_1)(\nL_j\xi^{(m)}_2)\nL_k\xi^{(n)}_3\right).
\nonumber\en
Here, ${\cal C}_{ij}\equiv\nabla^{-2}\nL_i\nL_j$ and $\nL_i\nL_j$ denotes
the second-order partial derivative $\partial^2/\partial a_i\partial a_j$.
The operator $\nabla^{-2}$ is the inverse of the (Lagrangian)
Laplacian: given a periodic function $f$ with zero spatial mean over the periodic
cell, $\nabla^{-2}f$ is defined as the unique periodic function $\psi$ with
zero spatial mean, solving $\nabla^2\psi=f$.

The explicit character of the recurrence relation \rf{bigXis}
permits us to derive bounds for the
time-Taylor coefficients, provided the initial condition
possesses a certain minimum smoothness, e.g., if the initial vorticity is
H\"older-continuous (i.e., belongs to the space C$^\alpha$ for some
$0<\alpha<1$) or belongs to the space of absolutely summable Fourier series
(henceforth, ASFS). This way one can prove various analyticity theorems
\cite{fz,zf}, which state that the Lagrangian gradient of the Taylor series
\rf{xiSeries} converges in the space to which the initial vorticity belongs,
provided $|t|<\Rlb$. Here, $|t|$ denotes the modulus of complex $t$ and $\Rlb$
\textit{a lower bound for the radius of analyticity}, which is typically
of the order of the inverse C$^\alpha$ or ASFS norm of the initial vorticity.
Within the disk of radius $\Rlb$, convergence of the time-Taylor series
\rf{xiSeries} and that of its spatial gradient is guaranteed and so is analyticity
in time. The interested reader will find elementary proofs in \cite{zf}. Of
course, this is a key result for the construction of numerical solutions.
Note that in 2D, by a theorem of Wolibner \cite{wo}, when the initial
vorticity is H\"older-continuous, the vorticity stays smooth
for all positive times. In 3D this problem is open.

Suppose that, by using the Taylor series around $t=0$ we have solved
\rf{caueq} in some real time interval $[0,t_1]$, where
$t_1<\Rlb$. If the solution at time $t_1$ has the required
smoothness, it can be used as a new initial condition
for \rf{caueq} at time $t_1$. For this, Lagrangian coordinates must be
reintroduced for fluid particles starting at $t_1$. At $t=t_1$ the new
Lagrangian coordinates coincide, of course, with Eulerian
coordinates. Hence we need to revert to Eulerian coordinates at $t_1$.
In principle, this can be done by composing the Lagrangian fields in
the original coordinates with the inverse of the Lagrangian map, but
there are other ways, as we will see.

This procedure can be extended
into a Weierstrass-type analytic continuation: one repeats the process at
times $t_1<t_2\ldots$, as long as the successively
constructed solutions have enough spatial smoothness to give time-Taylor
series with a finite non-vanishing radius of convergence. In 2D, when the
initial vorticity is H\"older-continuous, the required smoothness will persist
forever. In 3D loss of smoothness (blow-up) cannot be ruled out.

\subsection{The radius of convergence of the time-Taylor series}\label{ss:lpth}

When the initial vorticity is bounded in a suitable function space,
such as C$^\alpha$ or ASFS, we know that, for
$|t|<\Rlb$, the gradient of the Taylor series \rf{xiSeries} is guaranteed to
converge in that space for all $\ba$. Convergence may hold over a
larger time interval, but without any proven control over the spatial
smoothness of the Lagrangian map thus constructed. Knowing the radius
of convergence of this Taylor series is important for numerical
applications, even if $|t|$ does not approach the edge of the disk of
convergence: this parameter controls how fast the terms in the
series tend to zero at large orders $s$.

For a given starting point $\ba$, there exists a non-negative number $R(\ba,0)$,
called the radius of convergence. It is the largest number such that the series
\rf{xiSeries} is absolutely convergent for any $0\le|t|<R(\ba,0)$.
The second argument of the radius $R$ indicates the time
around which the Taylor expansion is performed.
Note that the radius of convergence can be infinite.
It is given by the Cauchy--Hadamard formula (see, e.g., \cite{lang})
\BE{1\over R(\ba,0)}=\limsup_{s\to\infty}|\bxi^{(s)}(\ba)|^{1/s}
\equiv\lim_{S\to\infty}\,\sup_{S\ge s}|\bxi^{(S)}(\ba)|^{1/S}.\EE{Rofa}

It might be natural to define a radius of convergence $R(0)$ by taking the
infimum of $R(\ba,0)$ over all initial points $\ba$.
How can we measure this infimum in practice, given a set of Taylor
coefficients $\bxi^{(s)}(\ba)$ up to a large value of $s$, and
assuming that the $\ba$ values are on a large set of collocation
points? Following the usual numerical analyst's approach,
one might try to just evaluate the
local radii of convergence $R(\ba,0)$ on a grid in the fluid domain
and then find the minimum. This approach would have the advantage
that it might indicate, what region in space is characterised by small
values of the local radii of convergence, eventually perhaps leading
to blow-up. However, the techniques for estimating the
radius of convergence are numerically quite demanding, as
discussed in \ref{a:conver}, suggesting that this approach is not
numerically optimal. Furthermore, when the radius of convergence $R(\ba,0)$
is discontinuous, such an approach is clearly inappropriate.

Even if $\bxi^{(s)}(\ba)$ are analytic in $\ba$,
we cannot guarantee the continuity of $R(\ba,0)$
(the series $\sum_{s\ge0}\,2^{-|\ba-\ba_\star|^2s^2}\,t^s$, unrelated to the Euler
equation, illustrates this difficulty).
If $R(\ba,0)$ takes small values at exceptional off-grid points,
we will miss them when measuring $R(0)$ numerically on a discrete grid.
We find it therefore more suitable to employ the quantity
\BE\Rei(0)={\rm ess\,inf}_{\ba}R(\ba,0),\EE{yesinf}
where ``ess'' stands for ``essential'', namely, small outlier values
on a zero-measure set are disregarded.

A simple way to evaluate $\Rei(0)$ makes use of the ordinary
power series obtained by replacing $\bxi^{(s)}(\ba)$ by its L$^p$ norm, namely,
\BE\sum_{s=1}^\infty\|\bxi^{(s)}(\ba)\|_p\,t^s.\EE{Tp}
We denote by $R_p$ the radius of convergence of the series \rf{Tp}.
We assume that the initial vorticity is space-periodic
and is either H\"older-continuous or has an
absolutely summable Fourier series. The same properties hold then for all the
gradients of the Taylor coefficients $\nL\bxi^{(s)}(\ba)$. As a consequence,
all the Taylor coefficients $\bxi^{(s)}(\ba)$ are in the Lebesgue
spaces L$^p$ of functions, such that
the $p$th power of their absolute value is integrable. In \ref{a:lpth}
we prove a theorem stating, roughly, that the radius
$\Rei(0)$ is also the radius of convergence of the series \rf{Tp}.
More precisely,\\
\textit{1. For any $p\ge1$, we have $\Rei(0)\ge R_p$.\\
2. If the displacement $\bxi(\ba,t)$, defined by \rf{xiSeries}, is
in {\rm L}$^p$ for any $0\le t<\Rei(0)$, then $\Rei(0)=R_p$.}\\
(The proof of this theorem does not involve the Euler equation.)

For certain solutions of the Euler equation, the assumption in
the second part of the theorem is known to be satisfied. For
example, for 2D flow whose initial vorticity is H\"older-continuous,
Wolibner's theorem \cite{wo} implies that the Lagrangian displacement is always
in L$^p$. In general, for 3D Euler flow, we do not know the spatial regularity
properties of the displacement beyond the time of guaranteed analyticity and
thus we only know that $\Rei(0)\ge R_p$.

To gain heuristic insight into why the theorem holds, we observe
that, since the radius of
analyticity depends only on the modulus of the Taylor coefficients, no
generality is lost by taking them to be scalar and non-negative functions of
$\ba$. For such a Taylor series, with radius of convergence $R(\ba,0)$, the
coefficient of order $s$ is roughly proportional to $R(\ba,0)^{-s}$ (this
would be exact for a geometric series). To calculate the L$^p$
norm of the coefficient of order $s$, we must integrate the $p$th power
of this Taylor coefficient over $\ba$ and take the $p$th root.
It seems plausible
that, for large $s$, these integrals will be dominated by those $\ba$ where
$R(\ba,0)$ is close to its infimum. We thus expect
that the L$^p$ norm of the Taylor coefficient of order $s$ is roughly
proportional to $({\rm ess\,inf}_{\ba}R(\ba,0))^{-s}$, thereby
yielding a radius of convergence given by $\Rei$, namely, \rf{yesinf}.

\section{Description of the Cauchy--Lagrangian algorithm}\label{s:CLalgorithm}

In Section~\ref{ss:overview} we show that the basic ideas of the
Cauchy--Lagrangian (CL) algorithm are quite simple, given our
understanding of the analytical properties of the Lagrangian map,
which were recalled in Section~\ref{ss:Lagfo}. Then, in
Sections~\ref{ss:init}--\ref{ss:optimal} we successively examine the
constraints on the initial conditions (in connection with the spatial
truncation errors), the temporal truncation errors, the reversion from
the Lagrangian to the Eulerian representations and the
minimisation of computational complexity.

\subsection{Overview: CL as a high-order semi-Lagrangian method}\label{ss:overview}

Our purpose is to solve the Euler equation \rf{euler} with the initial
condition $\bv(\bx,0)\equiv\vz(\bx)$ or $\bom(\bx,0)\equiv\oz(\bx)$ in a time
interval $[0,\,T]$. In principle, we can just require that $\vz$ or $\oz$
be in a function space that guarantees finite-time analyticity of the
Lagrangian map, such as (for the vorticity) C$^\alpha$
or the space of absolutely summable Fourier series ASFS.
For practical reasons, more regularity is desirable (see Section~\ref{ss:init}).

\begin{figure}[t]
\centerline{\psfig{file=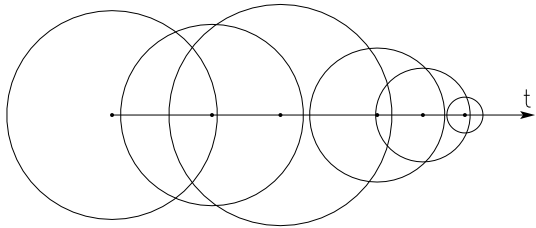,width=70mm}}
\vskip2mm
\caption{The successive disks of analyticity in the complex $t$ plane.}
\label{f:disks}\end{figure}

The simplest case is when $T$ is less than the radius of analyticity bound
$\Rlb(0)$ of the time-Taylor series of the Lagrangian map, starting at
$t=0$. We then use the time-Taylor expansion of the displacement
\rf{xiSeries}, together with the recurrence relation \rf{bigXis}
to calculate as many Taylor coefficients as needed (see
Sections~\ref{ss:temptrunc} and \ref{ss:optimal}). Here we just
observe that, for spatial periodicity conditions, it is straightforward
to calculate the Taylor coefficients using a pseudospectral method
with dealiasing. In 3D, we can then obtain the vorticity in Lagrangian
variables at time $T$ by applying Cauchy's vorticity formula \cite{cau}:
\BE\bom^{\rm L}(\ba,T)=\oz(\ba)\cdot\nL\bx(\ba,T).\EE{Cfrm}
(In 2D, this formula reduces to the constancy of the vorticity
$\bom^{\rm L}(\ba,T)=\oz(\ba)$.) Finally we can revert to Eulerian coordinates
by composing this Lagrangian vorticity with the inverse Lagrangian map
$\ba(T):~\bx\mapsto\ba(\bx,T)$, which maps the fluid particle position
$\bx$ at time $T$ to its initial position $\ba$, to obtain
\BE\bom(\bx,T)=\bom^{\rm L}(\ba(\bx,T),T).\EE{pedestrianbacktolabels}
Numerically, this is best done by interpolating from a Lagrangian
to an Eulerian grid (see Section~\ref{ss:interpolation}). Finally,
the Eulerian velocity at time $T$, if needed, can be obtained either
from the Eulerian vorticity at time $T$ (by inverting the curl, e.g.,
applying Fourier methods, or by interpolating to an Eulerian grid the
time derivative of the Lagrangian map, obtained by time
differentiation of the Taylor expansion \rf{xiSeries}).

Now, consider the case $T\ge\Rlb(0)$ when convergence of the Taylor series
is not guaranteed. Lack of convergence does not imply
loss of regularity at some {\it real} time $t$ between 0 and $T$. It may
just be due to complex-time singularities whose moduli are
less than or equal to $T$. Even if the Taylor series still converges, the sum
might be insufficiently smooth to allow starting a new Taylor
expansion (see end of Section~\ref{ss:Lagfo}). This is impossible in
2D with a H\"older-continuous initial vorticity, but cannot be
ruled out in 3D. The way to handle times $T$ that are beyond the
radius of analyticity bound $\Rlb(0)$, is the multistep method, which
essentially implements numerically the analytic continuation explained in
Section~\ref{ss:Lagfo}. We pick a time $t_1$ such that $0<t_1<\Rlb(0)$
and start the process all over at time $t_1$ using the Eulerian field
$\bv(\bx,t_1)$ as the initial velocity (or $\bom(\bx,t_1)$ as the initial
vorticity). We then have, in the complex time plane, a new disk of
guaranteed analyticity centred at $t_1$ with radius $\Rlb(t_1)$, which may well
extend on the real positive time axis beyond the furthest point reached
by the first disk of guaranteed analyticity. The process may be continued
(see Fig.~\ref{f:disks}) with analyticity disks centred at various
times $0<t_1<t_2<\ldots$, as long as the solution remains in a space
such as C$^\alpha$ or ASFS which guarantees that $\Rlb>0$.

\begin{wrapfigure}[30]{L}{.26\textwidth}
\vspace*{-1.4\baselineskip}
\begin{center}
\includegraphics[width=.28\textwidth]{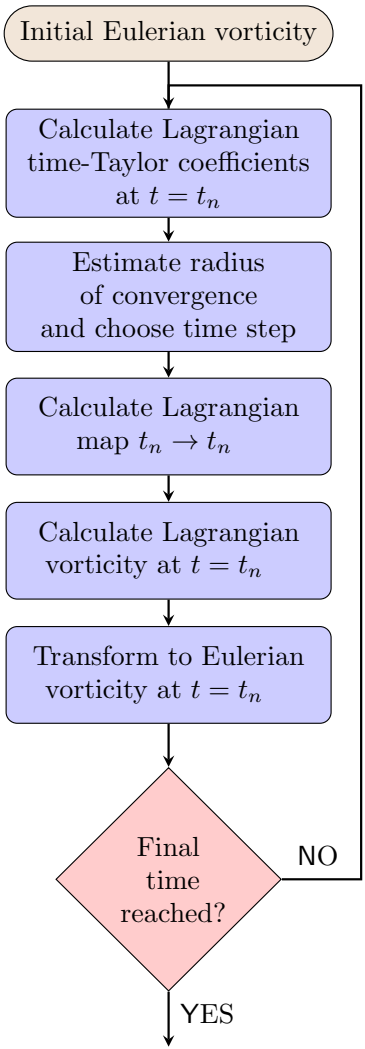}
\end{center}
\vspace*{-1.6\baselineskip}
\caption{A flowchart of the CL algorithm operation.}
\label{f:FC}\end{wrapfigure}

The algorithm outlined above is illustrated by the flow chart
in Fig.~\ref{f:FC}.

To finish this overview, we observe that our Cauchy--Lagrangian method
belongs to the family of semi-Lagrangian (SL)
methods, often employed, e.g., in meteorology, engineering, mechanics and
plasma physics \cite{bo,sc}. An SL method solves an evolution equation along
trajectories of particles for one time step and then reverts to the static
regular Eulerian grid, usually by interpolation. Unlike pure Eulerian ones,
SL methods are not subject to the CFL constraint
\cite{cfl}; consequently, their time steps are independent of the spatial resolution.
For time-analytic Lagrangian fluid particle trajectories, the time step must
just be smaller than the radius of convergence of the time-Taylor
series, which is typically the inverse of the maximum velocity gradient.

SL methods fall into two categories: upstream methods compute the so-called
departure points, i.e., the positions of Lagrangian particles which
are advected onto the regular grid during the time step, whereas for downstream
SL methods the departure points of trajectories are on the regular grid,
and the arrival points are on the distorted grid. The latter methods are
simpler, because they do not require solving auxiliary problems
to find departure points \cite{lp}. In this terminology, our CL method is
downstream. Standard SL methods are typically of second-order
accuracy in time; numerical schemes of higher orders are not widely used, being
too complicated. Indeed, to advance by one time step, a high-order SL method
of the traditional type would need information about the solution at
many time levels; this would require many
additional interpolations. In contrast, in the Cauchy--Lagrangian
approach, the recurrence relation \rf{bigXis} allows us to
construct very high-order schemes requiring only one
interpolation per time step.

\subsection{Spatial truncation errors and the choice of initial conditions}\label{ss:init}

As we will see in Sections~\ref{s:testing} and \ref{s:rounding}, our
CL method is particularly well adapted and efficient for problems where
high accuracy is important. Thus we want to keep all errors as (consistently)
small as possible. This includes of course spatial truncation errors, stemming
from the use of a finite number of modes. Actually, this
is much easier if we limit ourselves to \textit{initial conditions that are
analytic functions of the space variable.} Then, the Eulerian solution
remains analytic in space for at least a finite time (all times in 2D)
\cite{bb}. During that time interval, the spatial Fourier
transform of the solution falls off roughly as $\e^{-\delta(t)k}$,
where $k$ is the wave number and $\delta(t)$ is the radius of the tube
of analyticity, the distance between the real spatial domain
and the nearest complex-space singularity for the complexified spatial
variables. Hence, truncation errors fall off then as
$\e^{-\delta(t)k_{\max}}$, where $k_{\max}$ is the maximum wave number
present in the numerical truncation; thus they remain small as long as
$\delta(t)$ is significantly larger than the spatial numerical mesh \cite{ssl}.

In contrast, consider what happens when an initial condition has
merely the minimum regularity required for ensuring the time-analyticity
of the Lagrangian map, e.g., an initial vorticity is H\"older-continuous.
The spatial Fourier transforms of the velocity field during the time interval
$[0,T]$ will then fall off in a roughly algebraic way with the wave number
$k$. Thus spatial truncation errors will fall off only algebraically
with $k_{\max}$. Achieving high accuracy may then require spatial
resolutions that are beyond the reach of any existing computer. To avoid
this, in our test computations we limit ourselves to spatially
analytic initial data.

One of the simplest space-analytic initial conditions that can be used
for testing a CL code is the AB flow,
discussed in Section~2.6.1 of \cite{zf}, a steady-state solution of
the 2D Euler equation. Its vorticity is $\omega^{\rm(init)}=\sin x\cos y$.
Although its Eulerian-coordinates dynamics is trivial, its Lagrangian
dynamics is not completely trivial, but can still be determined
analytically using elliptic integrals. More examples of flows with
analytic initial data and truly non-trivial dynamics will be given in
Section~\ref{ss:flow-methods}.

\subsection{Temporal truncation errors}\label{ss:temptrunc}

In numerical applications the Lagrangian displacement is approximated by the time-Taylor series
\rf{xiSeries} truncated to some finite order $S$:
\BE\bxi_S(\ba,\Delta t)=\sum_{s=1}^S\bxi^{(s)}(\ba)(\Delta t)^s.\EE{xiSeriesS}
Which constraints on $\Delta t$ and on $S$ guarantee
that the remainder of the series is sufficiently small?

First, we demand that the series \rf{xiSeries} be convergent, i.e.,
$\Delta t$ should be less than the radius
of convergence for all starting points $\ba$, and hence
less than $\Rei$, the infimum over $\ba$
of these radii (ignoring sets of zero measure).

The determination of $\Rei$ proceeds as follows.
Using the material of Section~\ref{ss:lpth} and of \ref{a:lpth}, we
replace the full functional Taylor series by the series of the L$^p$ norms
of the Taylor coefficients. This gives a series with radius of
convergence $R_p$. In 2D, $R_p=\Rei$ and, in principle, in 3D we just have
$R_p\le\Rei$ --- however, since in 3D a blow-up of solutions
of the Euler equation in the absence of solid boundaries
perhaps never takes place,
it is conceivable, at least for space-analytic initial conditions, that $R_p=\Rei$.
The most convenient is to use the L$^2$ norm, which can be evaluated using
Parseval's theorem. We thus have to determine the radius of
convergence $R_2$ of a power series with positive coefficients.
How to do this in practice is discussed in \ref{a:conver}.

Second, we demand that the error due to truncation of the time-Taylor
series be small. To see what this implies, we consider the remainder
of the Taylor series and majorise its L$^2$ norm as follows:
\BE\|\bxi(\ba,t)-\bxi_S(\ba,t)\|\le
\sum_{s=S+1}^{\infty}\|\bxi^{(s)}\|(\Delta t)^s,\EE{majorL2}
where $\|\cdot\|$ denotes the L$^2$ norm. The series in the r.h.s.~has
positive coefficients and its sum has a singularity at $\Delta t=R_2$.
Thus the sum is typically well approximated (at large orders
and up to logarithmic factors) by a geometric series of ratio
$r\equiv\Delta t/R_2$, and the L$^2$ norm of the
remainder will typically be less than $\|\bxi^{(S)}\|(\Delta t)^Sr/(1-r)$.
For the values of the ratio $r$ that minimise the complexity of
the computation (see Section~\ref{ss:optimal}), $r/(1-r)\approx0.15$\,.
Hence, our criterion for keeping the truncation error small is
\BE\|\bxi^{(S)}\|(\Delta t)^S<\varepsilon.\EE{estet0}
That is, we demand that, in the time-Taylor series, the last term kept before
truncation be less, in L$^2$ norm, than a prescribed accuracy $\varepsilon$.
Note that taking larger values of $S$ allows us to advance forward with larger
time steps. The choice of the optimal truncation order
is discussed in Section~\ref{ss:optimal}.

\subsection{Reversion to an Eulerian grid by interpolation}\label{ss:interpolation}

\begin{wrapfigure}[16]{L}{.33\textwidth}
\vspace*{-1.6\baselineskip}
\begin{center}
\includegraphics[width=.33\textwidth]{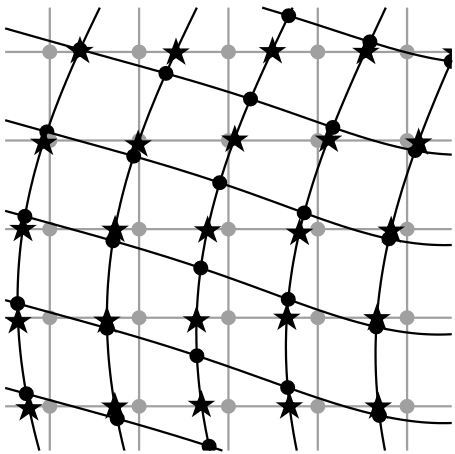}
\end{center}
\vspace*{-1.3\baselineskip}
\caption{Distorted grid (black), image of the uniform
grid (grey) by the Lagrangian map, and hybrid grid (stars).}
\label{f:grid}\end{wrapfigure}

In our implementation of the Cauchy--Lagrangian method, the calculation
of all the time-Taylor coefficients is done by pseudospectral techniques.
For this to be done efficiently, the data at the beginning of a new time step
must be known on a regular grid of $N^d$ collocation points, where $N$ is
the number of grid points per spatial period
and $d$ the dimension of space. This will be called here
the uniform grid. Note that after the Lagrangian
time-stepping is done, the new field, say, the vorticity, can either be
determined in Lagrangian coordinates on the uniform grid or
in Eulerian coordinates on a distorted grid, which is the
image of the original grid by the Lagrangian map during the time interval
$\Delta t$, shown in black in Fig.~\ref{f:grid}. Henceforth, we write just
``image'' for ``image by the Lagrangian map''.

To repeat the time-stepping process, we need to know the
new Eulerian field, say, the vorticity, on a uniform Eulerian grid.
We thus need to interpolate the Eulerian field
from the distorted grid to the uniform grid. As we will see, the distortion
of the grid remains moderate, so that
the interpolation can be done using the {\it cascade interpolation}
introduced in \cite{pl}, which, in dimension $d$, is a superposition
of $d$ one-dimensional interpolations. We now briefly describe it.

We begin with the 2D interpolation procedure, assuming $2\pi$-periodicity
and a time step from $t=0$ to $t=\Delta t$. Positions of fluid particles
at times 0 and $\Delta t$ are $\ba=(a,b)$ and $\bx=(x,y)$\break in the Lagrangian
and Eulerian coordinates, respectively. The Lagrangian map takes $\ba$ to
$\bx(\ba)=(x(a,b),y(a,b))$. The uniform grid is composed of the points
$\ba_{ij}=(a_i,b_j)$. Here, $a_i=2\pi i/N$ and $b_j=2\pi j/N$, where
the indices $i$ and $j$ run from 1 to
the number of grid points per spatial period $N$.
The same uniform grid is used for the Eulerian description after one time step.
The $N$ lines $(a,b_j)$, where $a$ varies continuously in the
interval $[0,2\pi]$, are called the horizontal grid
lines. Similarly, the $N$ lines $(a_i,b)$ are called the vertical grid lines.

The interpolation employs an intermediate hybrid
grid, formed of intersections of the horizontal grid lines with the images
of the vertical grid lines. This hybrid grid is shown
by stars in Fig.~\ref{f:grid}. For the cascade interpolation, it is
essential that each line of the former and latter type have a
single intersection. This follows from a monotonicity result, discussed below.

First, we determine the $x$-coordinates of
hybrid grid points. For each given pair of a horizontal and a vertical line this is
done as follows. The images of the discrete grid points on the vertical line
have coordinates $(x_n,y_n)$. Along the image of the vertical line, in the
neighbourhood of the intersection, a 1D interpolation of $x$ as a
function of $y$ is done; this gives the $x$ value where $y$ takes the value
$y_j$ corresponding to the given horizontal line. Second, the vorticity, which is
known at the same discrete locations $(x_n,y_n)$ along the image of the given
vertical line, is also 1D-interpolated in $y$ to the same values $y_j$. Finally,
along each horizontal line, the vorticity is 1D-interpolated from the points
on the hybrid grid to the points on the uniform grid.
It is straightforward to extend the procedure to 3D with grid lines
and grid planes and again only 1D interpolations (see \cite{pl}).

Now we turn to the issue of monotonicity of the map $b\mapsto
y(a_i,b)$. First a qualitative observation: given the time-analyticity of the
Lagrangian map and the fact that it is slightly smoother in space than C$^1$,
for small $\Delta t$ the map is close to the identity map and its Jacobian
matrix is close to the identity matrix.
Hence, $\partial y(a,b)/\partial b>0$ for sufficiently small $\Delta t$,
which implies the required monotonicity.
This observation can be made quantitative by using
the tools of \cite{zf} (Section 2.3). It is shown there (in 3D) that
analyticity holds if $\Gamma\Delta t<T_c=(8-5\sqrt2)/3$, where
$\Gamma$ is the sum of the moduli of the Fourier coefficients
of the initial vorticity.
It is easy to show that the required positivity then holds at least to the same time.
The positivity may cease to hold if $\Delta t$ is just restricted to
satisfy $\Delta t<\Rei$, but
in practice, we use values of $\Delta t$ significantly smaller
(see Section~\ref{ss:optimal}); in none of the simulations reported in
Section~\ref{s:testing} have we detected any lack of monotonicity.

In the current version of the code, the 1D interpolations are
polynomial and use eight points.
The interpolation procedure can be modified in a number of ways.
We can interpolate first in the $x$-variable and then in the
$y$-variable, and compare the results to estimate the accuracy of
the interpolation. Instead of polynomial interpolation, we can use
rational functions, splines, or any other 1D interpolation algorithm.
We can use different numbers of interpolation nodes in different
regions of the periodicity square, depending, e.g., on how large is
the vorticity gradient (note that in Fig.~\ref{figisov}
of Section~\ref{s:testing} the areas of large gradients are small for $t\ge3$).
For the interpolation we can use a mesh finer, than the one
used to compute the Lagrangian map. Spatial derivatives of the vorticity
can also be used for interpolation.
We will experiment with these modifications in future studies.

Errors of interpolations can be estimated on the fly as follows:
the interpolation yields the vorticity on the target regular grid;
from this, by FFT we can calculate its Fourier coefficients. By
direct summation of the Fourier series (``slow Fourier transform''), we
can then compute the vorticity on the distorted grid and find how closely it
matches the one we employed for the interpolation. Such
computations are CPU-intensive, but it suffices to perform them at those points
where the largest interpolation errors are suspected, e.g.,
where high-order spatial derivatives of the vorticity are large.

\subsection{Choosing the optimal truncation order and time step}\label{ss:optimal}

In computations, the displacement is given by the time-Taylor series
\rf{xiSeriesS} truncated at order $S$. For a given accuracy
$\varepsilon$, how should we choose $S$ and
the time step $\Delta t$ to minimise the complexity of computations
when integrating in a time interval $[0,T]$? We assume that the
spatial resolution and thus the complexity of the spatial FFTs are fixed.
The order of the interpolations needed for
reversion to Eulerian coordinates is held fixed as well.
Furthermore, we assume that the radius of convergence $\Rei(t)$
of the time-Taylor series around time $t$ does not vary significantly
in the interval $[0,T]$ and we denote it $R$.
As we will see in Section~\ref{s:depletion}, at least in two
dimensions, $\Rei(t)$ does not decrease in time dramatically, and hence
this is a reasonable approximation.

We prescribe the accuracy $\varepsilon$ and apply the condition \rf{estet0}.
The L$^2$-norm series $\sum_{s=0}^\infty\|\bxi^{(s)}\|(\Delta t)^s$,
which has positive coefficients and a radius of convergence $R$, can be
approximated, at least for large orders, by the geometric series
$A\sum_{s=0}^\infty(\Delta t/R)^s$, where $A$ is a
positive constant. This exponential dependence on the order,
\BE\|\bxi^{(s)}\|\approx AR^{-s},\EE{expoorder}
is well supported by the numerical results on high-order
time-Taylor expansions of Section~\ref{s:rounding} (see,
in particular, Fig.~\ref{coeft0135}).

For large truncation orders $S$, given the nonlinearity of order $d$ stemming
from the determinant in \rf{dgra}, the number of operations is O$(S^d)$. Here
$d$ is the space dimension.
The number of time steps for integration to time $T$ is $T/\Delta t$;
\rf{estet0} and \rf{expoorder} yield $\Delta t=R(\varepsilon/A)^{1/S}$.
Hence, the number of operations needed to reach a given time
is proportional to $S^d(\varepsilon/A)^{-1/S}$. Minimising this expression
over $S$, we find the optimal order
$S\approx-(1/d)\ln(\varepsilon/A)$ and $\Delta t=R\,\e^{-d}$.
For small $S$ and large $N$, the complexity $SN^d\ln N$ of FFTs needed for
applying the Calderon--Zygmund operators involved in the expression of
$\bxi^{(s)}$ becomes dominant. By similar arguments, the optimal order is then
$S\approx-\ln(\varepsilon/A)$. When both complexities
are essential, the optimal order of expansion, $S$, satisfies the inequality
\BE-(1/d)\ln(\varepsilon/A)\le S\le-\ln(\varepsilon/A).\EE{opS}

\section{Testing the Cauchy--Lagrangian numerical method in two dimensions:\\
Comparison with Eulerian simulations}
\label{s:testing}

In Section~\ref{ss:flow-methods} we specify the test flows and
describe three Eulerian algorithms that have been
used for comparison with the Cauchy--Lagrangian algorithm.
Validation, with emphasis on accuracy, is presented in
Section~\ref{ss:validation}. Efficiency of the CL algorithm is
discussed in Section~\ref{ss:efficiency}. In Section~\ref{ss:tygers}
we discuss spatial truncation artefacts. All computations reported in this
section are in double precision.

\subsection{Flows and numerical methods used for validation}
\label{ss:flow-methods}

For the reason explained in Section~\ref{ss:init}, all our tests of the CL
method have been done using 2D flows with analytic initial data having
non-trivial dynamics in both Eulerian and Lagrangian coordinates.
Two different flows, called the \textit{test flows},
were used as initial condition. The first one is a very simple
flow, here called the ``4-mode'' flow, with the initial vorticity
\BE\omega^{\rm(init)}=\cos x+\cos y+0.6\cos 2x+0.2\cos 3x.\EE{4mode}
The second one is a particular realisation of the so-called ``random'' flow
\cite{takesh}, used in Section~II.C of \cite{tyger},
where the time evolution of the latter flow is presented in detail.
The characterisation of the random flow is best done in the Fourier space,
where each wave vector consists of a pair of signed integers
$\bk\equiv(k_1,\,k_2)$, conventionally combined into shells
$K\le|\bk|<K+1$ where $K$ is integer. Each such shell involves $N(K)$ Fourier
harmonics. For all $\bk$ in the $K$th shell, the Fourier coefficients
$\widehat\omega_{\bk}$ of the initial vorticity are assigned the same
modulus $2K^{7/2}\exp(-K^2/4)/N(K)$ and phases that are uniformly and
independently distributed in the interval $[0,2\pi]$, except that
the coefficients for opposite wave vectors are complex conjugate,
so that the flow is real.

\begin{figure}[!t]
\begin{picture}(132,36)(0,0)
\put(10,0){\psfig{file=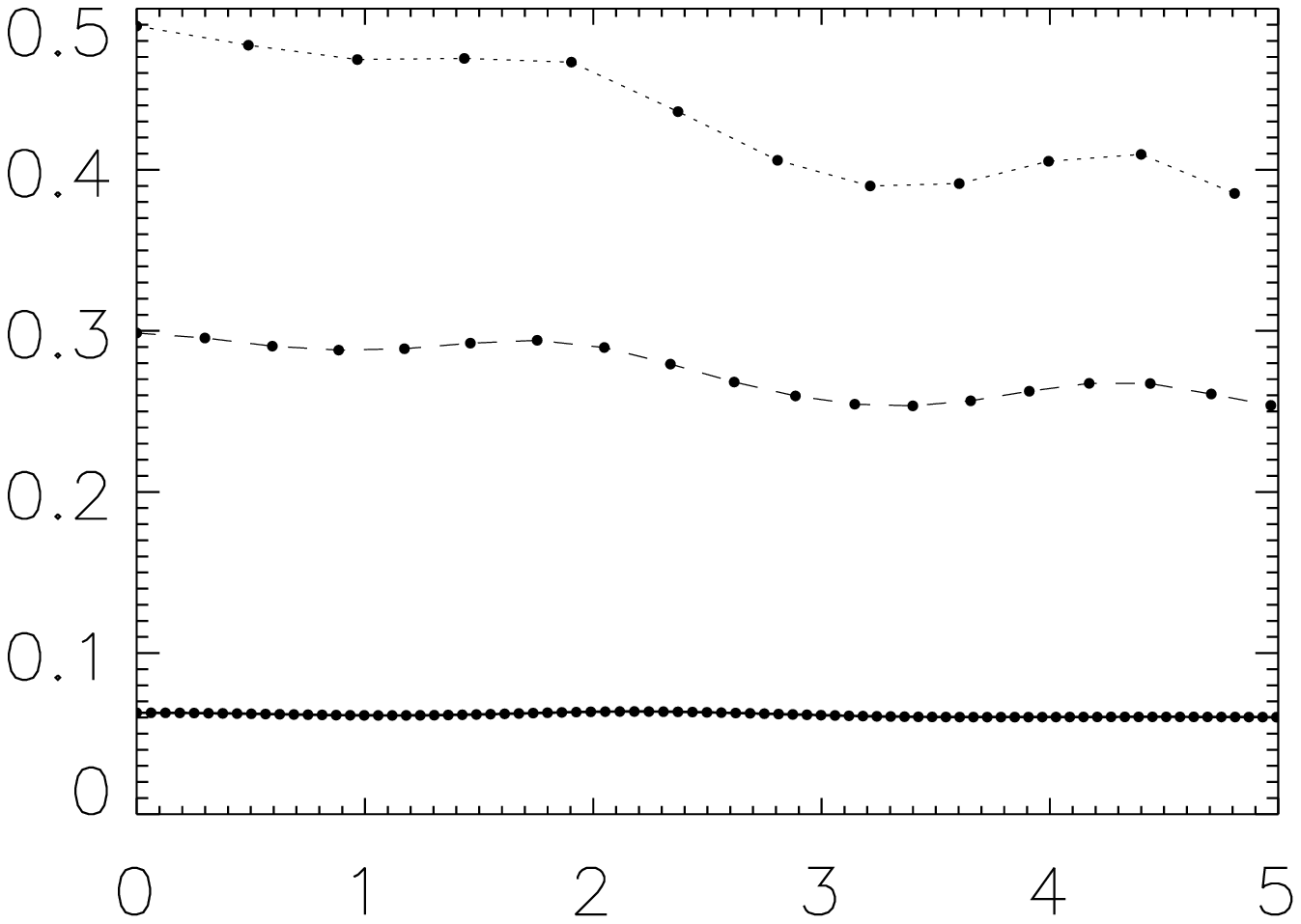,width=67mm}}
\put(4,30){(a)}
\put(76,0){\psfig{file=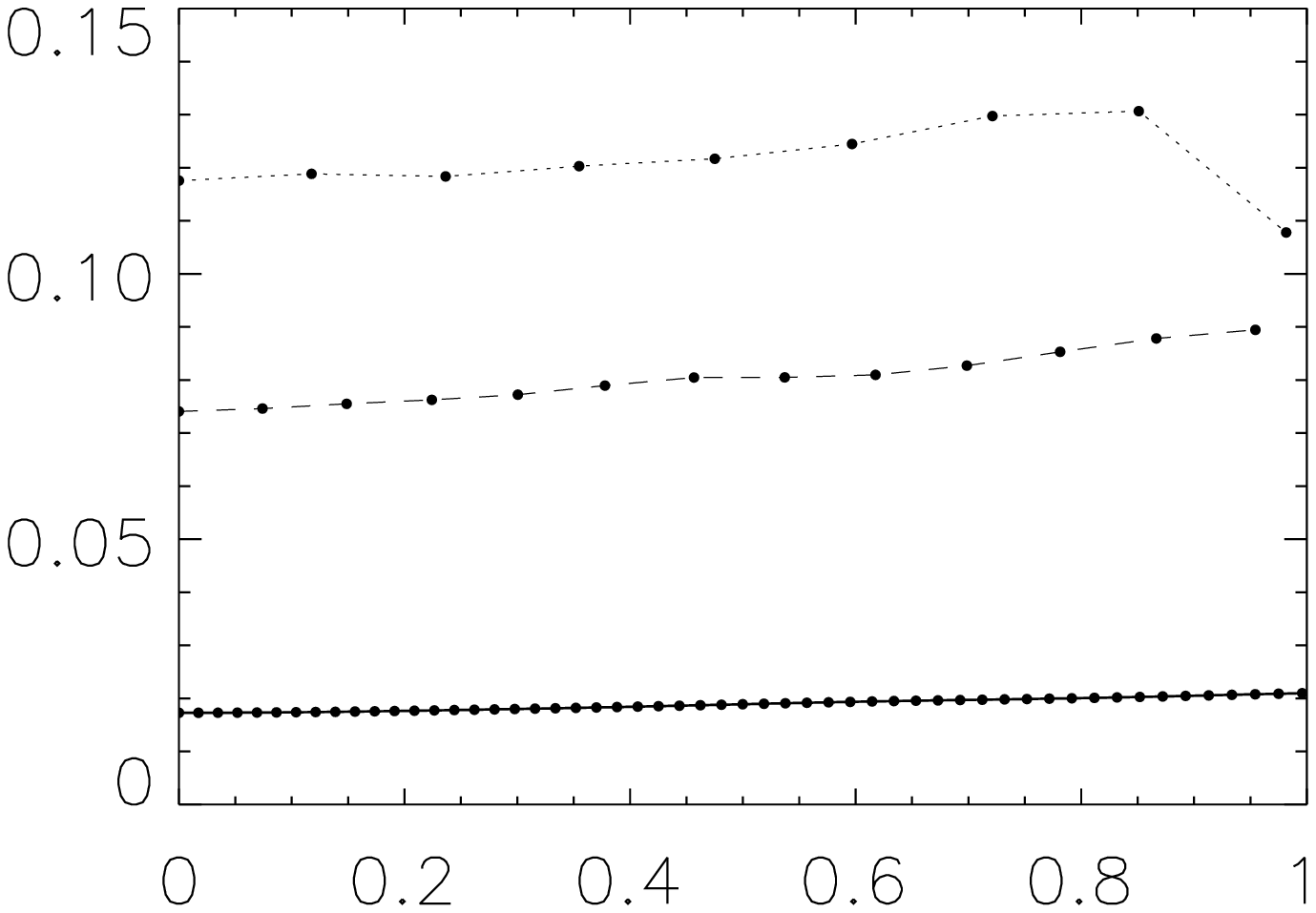,width=67mm}}
\put(70,30){(b)}
\put(17.5,24){$\Delta(t)$}
\put(83.5,23){$\Delta(t)$}
\put(42,29){CL24}
\put(42,20.5){CL16}
\put(42,8.5){CL8}
\put(106,30.5){CL24}
\put(106,21.5){CL16}
\put(106,9){CL8}
\put(54,5){$t$}
\put(120,5){$t$}
\end{picture}
\caption{Temporal variation of the time step of the CL code for the 4-mode (a)
and random initial conditions (b).}\label{steps}\end{figure}

We have used three programs to solve the 2D Euler equation, realising
the CL algorithm with the Taylor series truncated to order $S$ (denoted CL$S$),
described in Section~\ref{s:CLalgorithm}, and three algorithms working
exclusively in Eulerian variables, namely, the Runge--Kutta algorithms of
orders two and four (denoted RK2 and RK4, respectively) and the algorithm
relying on the Eulerian time-Taylor expansion truncated to order $S$ (denoted
ET$S$). The CL algorithms used in this section are CL8, CL16 and CL24 with time steps,
chosen as explained in Section~\ref{ss:optimal}, and with the accuracy
in \rf{estet0} set to $\varepsilon=10^{-12}$. Figure~\ref{steps} shows the
evolution of the time steps for the test flows. In principle, they are
allowed to vary in time, but, as seen, they change very moderately. This is
apparently related to the very slow temporal variations of the radii of
convergence of the time-Taylor series, discussed in Section~\ref{s:depletion}.

Our calculations have mostly been performed with spatial resolutions of
$1024^2$ and $2048^2$ Fourier harmonics before dealiasing (performed using the
standard 2/3-rule applying to quadratic problems). One calculation
had a resolution of $8192^2$ harmonics.

Concerning the standard RK2 and RK4 algorithms,
we only comment on the values of the time steps, taken constant
in our calculations. For RK2, the smallest time step that we were able
to run in reasonable CPU time was $\Delta t=10^{-4}$. As we will see,
this has yielded accuracies in
the range $10^{-6}-10^{-9}$, significantly worse than the ones obtained
with the other algorithms, requiring comparable or smaller CPU time.
For RK4, the time steps used are listed in Table~\ref{tb1}. They have been
chosen as the largest value that gives stable integration,
so that doubling of the time step causes numerical
overflow. The corresponding Courant numbers (based on the maximum
velocity and the maximum wave number after dealiasing) are
indeed slightly in excess of unity.

\begin{table}[b!]
\caption{Time step $\Delta t$ for RK4 and ET8 algorithms.}
\center\begin{tabular}{|c|cc|}\hline
Resolution\vphantom{$|^|_|$}&$1024^2$&$2048^2$\\\hline
4-mode initial flow\vphantom{$|^|_|$}&0.0025&0.00125\\
random initial flow&0.002 &0.001\\\hline
\end{tabular}\label{tb1}\end{table}

The Eulerian Taylor (ET) algorithm deserves more explanations.
By analogy with time integration exploiting analyticity of Lagrangian
fluid particle trajectories,
one can perform Eulerian time integration exploiting the analyticity
in time of the velocity and the vorticity when the initial data are
analytic in space and remain so for at least a finite time. The ET algorithm
was actually applied long ago to study real
and complex time singularities of the Taylor--Green
3D vortex flow \cite{br}. This was done with a single time step from 0
to $t$, followed by attempted analytic continuation in $t$ beyond the disk of
convergence using, e.g., Pad\'e approximants.
Unfortunately, this method failed to give clear evidence for or
against finite-time blow-up in 3D.
Here, we are interested in the 2D multistep form of the ET algorithm
applied to the vorticity formulation of the Euler equation
\BE{\partial\omega\over\partial t}+(\bv\cdot\nabla)\,\omega=0,\qquad
\nabla\times\bv=\omega\,{\bf e}_z,\EE{eulvort}
where $\omega$ denotes the only component of the vorticity and ${\bf e}_z$ is the unit vector in the z-direction. We assume that
the vorticity field $\omega$ is known at time $t$, and Taylor-expand
$\omega(t+\Delta t)$ in powers of $\Delta t$:
\BE\omega(t+\Delta t)=\sum_{s=0}^{\infty}\omega_s(\Delta t)^s.\EE{eults}
After substitution in \rf{eulvort}, we obtain the following recurrence
relations (where $\omega_0\equiv\omega(t)$):
\BE(s+1)\,\omega_{s+1}+\sum_{m=0}^s(\bv_m\cdot\nabla)\,\omega_{s-m}=0,\qquad
\nabla\times\bv_m=\omega_m\,{\bf e}_z,\quad s=0,1,\ldots.\EE{ssterm}
These recurrence relations look significantly simpler than the Lagrangian
recursion relation \rf{bigXis}. Actually, the Lagrangian one is much
better conditioned than the Eulerian ones, as we shall see in Section~\ref{s:rounding}.

Time-stepping is done, as in the CL method, by applying the Taylor series
truncated to a finite order $S$ and choosing an appropriate
time step $\Delta t$. Here, we report results obtained with
an 8-term truncation, called ET8. Higher-order truncations have special
problems, that will be discussed in Section~\ref{s:rounding}.

Given the analyticity in time of the Eulerian solution (for all times in 2D
and a finite time~in 3D \cite{bb}),
we certainly must require that the time step be smaller than the (Eulerian)
radius of convergence of the Taylor series. However, no Eulerian analogue
of the condition \rf{estet0} is available,
and much smaller time steps must be chosen; typically Courant numbers of order
unity are needed. The reason, however, is not stability, as in low-order
schemes, but the rounding errors, as explained in Section~\ref{s:rounding}.
In practice, for ET, we use the time steps of the same length
as for RK4, given in Table~\ref{tb1}.

Finally, a few words about spatial truncation errors. Since we restrict
ourselves to analytic initial data, the solutions are also analytic in both
Lagrangian and Eulerian coordinates. The Fourier transforms of such
functions decrease exponentially with the wave number $k$ as
$\e^{-\delta k}$, where $\delta$ is the radius of the cylinder of analyticity
(also sometimes called ``width of the analyticity strip''). As shown
in \cite{ssl}, $\delta$ can be estimated from the high-$K$ asymptotics
of the vorticity spectrum (also called ``enstrophy spectrum''),
\BE E_\omega(K)\equiv{1\over2}\sum_{K\le|\bk|<K+1}|\widehat\omega_{\bk}|^2,\EE{defomegaspec}
where $\widehat\omega_{\bk}$ are the spatial Fourier coefficients of
the vorticity $\omega$. The large-wave-number behaviour of the vorticity
spectrum is typically
\BE E_\omega(K)\sim K^\alpha\e^{-2\delta K}.\EE{asenspec}
Here, the distance $\delta$ tells us how far away from the real domain
the nearest complex singularities are, while the exponent $\alpha$ contains
information about the type of singularities \cite{pa}. As a consequence
of \rf{asenspec}, the relative spatial truncation errors are typically estimated as
$\e^{-\delta k_{\max}}$, where $k_{\max}$ is the maximum wave number
after dealiasing, $k_{\max}=N/3$, where $N$ is the number of grid points
per spatial period. The numerical determination of $\delta$ from
the vorticity spectrum is done by the same fitting method as for the
numerical determination of radii of convergence (see \ref{a:conver}). When
using the Cauchy--Lagrangian method, at each new time step, the Lagrangian
and Eulerian descriptions coincide and, as a consequence, the Lagrangian and
Eulerian $\delta$ are close.

\subsection{Validation of the CL algorithm and accuracy of agreement}
\label{ss:validation}

\begin{figure}[t]
\begin{picture}(132,36)(0,0)
\put(10,0){\psfig{file=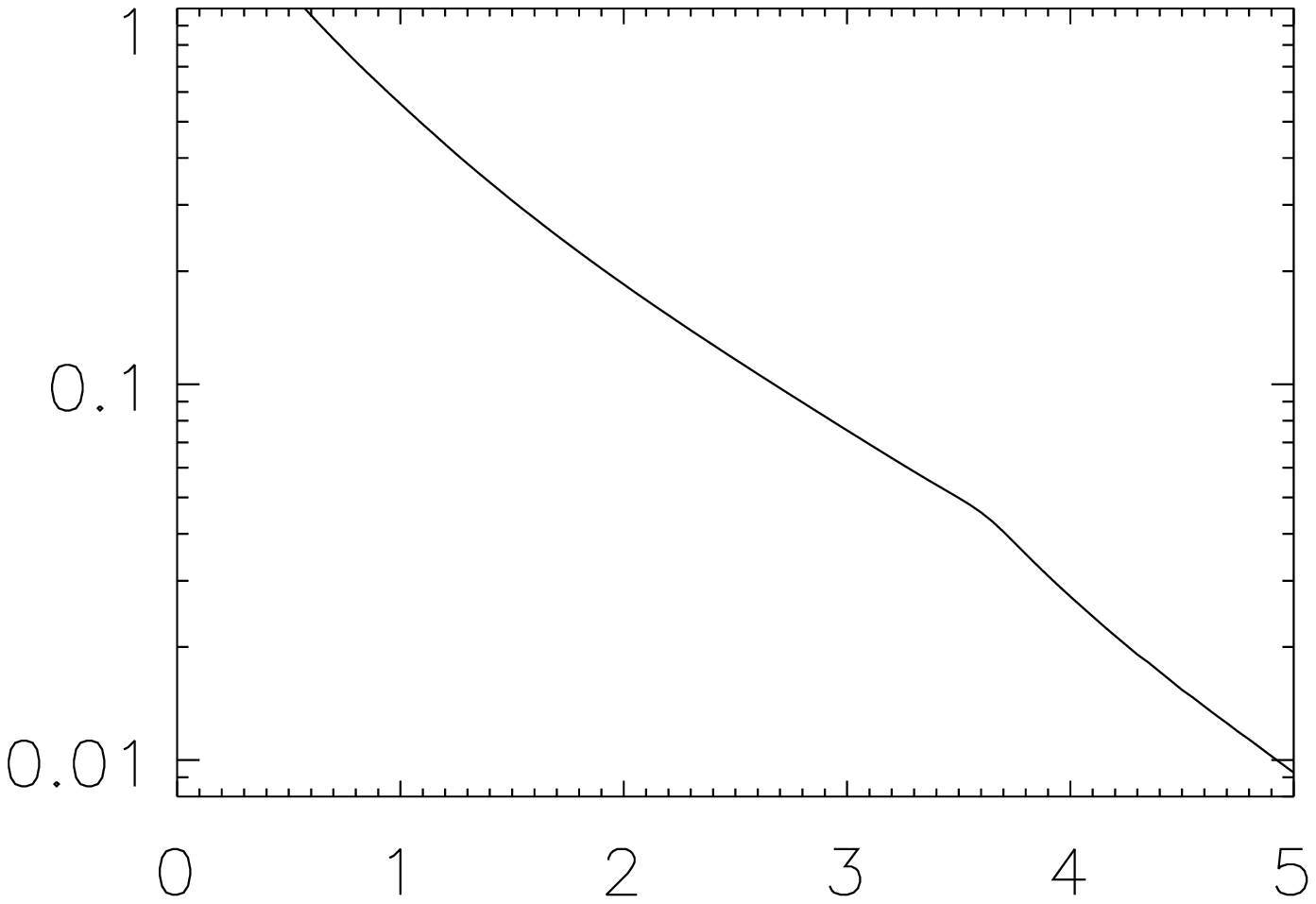,width=67mm}}
\put(4,30){(a)}
\put(76,0){\psfig{file=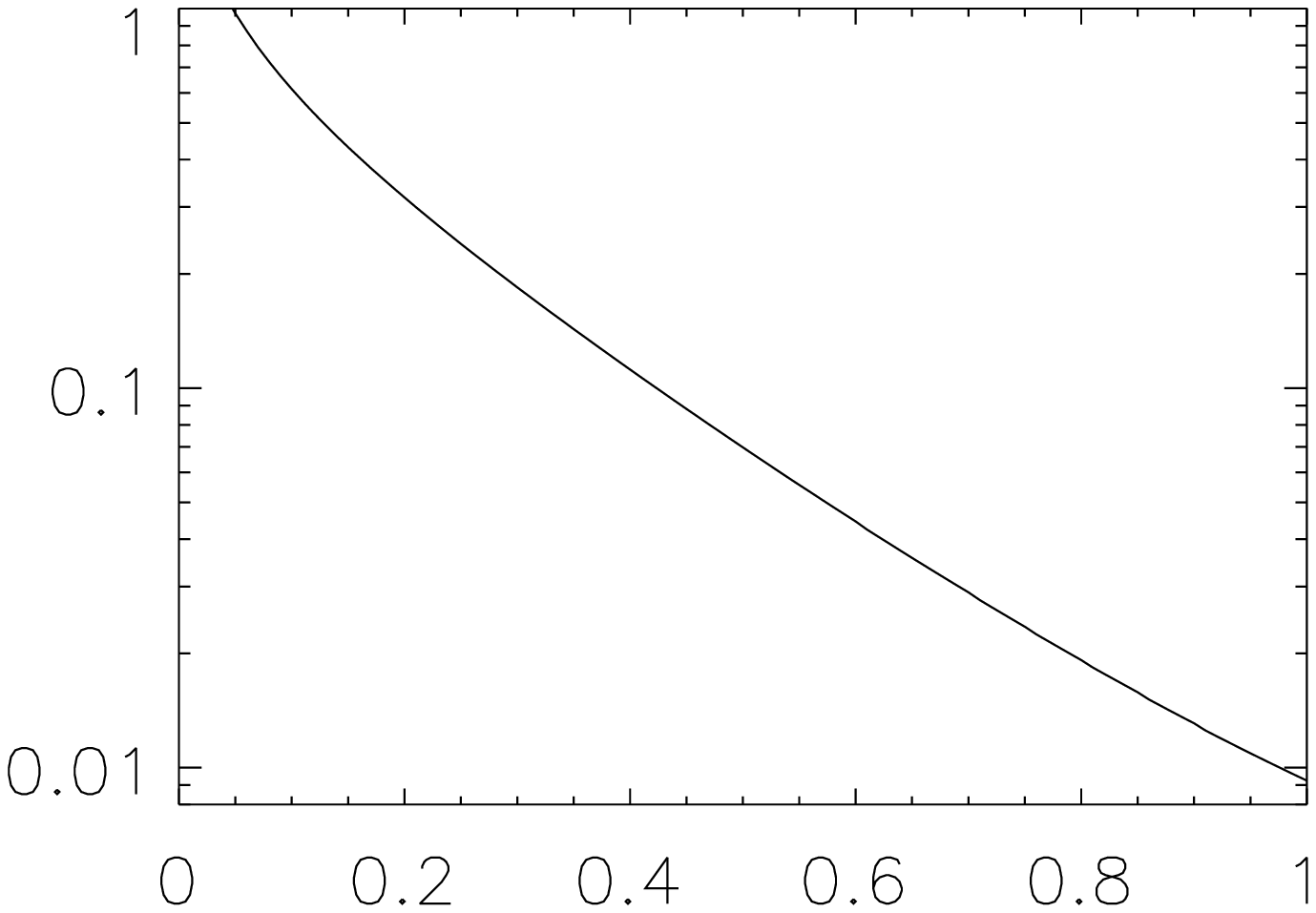,width=67mm}}
\put(70,30){(b)}
\put(18,28){$\delta$}
\put(84,28){$\delta$}
\put(51,5){$t$}
\put(117,5){$t$}
\end{picture}
\caption{Radius $\delta(t)$ of the analyticity cylinder for the 4-mode
(a) and the random (b) initial condition. CL8 algorithm.
Resolution: $8192^2$ harmonics.}\label{fdel}

\begin{picture}(132,36)(0,0)
\put(10,0){\psfig{file=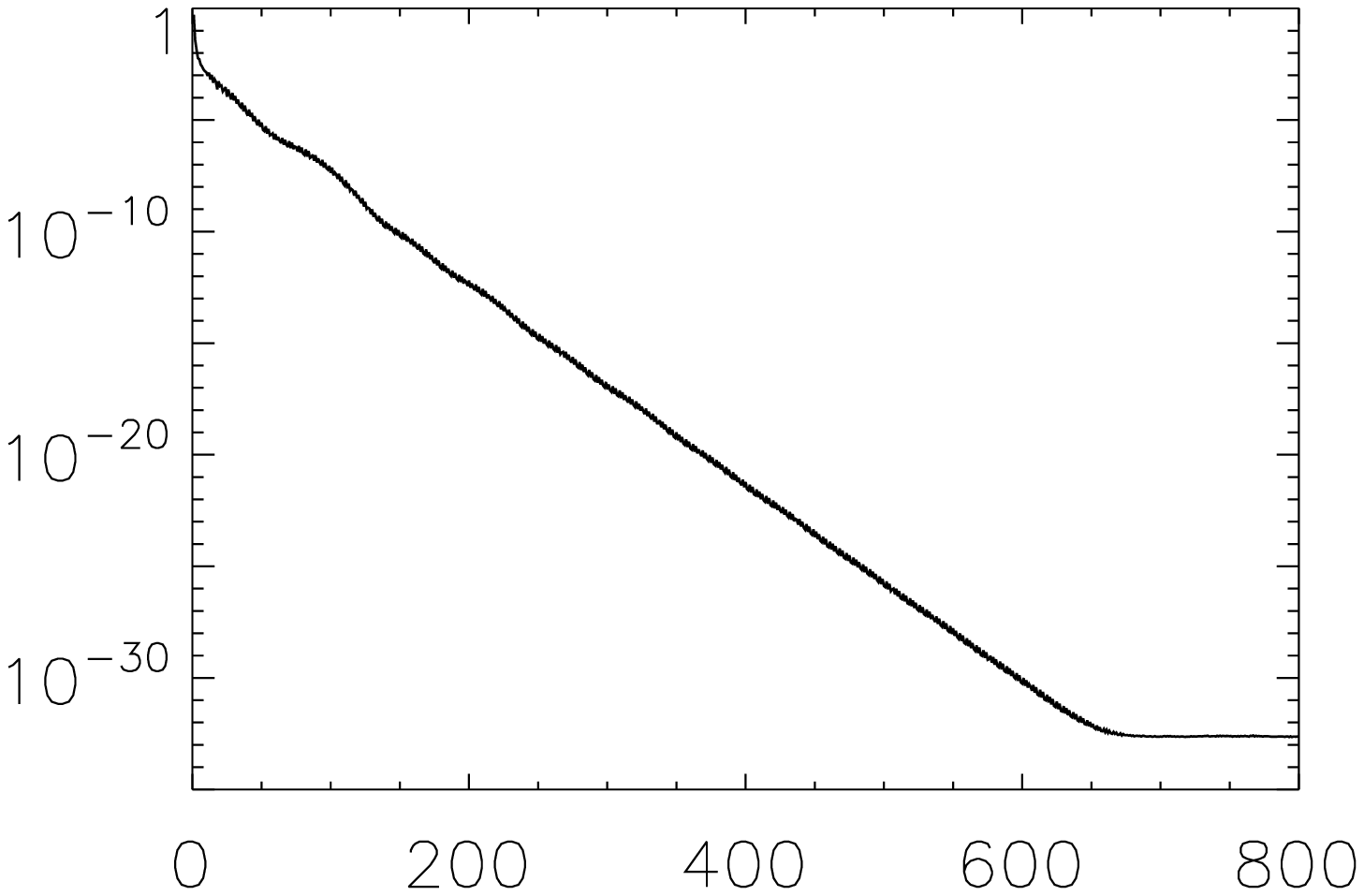,width=67mm}}
\put(4,30){(a)}
\put(76,0){\psfig{file=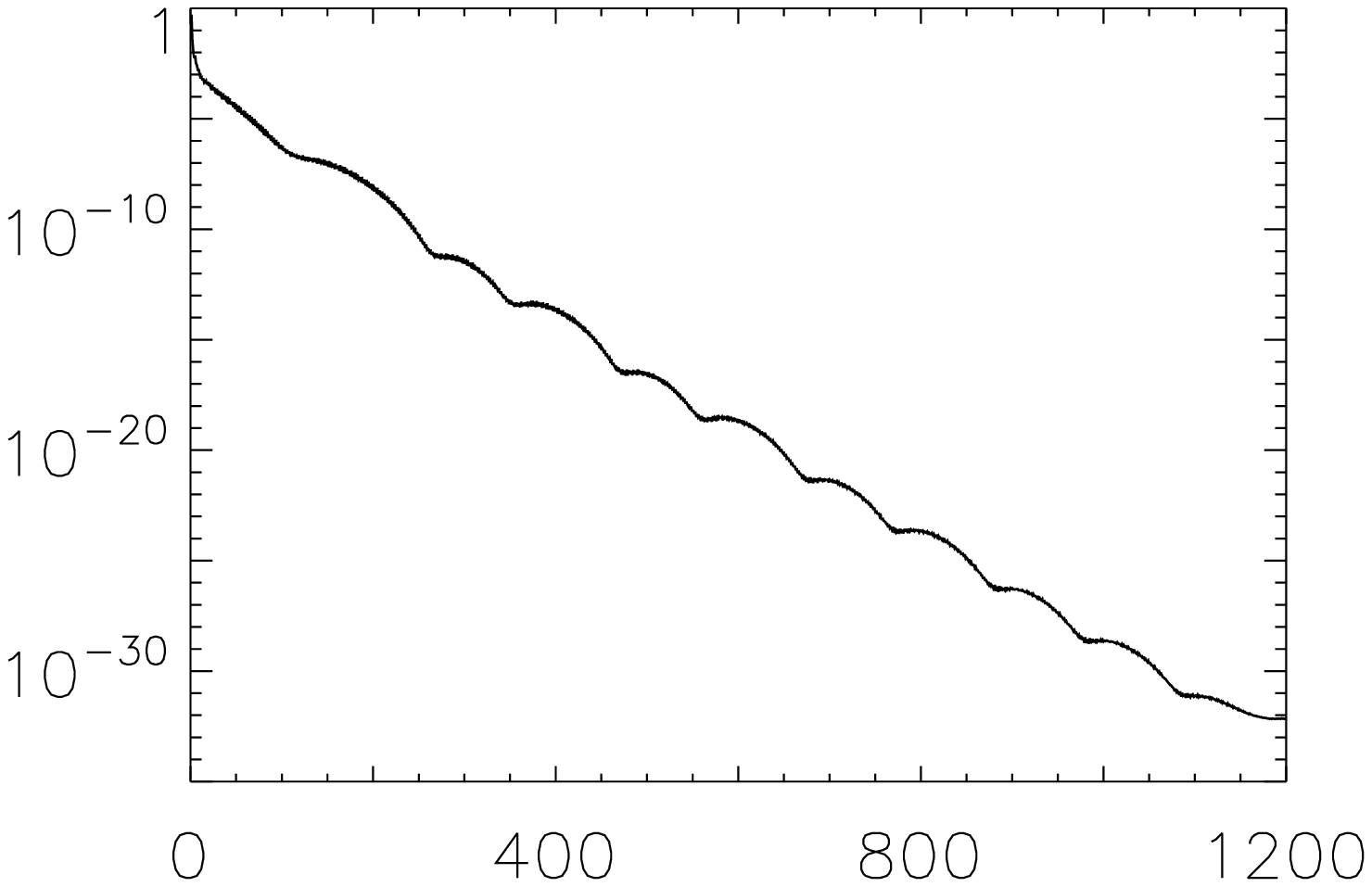,width=67mm}}
\put(70,30){(b)}
\put(18,18){$E_\omega(K)$}
\put(84,18){$E_\omega(K)$}
\put(42,5){$K$}
\put(108,5){$K$}
\put(42,26){$t=3.5$}
\put(108,26){$t=4$}
\end{picture}
\caption{Spectra of the vorticity for the 4-mode initial condition
at several times, as labelled. CL8 algorithm. Resolution: $8192^2$ harmonics.}
\label{specs}\end{figure}

For each test flow (see Section~\ref{ss:flow-methods}), we have to
determine up to what maximum time $\tf$ we can integrate the
2D Euler equation without encountering excessive loss of
accuracy. This will of course depend on the resolution and on how
rapidly the Fourier coefficients fall off with the wave number $k$.
As explained at the end of Section~\ref{ss:flow-methods}, this fall-off is
controlled by the radius of spatial analyticity, measuring the distance
$\delta(t)$ of the nearest complex-space singularities to the real-space
domain. To monitor $\delta(t)$, we have used the Cauchy--Lagrangian method
to integrate the Euler equation with the resolution of $8192^2$ Fourier
harmonics. The measured $\delta(t)$ for the two test flows are shown in
Fig.~\ref{fdel}. Each point on these graphs is obtained by processing
the vorticity spectrum $E_\omega(K)$ at the corresponding time by
the fitting technique of \ref{a:conver}. Examples
of vorticity spectra are shown in Fig.~\ref{specs}.

We determine from Fig.~\ref{fdel} that, at the largest times shown for
both flows, $\delta$ is about $10^{-2}$. Since $k_{\max}=8192/3\approx2731$,
the resulting relative truncation error is about
$\e^{-\delta(t)k_{\max}}\approx 2\times10^{-12}$. (The product
$\delta(t)k_{\max}$ is sufficiently large to use this asymptotic wave number
formula.) We also note that for the lowest resolution used in our simulations,
namely, $1024^2$ harmonics, the 4-mode flow achieves comparably
small spatial truncation errors only up to $t=2$. If, however, we just request
a level of spatial truncation error of about $10^{-4}$, which is not visually
detectable on contours of the vorticity and of its Laplacian, we can
extend the computations to about $t=4$. We will come back to this
in Section~\ref{ss:tygers}, when discussing truncation artefacts.

Before comparing the results for the different methods, it is of
interest to display the evolution of the test flows graphically. For
the random flow, this was done in Section~II.C of \cite{tyger}. Here we focus
on the 4-mode flow.

\begin{figure}[t!]
\centerline{\psfig{file=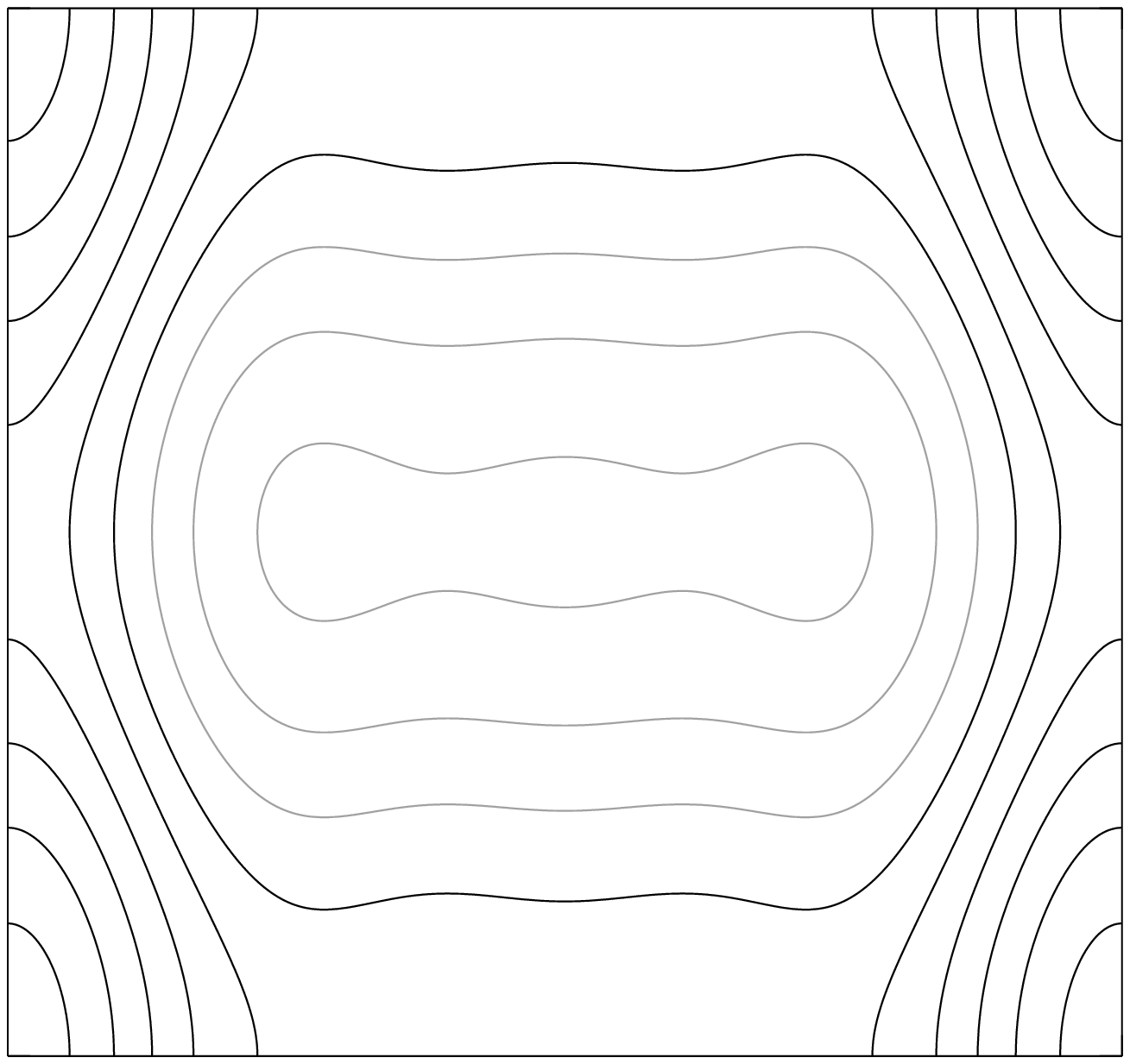,width=45mm}~\psfig{file=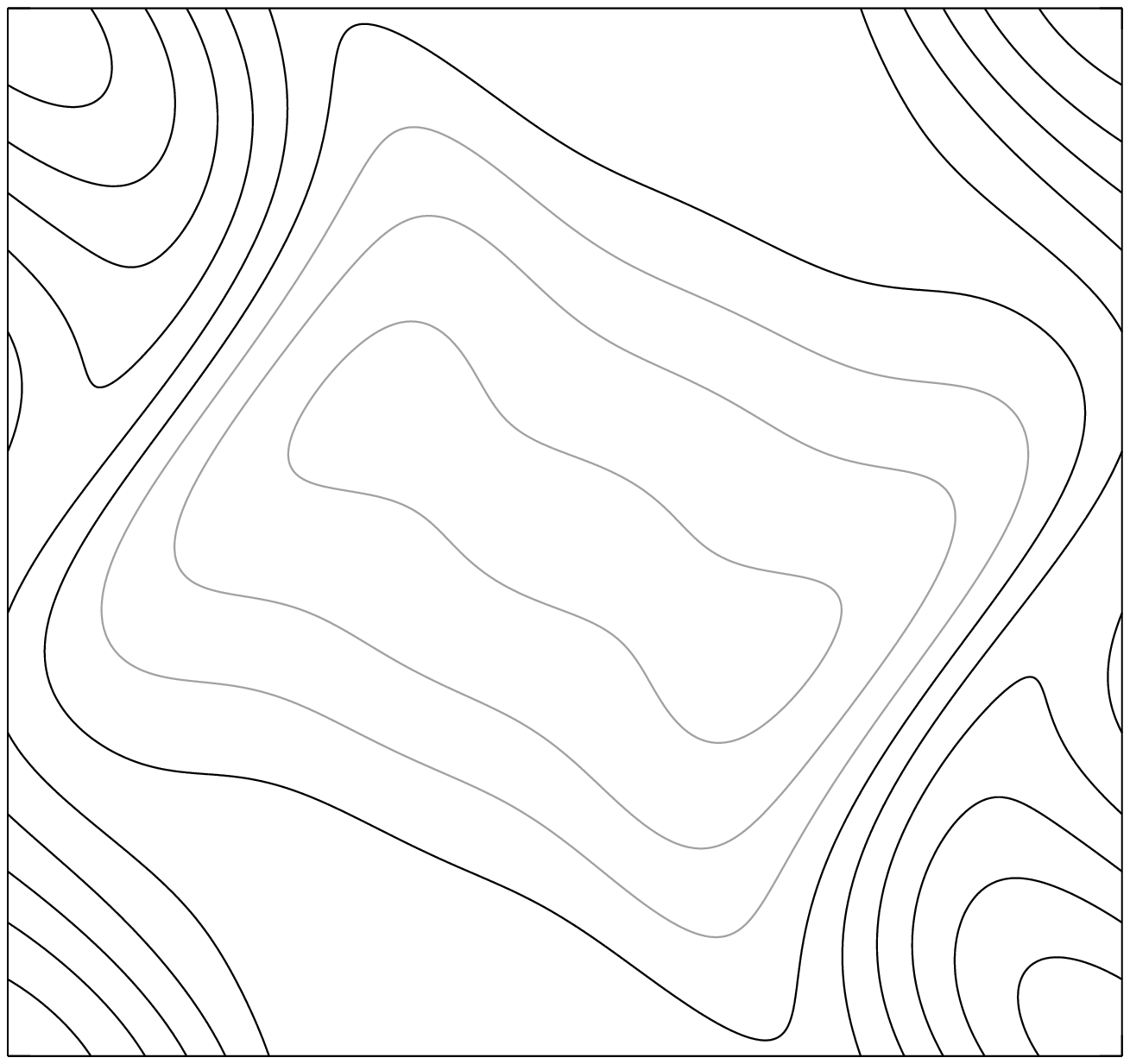,width=45mm}~\psfig{file=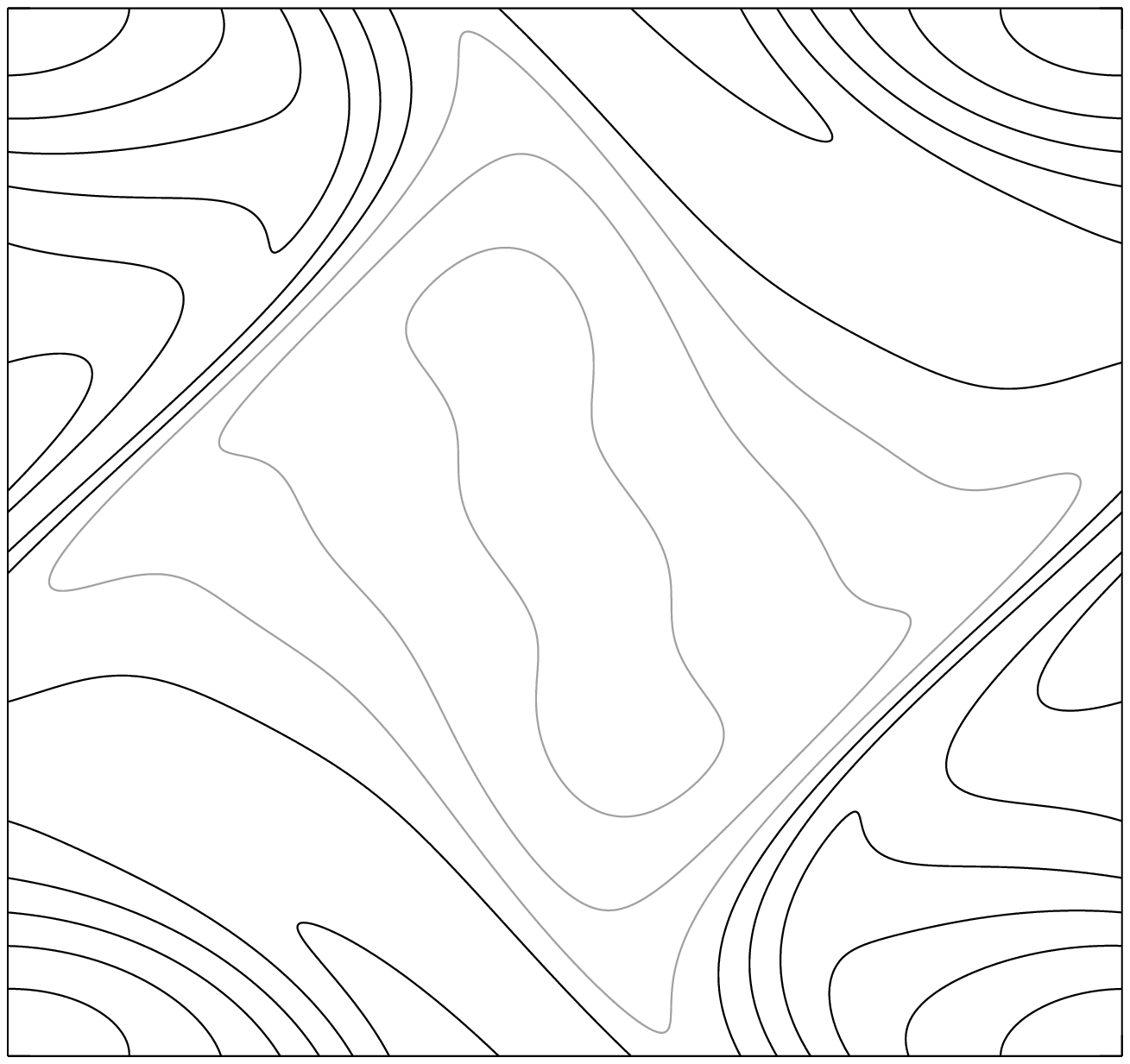,width=45mm}}

\vskip7mm
\centerline{\psfig{file=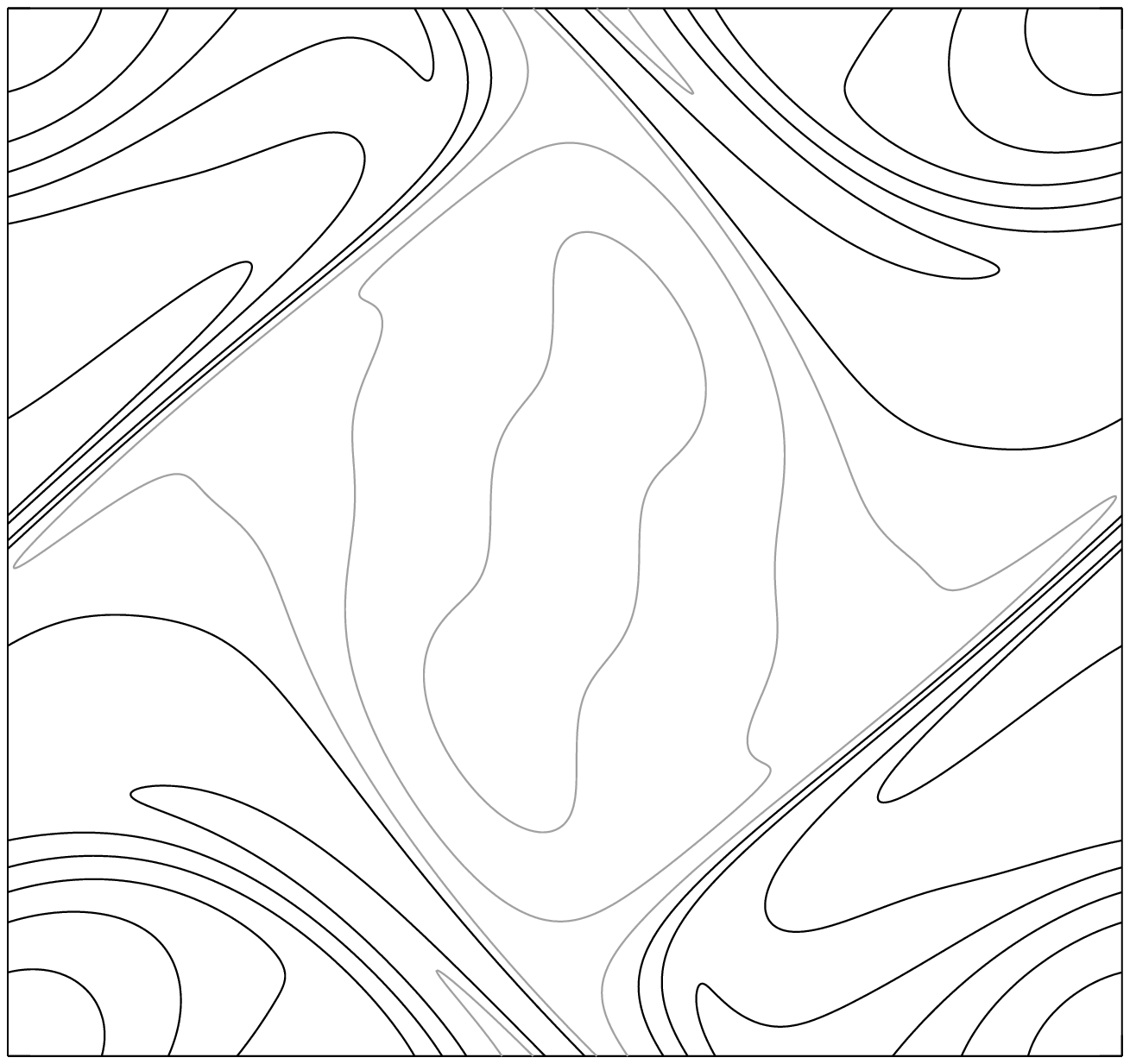,width=45mm}~\psfig{file=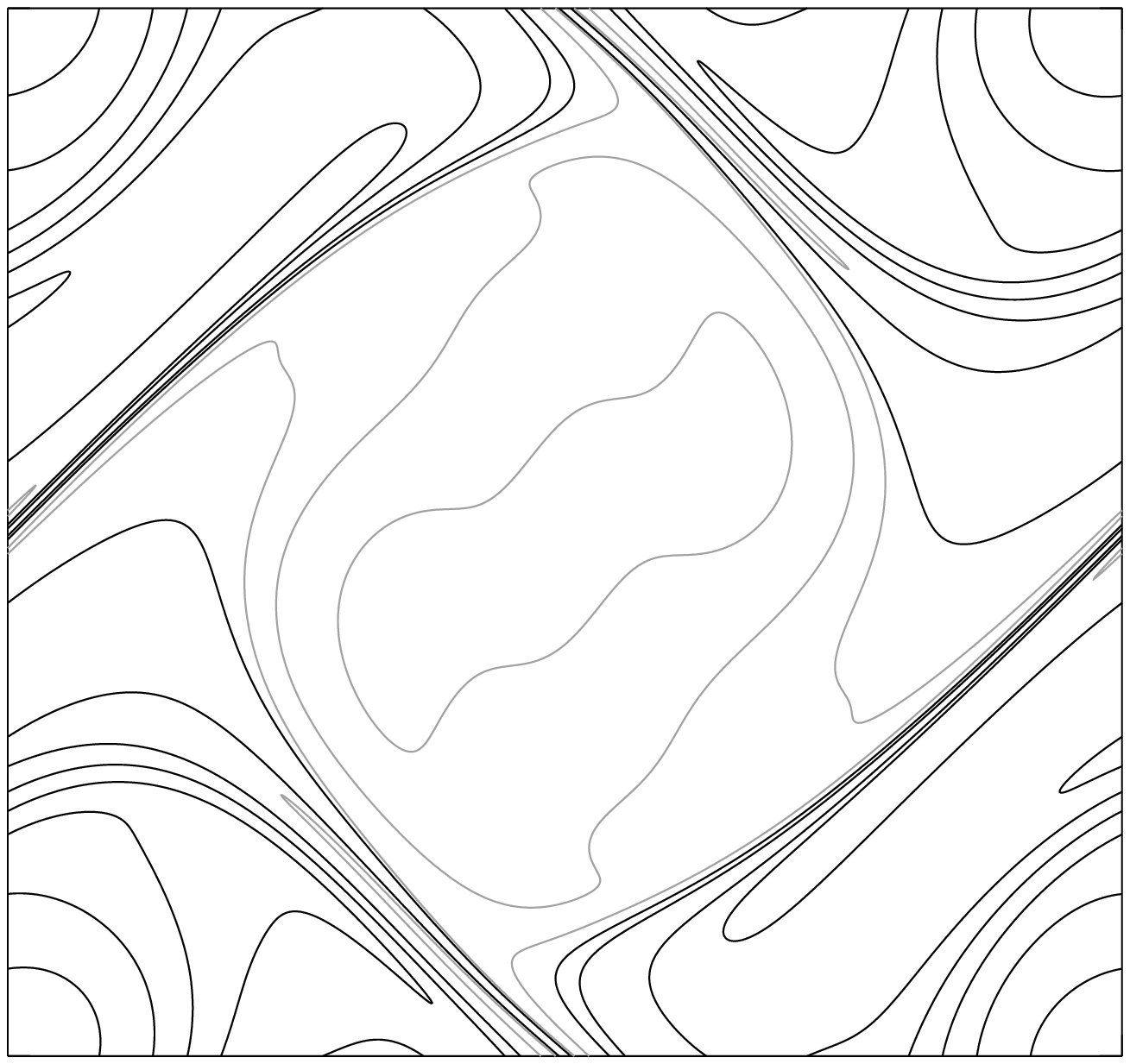,width=45mm}~\psfig{file=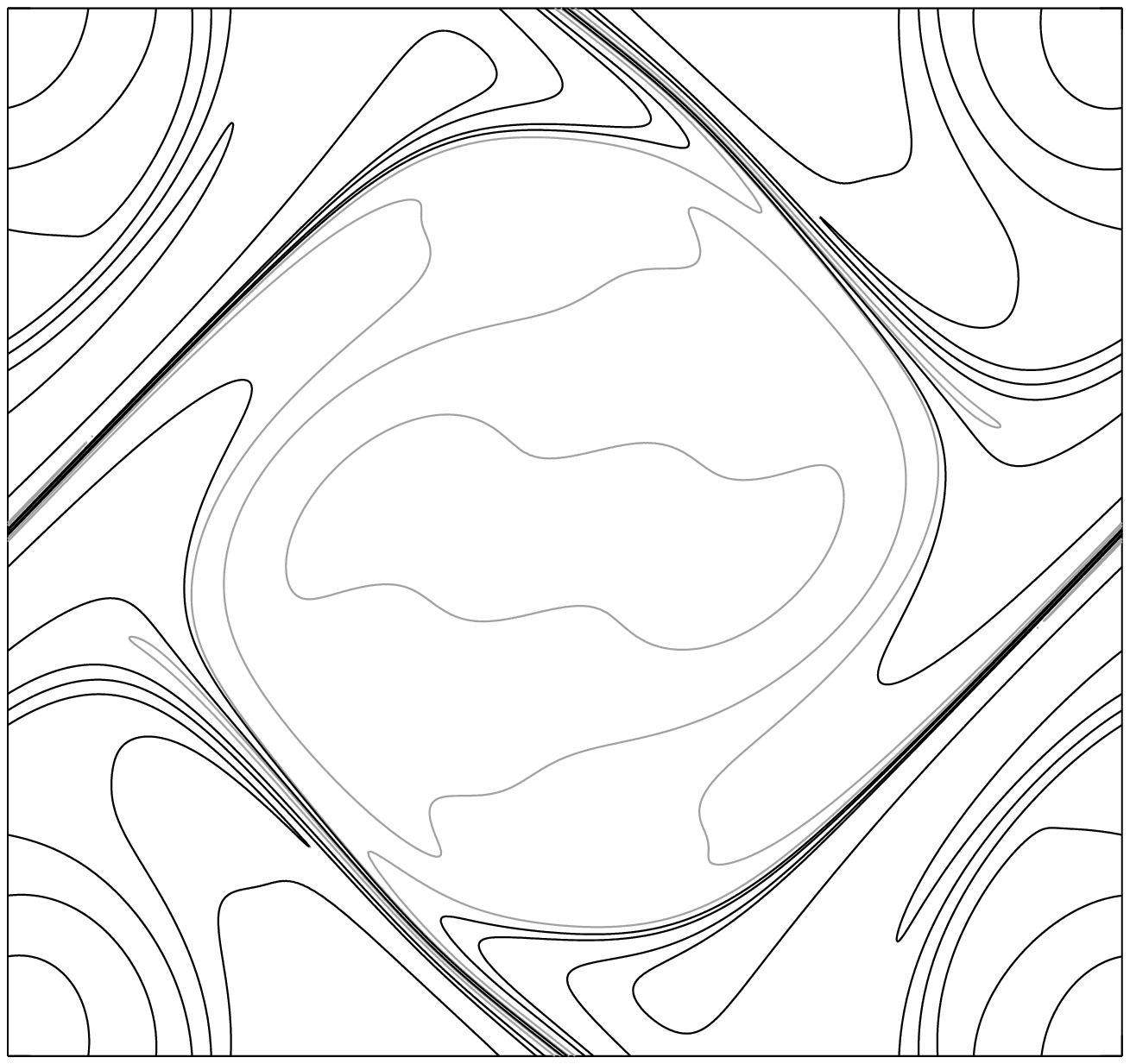,width=45mm}}

\vspace*{-98.5mm}\hskip26mm$t=0$\hskip37mm$t=1$\hskip37mm$t=2$

\vspace*{44.5mm}\hskip26mm$t=3$\hskip37mm$t=4$\hskip37mm$t=5$

\vskip42mm
\caption{Isolines of the vorticity (step 0.5) computed by CL8 with the
resolution of $1024^2$ harmonics for the 4-mode initial condition.}
\label{figisov}

\vspace*{1.5\baselineskip}
\centerline{\psfig{file=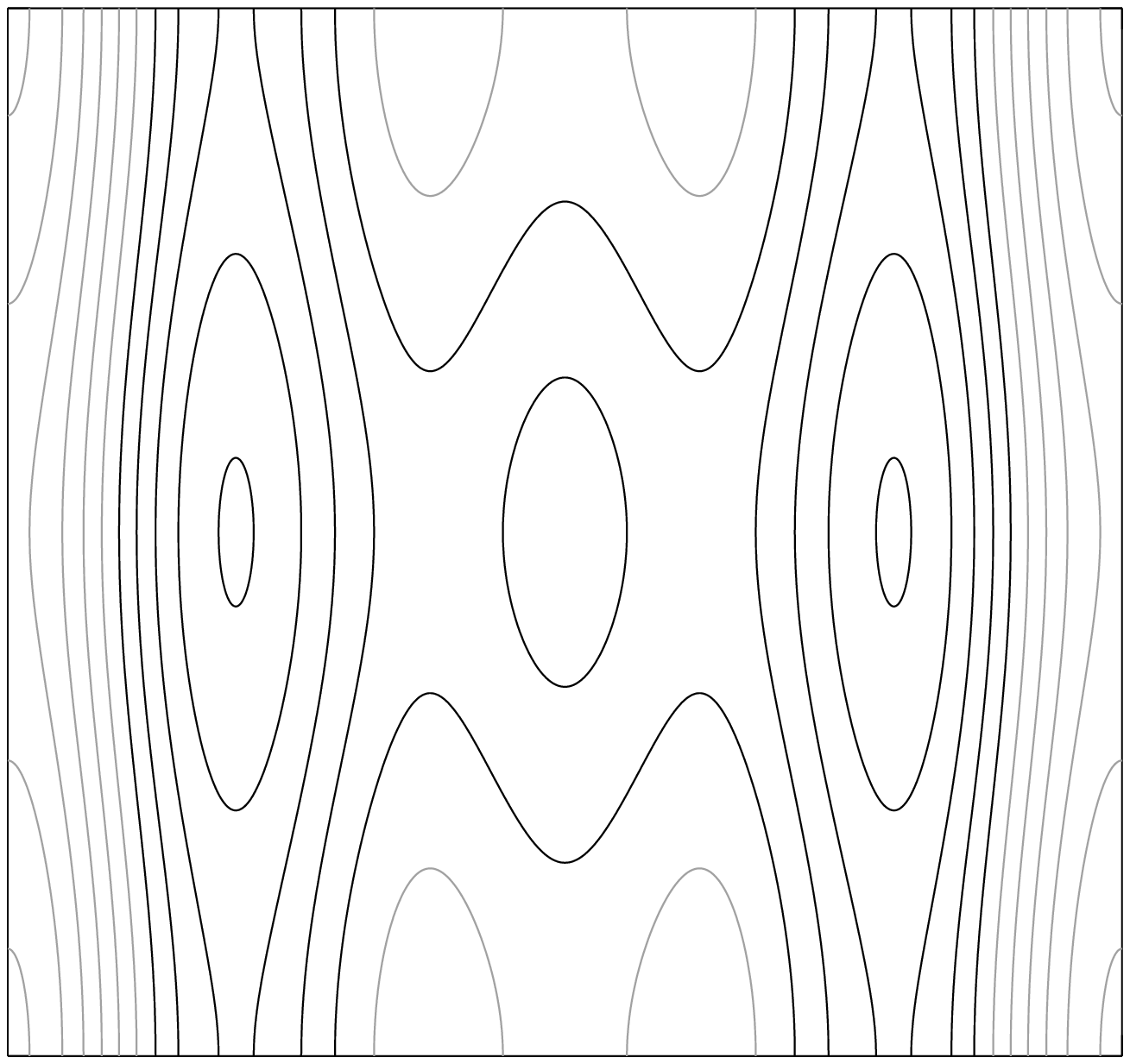,width=45mm}~\psfig{file=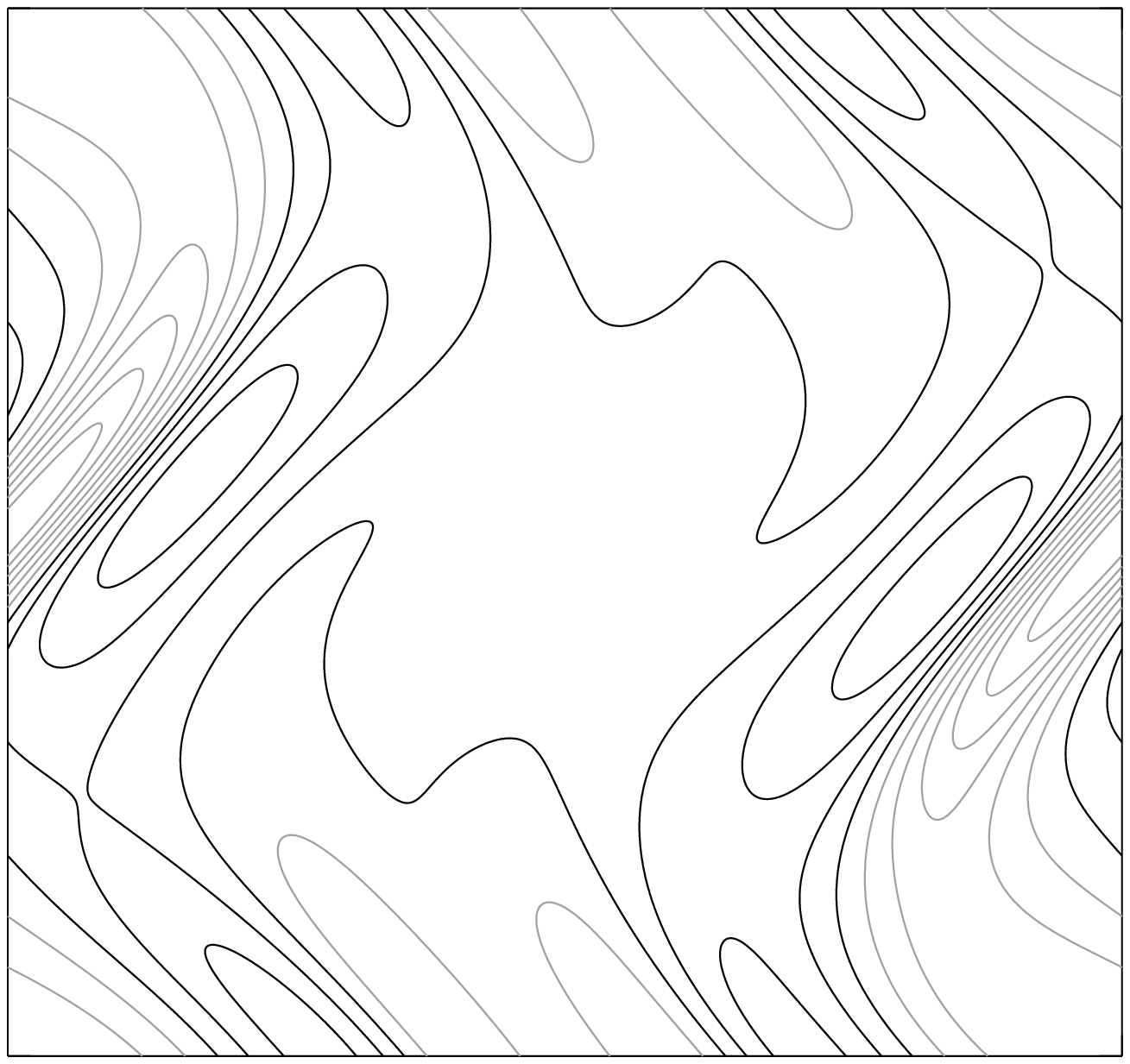,width=45mm}~\psfig{file=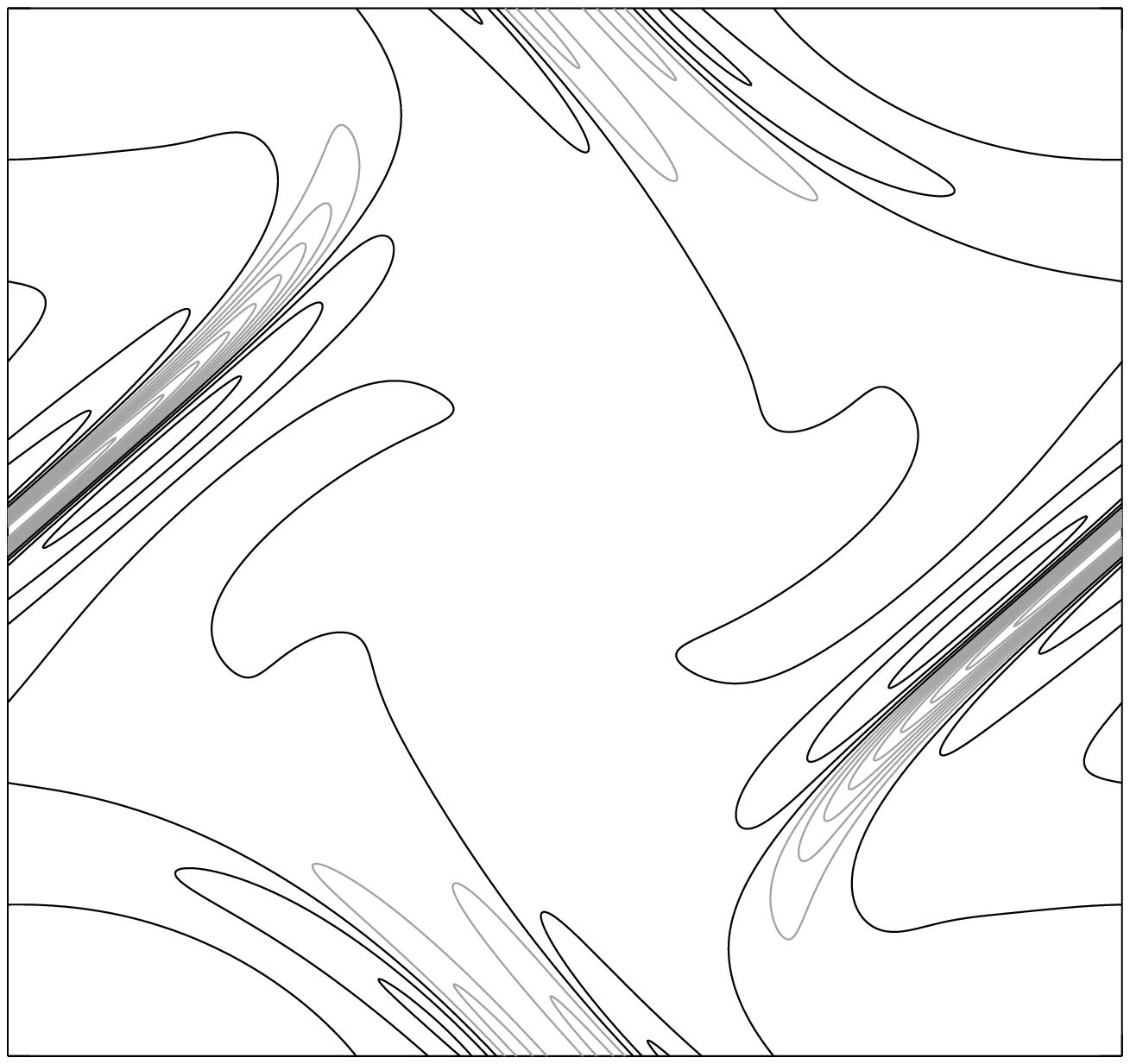,width=45mm}}

\vskip7mm
\centerline{\psfig{file=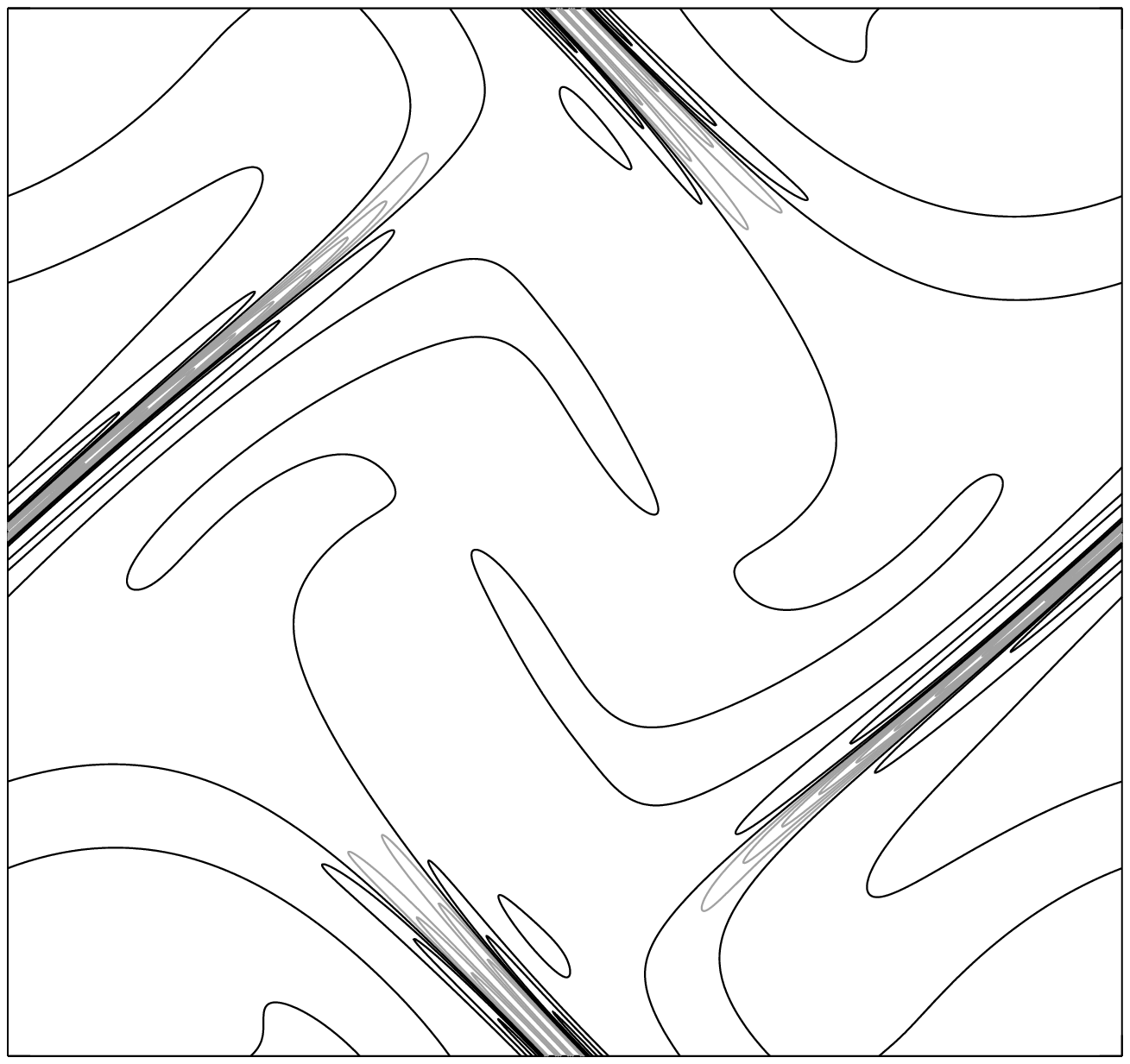,width=45mm}~\psfig{file=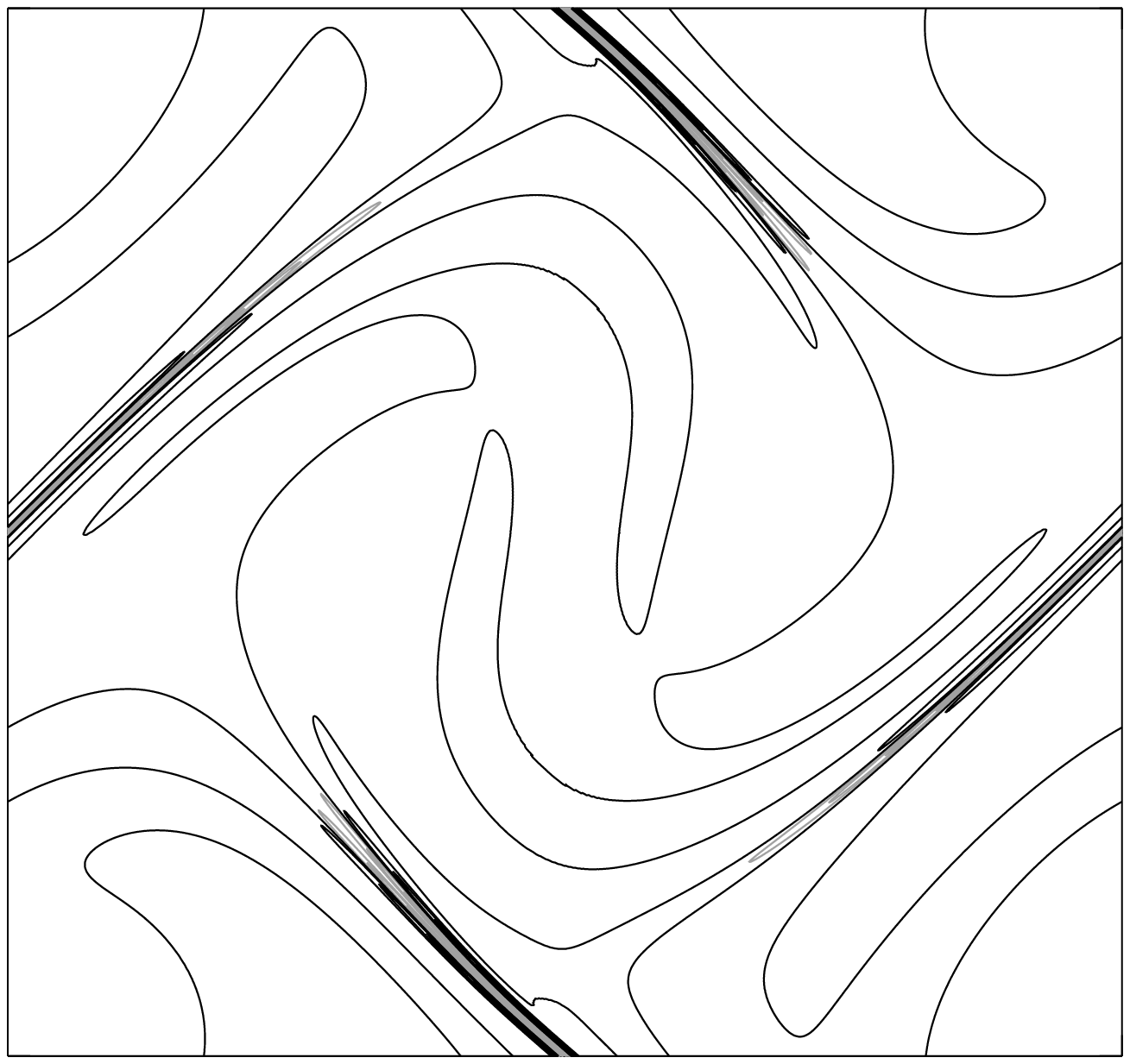,width=45mm}~\psfig{file=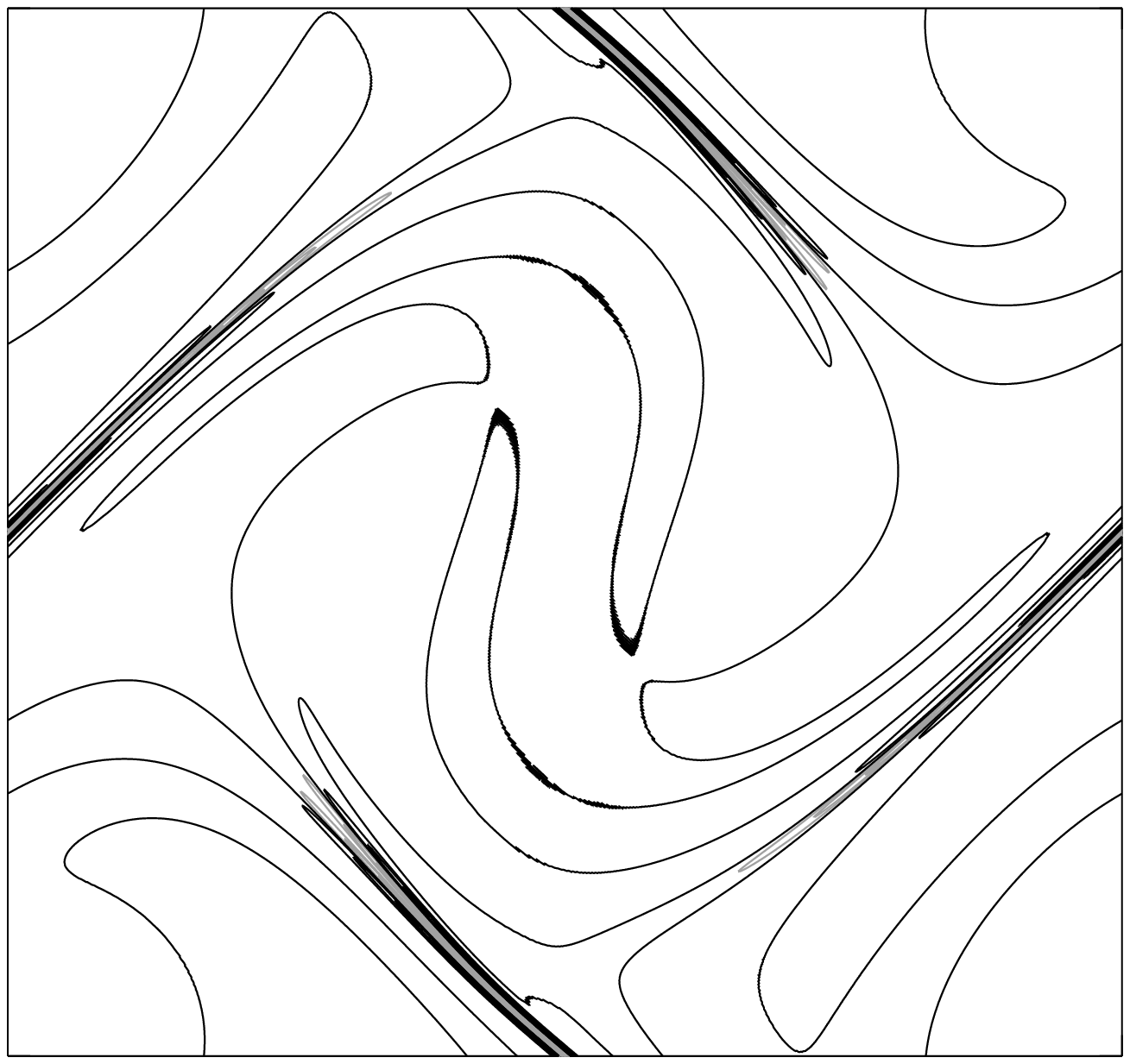,width=45mm}}

\vspace*{-98.5mm}\hskip26mm$t=0$\hskip37mm$t=1$\hskip37mm$t=2$

\vspace*{44.5mm}\hskip26mm$t=3$\hskip37mm$t=4$\hskip37mm$t=4.1$

\vskip42mm
\caption{Isolines of the Laplacian of the vorticity (step 1 for $t=0$;
2 for $t=1$; 10 for $t=2$; 60 for $t=3$; 500 for $t=4$ and 4.1), computed
by CL8 with the resolution of $1024^2$ harmonics for the 4-mode initial condition.}
\label{figisol}\end{figure}

The evolutions of the vorticity and of its Laplacian for the 4-mode initial
flow are shown in Figs.~\ref{figisov} and~\ref{figisol}, respectively. We also
performed high-resolution ($8192^2$) simulations for the same initial flow.
The flow organises itself into a periodic array of large eddies
centred at the centre and corners of the displayed periodicity
box. Fine structures involving large vorticity
gradients develop at the boundaries between these large eddies.

At $t=2$, we see in plots of vorticity isolines, and more clearly
in the isolines of its Laplacian, two thin vortical sheets
(related by the central symmetry of the flow), which are roughly inclined
at +45$^\circ$ with
respect to the x-axis and two other ones, roughly inclined at -45$^\circ$,
which are not quite as fine. Around $t=3.7$, the former and the latter have
thinned to
comparable scales. The thickness of such structures can be interpreted
as a rough measure of the distance to complex-space singularities,
because it determines the strength of spatial derivatives in the
real-space domain. So, here we have a pair of singularities (actually,
singular manifolds) at complex-conjugate spatial locations that
competes with another similar pair, with eventually the one pair
furthest from the real domain catching up with the one which was
closer. This singularity catch-up explains the bending of the plot in
Fig.~\ref{fdel} around $t=3.7$. Another feature shown in
Fig.~\ref{figisol}(b), after the singularity catch-up, are oscillations
superimposed on a roughly exponential energy spectrum. Detailed
examination of contours of the Laplacian of vorticity at the highest
resolution ($8192^2$, not shown) reveals that the thin sheets seen in
Fig.~\ref{figisol}, which are inclined at -45$^\circ$ have a doublet structure, while those inclined at +45$^\circ$ have a singleton structure. The doublet
structure gives rise to oscillatory interference in the Fourier transform.

We now show that the CL results
agree with those obtained by traditional methods for both
test flows within the time intervals where spatial truncations effects
are negligibly small. Table~\ref{tb5} shows the vorticity discrepancies between
results of different algorithms, namely, the maximum over space of the difference
of the vorticity calculated with two different codes. This is done
for the two test flows at different spatial
resolutions and different output times. We also performed energy and enstrophy
conservation tests, whose results
are presented in Tables~\ref{tb3} and \ref{tb4},
respectively. (We will return to conservation issues in Section~\ref{ss:tygers}.) The main
result is that the vorticity discrepancies between the CL results and those
obtained by the standard Eulerian RK4 method are very small, around
$10^{-10}$\,--\,$10^{-12}$, except for later times, when spatial truncation
errors become important. Note that the discrepancy with the RK4 results
is roughly the same for CL8, CL16 and CL24. Choosing
between these methods is thus only a matter of efficiency,
as discussed in the next section.

\begin{table}[p]
\caption{Vorticity discrepancies between solutions (maximum over
the periodicity box of the absolute value of the difference), computed by
different methods for the test flows at various times and resolutions.}
\center\begin{tabular}{|c|c|c|c|c|c|c|}\hline
Run\vphantom{$|^|_|$}&$t$ &CL8-RK4 &CL16-RK4 &CL24-RK4& ET8-RK4 &RK2-RK4\\\hline
4-mode,\vphantom{$|^|_|$}&1&$1.5\!\cdot\!10^{-12}$ & $1.5\!\cdot\!10^{-12}$ &$1.5\!\cdot\!10^{-12}$ &
$1.5\!\cdot\!10^{-12}$ &$10^{-9}$\\
$1024^2$&3&$5\!\cdot\!10^{-11}$ & $2\!\cdot\!10^{-11}$ & $3\!\cdot\!10^{-11}$ & $5\!\cdot\!10^{-12}$ &$1.5\!\cdot\!10^{-9}$\\
&5&$6\!\cdot\!10^{-5}$ & $5\!\cdot\!10^{-5}$ &$5\!\cdot\!10^{-5}$ & $2\!\cdot\!10^{-10}$ &$3\!\cdot\!10^{-9}$\\\hline
4-mode,\vphantom{$|^|_|$}&1&$3\!\cdot\!10^{-12}$ & $5\!\cdot\!10^{-12}$ & $8\!\cdot\!10^{-12}$
&$2\!\cdot\!10^{-12}$ &$10^{-8}$\\
$2048^2$&3&$1.5\!\cdot\!10^{-10}$ & $7\!\cdot\!10^{-11}$& $10^{-10}$&$3\!\cdot\!10^{-11}$&$6\!\cdot\!10^{-8}$\\
&5& $8\!\cdot\!10^{-4}$ & $2\!\cdot\!10^{-4}$ & $4\!\cdot\!10^{-4}$ &$4\!\cdot\!10^{-8}$ &$4\!\cdot\!10^{-6}$\\\hline
random,\vphantom{$|^|_|$}&0.2&$3\!\cdot\!10^{-10}$ & $2\!\cdot\!10^{-10}$ &$2\!\cdot\!10^{-10}$ &
$2\!\cdot\!10^{-10}$ &$3\!\cdot\!10^{-7}$\\
$2048^2$&0.6&$7\!\cdot\!10^{-8}$ & $2\!\cdot\!10^{-8}$ & $2\!\cdot\!10^{-9}$ & $2\!\cdot\!10^{-9}$ &$10^{-6}$\\
&1&$2\!\cdot\!10^{-3}$ & $2\!\cdot\!10^{-3}$ &$2\!\cdot\!10^{-3}$ & $2\!\cdot\!10^{-6}$ &$4\!\cdot\!10^{-5}$\\\hline
\end{tabular}\label{tb5}

\vspace*{\baselineskip}
\caption{Energy conservation errors.}
\center\begin{tabular}{|c|c|c|c|c|c|c|c|}\hline
Run\vphantom{$|^|_|$}& $t$ &CL8 & CL16 & CL24 & RK2 & RK4 &ET8\\\hline
4-mode,\vphantom{$|^|_|$}&1&$-3\!\cdot\!10^{-14}$ & $10^{-15}$ &$<10^{-15}$& $10^{-11}$ &$2\!\cdot\!10^{-14}$&$<10^{-15}$\\
$1024^2$&3&$-2\!\cdot\!10^{-14}$ & $7\!\cdot\!10^{-14}$& $2\!\cdot\!10^{-13}$& $-2\!\cdot\!10^{-11}$ &
$4\!\cdot\!10^{-14}$&$<10^{-15}$\\
&5&$-5\!\cdot\!10^{-9}$ & $-10^{-8}$& $-10^{-8}$& $-4\!\cdot\!10^{-11}$ &$7\!\cdot\!10^{-14}$&$2\!\cdot\!10^{-15}$\\\hline
4-mode,\vphantom{$|^|_|$}&1& $-3\!\cdot\!10^{-14}$& $10^{-15}$ &$10^{-15}$ &$10^{-11}$ &$10^{-15}$&$<10^{-15}$\\
$2048^2$&3& $-2\!\cdot\!10^{-14}$& $10^{-15}$ &$10^{-15}$ &
$-2\!\cdot\!10^{-11}$ &$<10^{-15}$&$-2\!\cdot\!10^{-15}$\\
&5& $-6\!\cdot\!10^{-11}$& $-10^{-10}$& $-10^{-10}$&
$-4\!\cdot\!11^{-11}$ &$10^{-15}$&$-5\!\cdot\!10^{-15}$\\\hline
random,\vphantom{$|^|_|$}&0.2& $-10^{-15}$& $-2\!\cdot\!10^{-15}$ &$-2\!\cdot\!10^{-15}$ &
$-7\!\cdot\!10^{-11}$ &$10^{-14}$&$-2\!\cdot\!10^{-15}$\\
$2048^2$&0.6& $-4\!\cdot\!10^{-14}$& $-4\!\cdot\!10^{-14}$ &$-3\!\cdot\!10^{-14}$ &
$2\!\cdot\!10^{-11}$ &$3\!\cdot\!10^{-14}$&$-3\!\cdot\!10^{-14}$\\
&1& $-4\!\cdot\!10^{-11}$& $10^{-10}$& $-2\!\cdot\!10^{-10}$&
$6\!\cdot\!11^{-10}$ &$-3\!\cdot\!10^{-13}$&$-3\!\cdot\!10^{-13}$\\\hline
\end{tabular}\label{tb3}

\vspace*{\baselineskip}
\caption{Enstrophy conservation errors.}
\center\begin{tabular}{|c|c|c|c|c|c|c|c|}\hline
Run\vphantom{$|^|_|$}& $t$ &CL8 & CL16 & CL24 & RK2 & RK4 &ET8\\\hline
4-mode,\vphantom{$|^|_|$}&1&$-4\!\cdot\!10^{-14}$ &$<10^{-15}$&$<10^{-15}$&$3\!\cdot\!10^{-11}$ &$2\!\cdot\!10^{-14}$&$<10^{-15}$\\
$1024^2$&3&$-2\!\cdot\!10^{-12}$ &$-10^{-12}$& $-3\!\cdot\!10^{-12}$&
$6\!\cdot\!10^{-11}$ &$3\!\cdot\!10^{-15}$&$-3\!\cdot\!10^{-15}$\\
&5&$-2\!\cdot\!10^{-6}$ & $-10^{-6}$ & $-2\!\cdot\!10^{-6}$&
$4\!\cdot\!10^{-10}$ &$-2\!\cdot\!10^{-12}$&$-5\!\cdot\!10^{-15}$\\\hline
4-mode,\vphantom{$|^|_|$}&1& $-4\!\cdot\!10^{-14}$&$<10^{-15}$&$<10^{-15}$&$3\!\cdot\!10^{-11}$ &$10^{-15}$&$<10^{-15}$\\
$2048^2$&3& $-3\!\cdot\!10^{-14}$&$-8\!\cdot\!10^{-14}$& $-10^{-14}$&
$6\!\cdot\!10^{-11}$ &$5\!\cdot\!10^{-15}$&$-6\!\cdot\!10^{-15}$\\
&5& $-2\!\cdot\!10^{-8}$& $-5\!\cdot\!10^{-9}$& $-7\!\cdot\!10^{-9}$&
$4\!\cdot\!10^{-10}$ &$10^{-13}$&$-4\!\cdot\!10^{-14}$\\\hline
random,\vphantom{$|^|_|$}&0.2& $10^{-14}$&$-2\!\cdot\!10^{-14}$&$-2\!\cdot\!10^{-14}$&
$-3\!\cdot\!10^{-10}$ &$10^{-13}$&$-2\!\cdot\!10^{-14}$\\
$2048^2$&0.6& $-6\!\cdot\!10^{-13}$&$-2\!\cdot\!10^{-13}$& $-2\!\cdot\!10^{-13}$&
$4\!\cdot\!10^{-9}$ &$10^{-13}$&$-2\!\cdot\!10^{-13}$\\
&1& $-10^{-8}$& $-10^{-8}$& $2\!\cdot\!10^{-9}$&
$3\!\cdot\!10^{-9}$ &$-8\!\cdot\!10^{-12}$&$-2\!\cdot\!10^{-12}$\\\hline
\end{tabular}\label{tb4}\end{table}

\begin{table}[h]
\caption{CPU times needed to reach various output times for
the codes applied to the test flows.}
\center\begin{tabular}{|c|c|c|c|c|c|c|c|}\hline
Run\vphantom{$|^|_|$}& $t$ & CL8& CL16& CL24& RK2 & RK4 & ET8\\\hline
4-mode,\vphantom{$|^|_|$}&1& 1m1s & 27s & 32s& 92m50s & 7m33s& 16m22s\\
$1024^2$&3& 2m56s & 1m15s & 1m15s & 278m20s & 22m26s & 48m1s\\
&5& 4m55s& 2m12s & 2m08s & 464m0s& 37m22s& 81m50s\\\hline
4-mode,\vphantom{$|^|_|$}&1& 4m20s & 1m59s & 2m17s & 402m20s& 65m01s & 140m10s\\
$2048^2$&3& 12m30s & 5m24s & 5m18s & 1206m50s & 194m20s& 420m20s\\
&5& 20m56s& 9m23s& 9m07s & 2011m20s & 323m40s &700m10s\\\hline
random,\vphantom{$|^|_|$}&0.2&3m05s & 1m29s & 1m32s & 81m40s& 16m19s & 34m53s\\
$2048^2$&0.6&8m42s & 3m57s & 4m32s & 244m50s & 48m20s& 103m40s\\
&1& 13m48s & 6m25s & 6m51s & 408m20s & 81m10s& 163m40s\\\hline
\end{tabular}\label{tb2}\end{table}

\subsection{Efficiency of the CL method}\label{ss:efficiency}

\begin{table}[b]
\caption{CPU times for different choices of the accuracy $\varepsilon$.}
\center\begin{tabular}{|c|c|c|c|c|c|c|}\hline
$\varepsilon$ & $10^{-12}$ & $10^{-11}$ & $10^{-10}$ & $10^{-9}$ &
$10^{-8}$\\\hline
CPU time & 295s & 224s & 169s & 135s &103s\\\hline
\end{tabular}\label{tchangeeps}
\end{table}

Before discussing the relative performances of the various algorithms,
let us note that most of the computations were made on a computer -- marketed
in 2011 -- that has an Intel Core i7-3960X processor (64 bit data width,
frequency 3.3 GHz up to 3.9 GHz in the Turbo regime, 15 Mbyte L3 cache).
It has 6 cores with hyper-threading, which does not actually matter for our
application, since we have so far experimented only with sequential versions
of the codes. Parallelising such a code is standard -- at least for the
Lagrangian time-stepping phase -- and
requires the same tools as for the parallelisation of any spectral code.

The CPU timings for the integration to various output times for our codes
are shown in Table~\ref{tb2}. The timings for the Lagrangian CL codes are way
below those for the best Eulerian code: at the resolution of $2048^2$ harmonics,
CL24 needs at least
an order of magnitude less CPU time than the commonly used RK4; the ET8
algorithm is at least 2 times slower than RK4, and RK2 is much slower.

We also observe that CL16 and CL24 have comparable efficiency and are
significantly faster than CL8.
Hence the optimal number of terms in the time-Taylor series is around $S=20$.
This is consistent with the arguments of Section~\ref{ss:optimal}
that give bounds for the optimal order, involving an unknown constant $A$.

All the results given so far are for the choice $\varepsilon=10^{-12}$
for the accuracy parameter that controls the truncation
errors on the time-Taylor series (see \rf{estet0}). It is also of
interest to measure the CPU time dependence on $\varepsilon$. For the
four-mode initial conditions at resolution $1024^2$ integrated with
CL8 from $t=0$ to $t=5$ the CPU time is shown in Table~\ref{tchangeeps}.
The results are roughly consistent with an eighth-power dependence of the
accuracy on the time step $\Delta t$, expected for CL8.

Let us finally discuss the durations of the runs for the extreme
resolution of $8192^2$ harmonics, used to calculate the radii $\delta$
of the cylinder of analyticity, shown in Fig.~\ref{fdel}. Since in a
Lagrangian code the time step does not depend on the spatial
resolution, we expect even larger CPU gains than reported above. This
is indeed the case: at the resolution of $8192^2$ harmonics and in
double precision, the integration from $t=0$ to $t=5$ with the
four-mode initial condition and our CL8 algorithm took 5.5 hours of
CPU on the computer mentioned at the beginning of this section,
whereas on the same computer, the standard RK4 runs for 564 hours, a
little over one hundred time slower. The vorticity discrepancies
between the two outputs are less than $10^{-13}$ up to $t=4.5$, as
expected from just rounding errors, and increase to $4.5\times
10^{-12}$ at $t=5$, which is consistent with the spatial truncation
errors estimated in Section~\ref{ss:flow-methods}. Note that the
truncation order $S=8$, used for the Lagrangian time-Taylor series is
well below the optimal order, which is around $S=20$. Using another
computer with a larger memory, allowing to run CL20, we timed CL20 and
RK4 over identical time intervals and found that the former runs about
200 faster than the latter. Of course the exact ratios depend not just
on memory size, but also on the presence or not of high-wavenumber
modes in the initial condition and on the precise method of
parallelisation. Such issues will be discussed elsewhere.

\subsection{Galerkin truncation artefacts}\label{ss:tygers}

\begin{figure}[t]
\centerline{\psfig{file=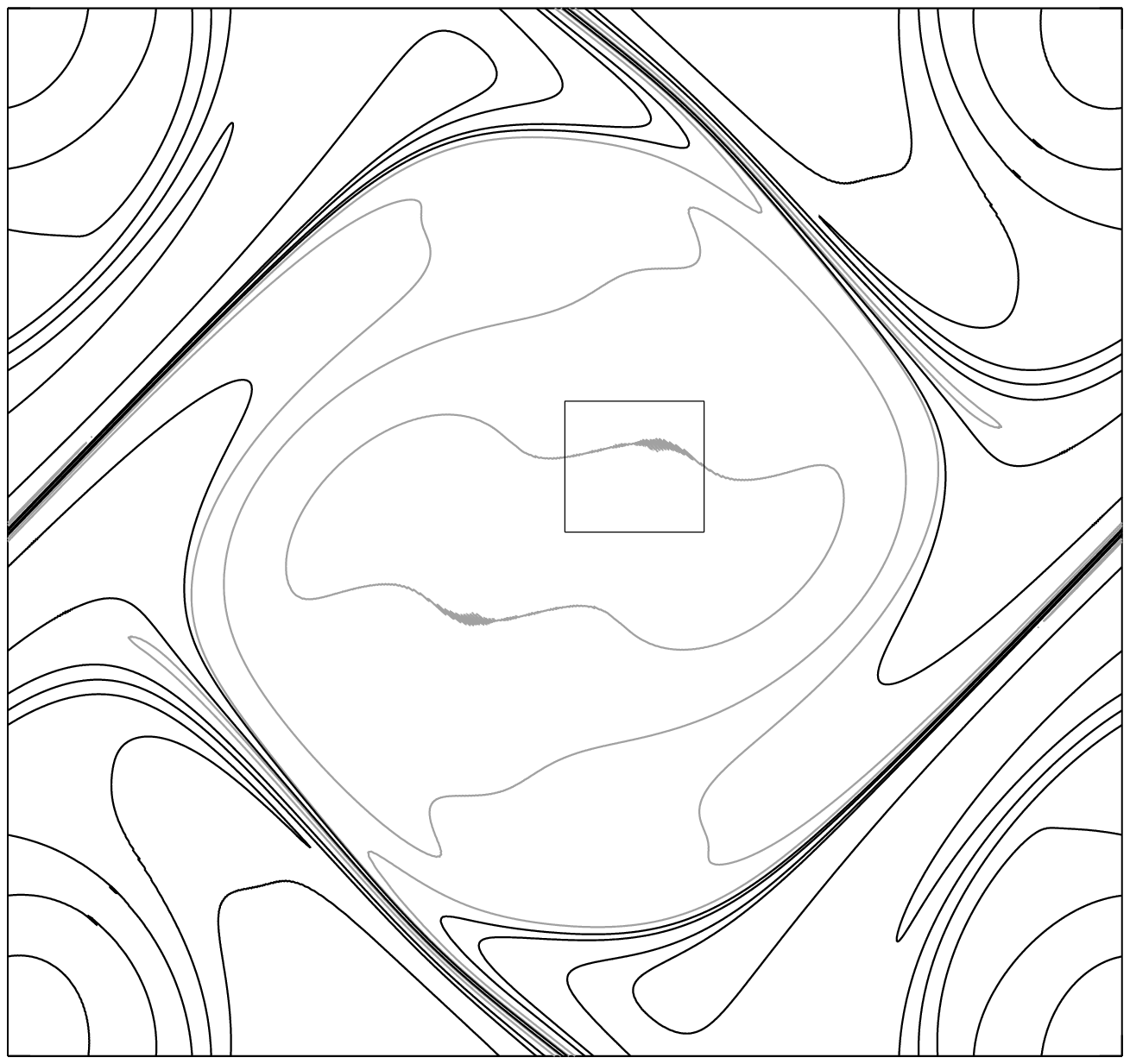,width=45mm}\hspace*{27mm}\psfig{file=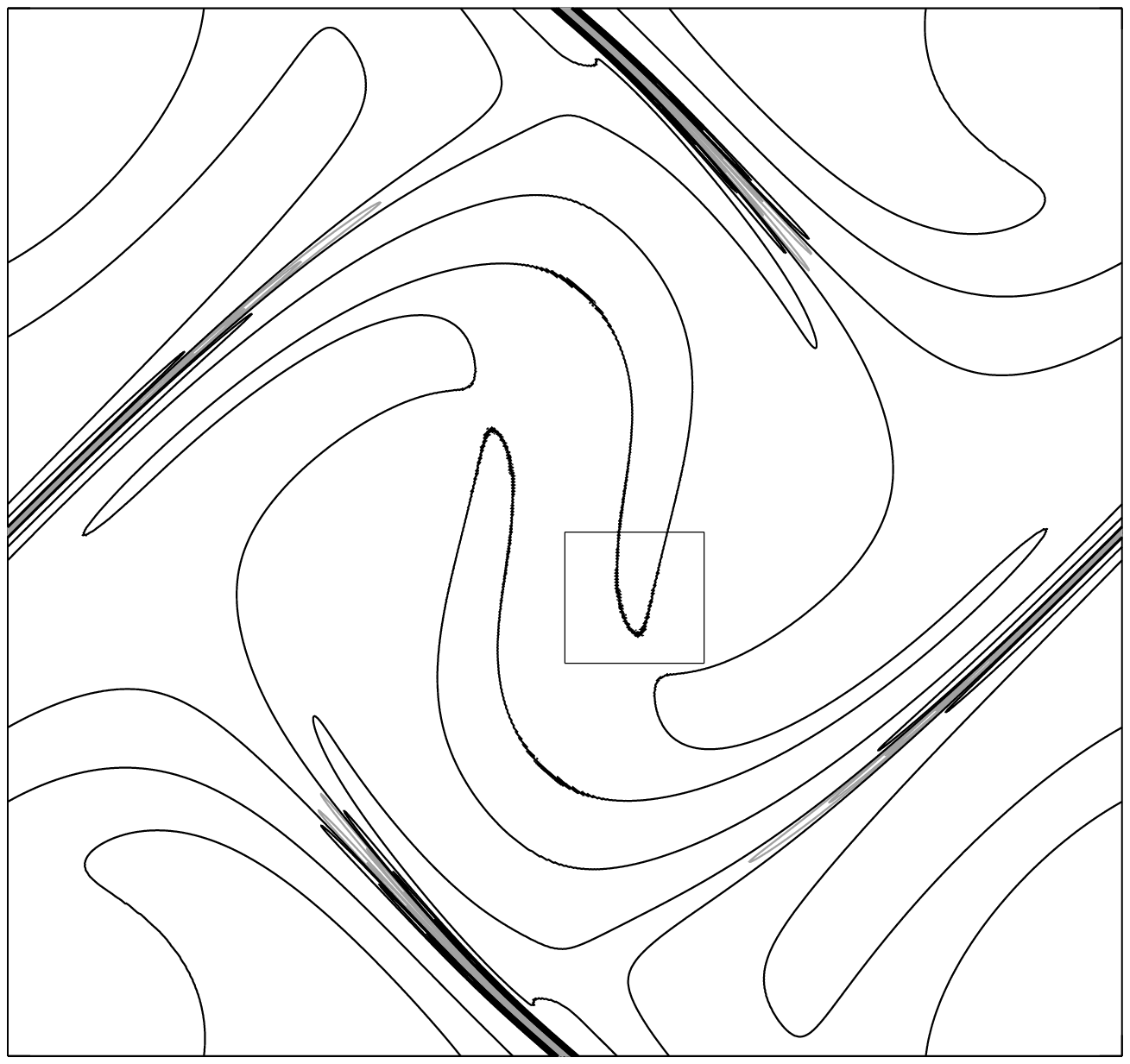,width=45mm}}

\vspace*{-37mm}\hskip2cm(a)\hskip67mm(b)

\vspace*{37mm}
\centerline{\psfig{file=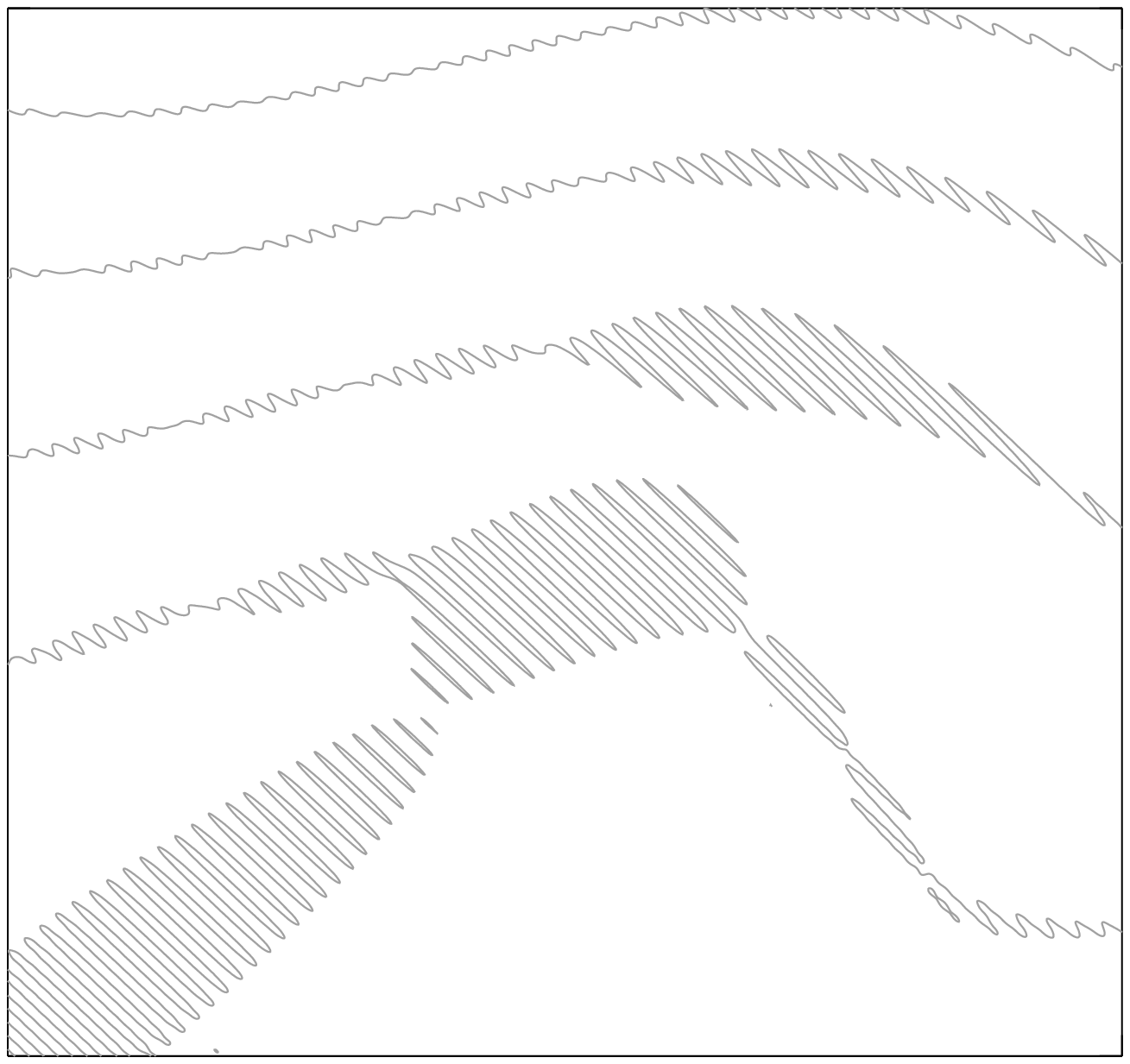,width=45mm}\hspace*{27mm}\psfig{file=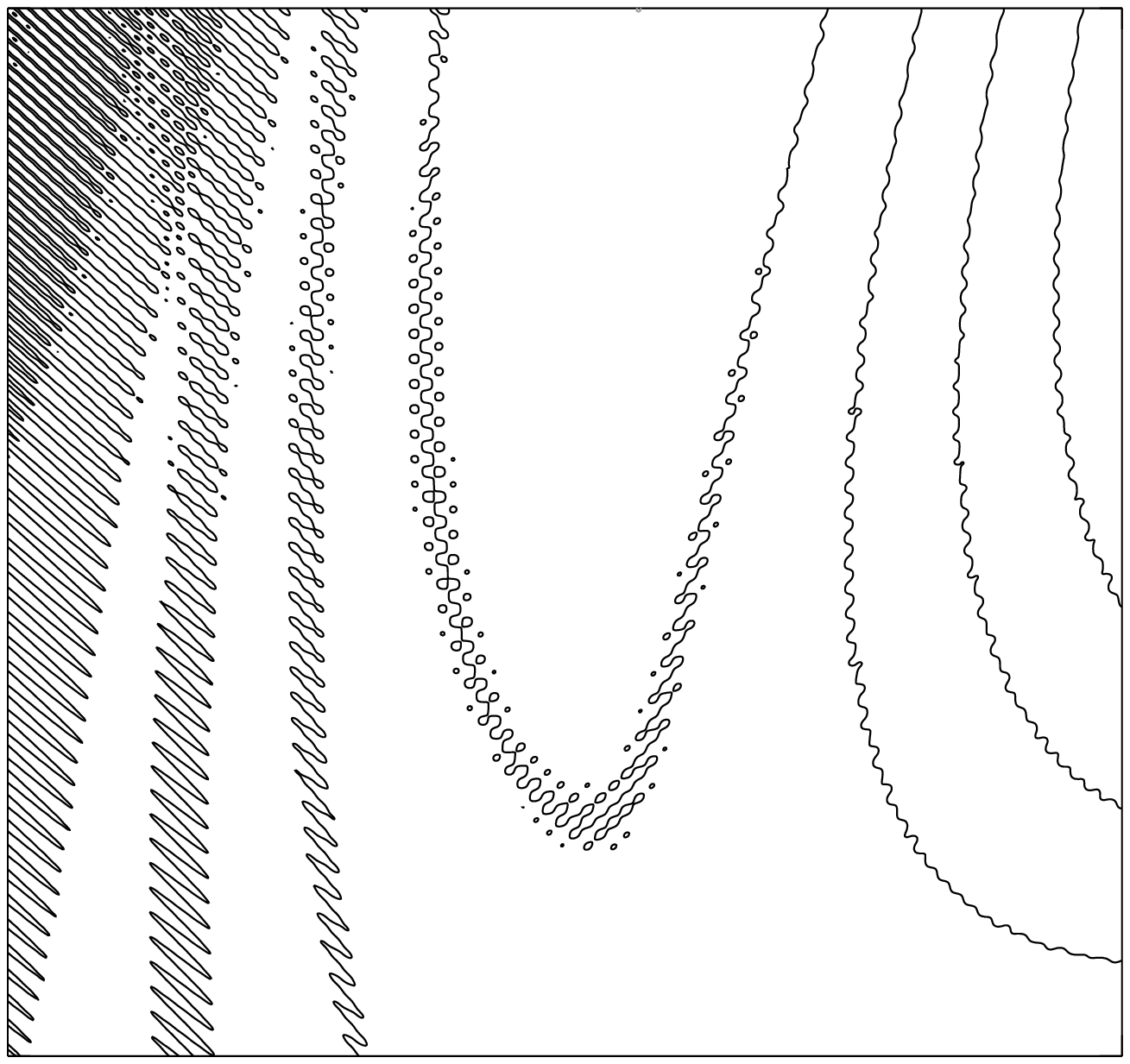,width=45mm}}

\vspace*{-37mm}\hskip2cm(c)\hskip67mm(d)

\vspace*{-15mm}\hspace*{24mm}{\Large${5\over4}\pi$\hskip69mm$\pi$}

\vspace*{34mm}\hspace*{27mm}{\Large$\pi$\hskip65.5mm${3\over4}\pi$}

\vspace*{-2mm}
\hspace*{3cm}{\Large$\pi$\hskip4cm${5\over4}\pi$
\hspace*{2cm}$\pi$\hskip4cm${5\over4}\pi$}

\vspace*{4mm}
\caption{Fourier truncation artefacts (tygers) in Eulerian simulations for the
4-mode initial condition with $1024^2$ harmonics by RK4. Isolines
of the vorticity (step 0.5) at $t=5$ (a) and of its
Laplacian (step 500) at $t=4$ (b). Blow-ups of these figures:
vorticity (step 0.05) at $t=5$ (c) and its Laplacian (step 0.5) at $t=4$ (d).}
\label{isoran}\end{figure}

In spectral simulations of ideal flow, one typically runs out of
spatial resolution after some time of integration. As pointed out in
Section~\ref{ss:init}, for analytic initial data, spatial truncation
errors roughly behave as $\e^{-\delta(t)k_{\max}}$, which ceases to be
small when complex-space singularities approach the real domain to
within a few spatial grid mesh sizes. In Eulerian simulations,
localised artefacts, in the form of spurious oscillations of contour
lines, called ``tygers,'' are then appearing in unexpected
locations \cite{tyger}. These are produced by truncation-generated
waves interacting resonantly at places where the fluid velocity
matches that of developing fine-scale structures. For the 4-mode
initial condition, the phenomenon is illustrated in Fig.~\ref{isoran},
which shows Eulerian RK4 simulations with the resolution of $1024^2$
harmonics at the earliest times when tygers become visible. As is
typical for small-scale objects, the tygers become at first
conspicuous in contour levels of the Laplacian of the vorticity,
around $t=4$ when $\delta$ is about 5 spatial mesh sizes, before they
are seen in the vorticity itself, around $t=5$ when $\delta$ is about
1.5 spatial mesh sizes. For the random initial condition at the same
resolution (not shown) similar phenomena happen around $t=0.71$ and $t=1$.

Do analogous artefacts emerge in our Cauchy--Lagrangian simulations?
Contours of the Laplacian of the vorticity show indeed similar spurious oscillations
(for the same flow and the same spatial resolution) around $t=4.1$ (see
the last panel in Fig.~\ref{figisol}). However, the vorticity shows nothing
comparable until $\delta$ drops well below one spatial mesh.

Hence, Galerkin truncation artefacts are sharply reduced in our CL simulations,
compared to Eulerian ones with the same resolution. This is surprising given
that the Eulerian and Lagrangian coordinates are resynchronised at each time step,
so that not much difference in behaviour and deterioration in time of solutions
would be expected. Part of this reduction could be due to the smoothing effect
of the interpolation stage, also observed in the semi-Lagrangian approach
of \cite{mc}. There is, however, another possibility: solutions truncated
to a finite number of spatial Fourier modes (generally referred to as Galerkin
truncated) behave differently in Eulerian and Lagrangian
coordinates. In Eulerian coordinates, it is easily shown that quadratic
invariants (such as energy and enstrophy, in 2D, or energy and helicity,
in 3D) remain invariant after Galerkin truncation. Hence, as observed
in \cite{tyger}, energy and enstrophy, which
in the absence of truncation would be transferred to smaller scales,
have to be transferred somewhere else and get pushed into the
tygers. In Lagrangian coordinates, the behaviour is different:
The energy is $E=(1/2)\int|\dot\bx(\ba,t)|^2\,d^2\ba$, where the integral is
over the $[0,2\pi]\times[0,2\pi]$ periodicity domain. Before truncation,
the energy is exactly conserved in any time interval $\Delta t$. After
Galerkin truncation, it is only conserved up to second order in
$\Delta t$, and the same holds for enstrophy.

\section{Rounding noise in high-order Lagrangian vs Eulerian time-marching}
\label{s:rounding}

Some challenging open problems in hydrodynamics, such as the issue of
blow-up for 3D flow \cite{gi} or for axisymmetric 2D
flow \cite{lh}, may require extremely accurate simulations in order to
distinguish genuine blow-up from simulation artefacts. This means that
spatial and temporal truncation errors must be kept very small and the
precision very high.

For flow in a domain of sufficiently simple
geometry (foremost space-periodic flow), the most accurate methods
for handling the spatial discretisation are of the spectral or pseudospectral
type~\cite{go}. The reason is that, for spatially analytic flow, the Fourier
coefficients fall off exponentially at large wave numbers \cite{Fhouches83,ssl}
and, thus, the truncation error drops exponentially when increasing
the resolution. Unfortunately, (pseudo-)spectral simulations hardly ever use
temporal discretisation (time-marching) of order higher than 4, the most
frequent being Runge--Kutta of fourth order (RK4). We have therefore a striking
imbalance between temporal and spatial discretisation. Of course,
when using Eulerian methods, one generally has to satisfy
the Courant--Friedrichs--Lewy \cite{cfl} condition:
\BE k_{\max}U_{\max}\Delta t<C_{\rm s},\EE{CFL}
where $k_{\max}$ is the maximum wave number (with dealiasing taken into
account), $U_{\max}$ is the maximum velocity in physical space,
and $C_{\rm s}$ is an order unity constant, depending on the numerical scheme.
As a consequence, high spatial resolution requires a very small time
step, so that a fourth-order scheme might be deemed sufficiently accurate.

\begin{figure}[t!]
\begin{picture}(132,36)(0,0)
\put(10,0){\psfig{file=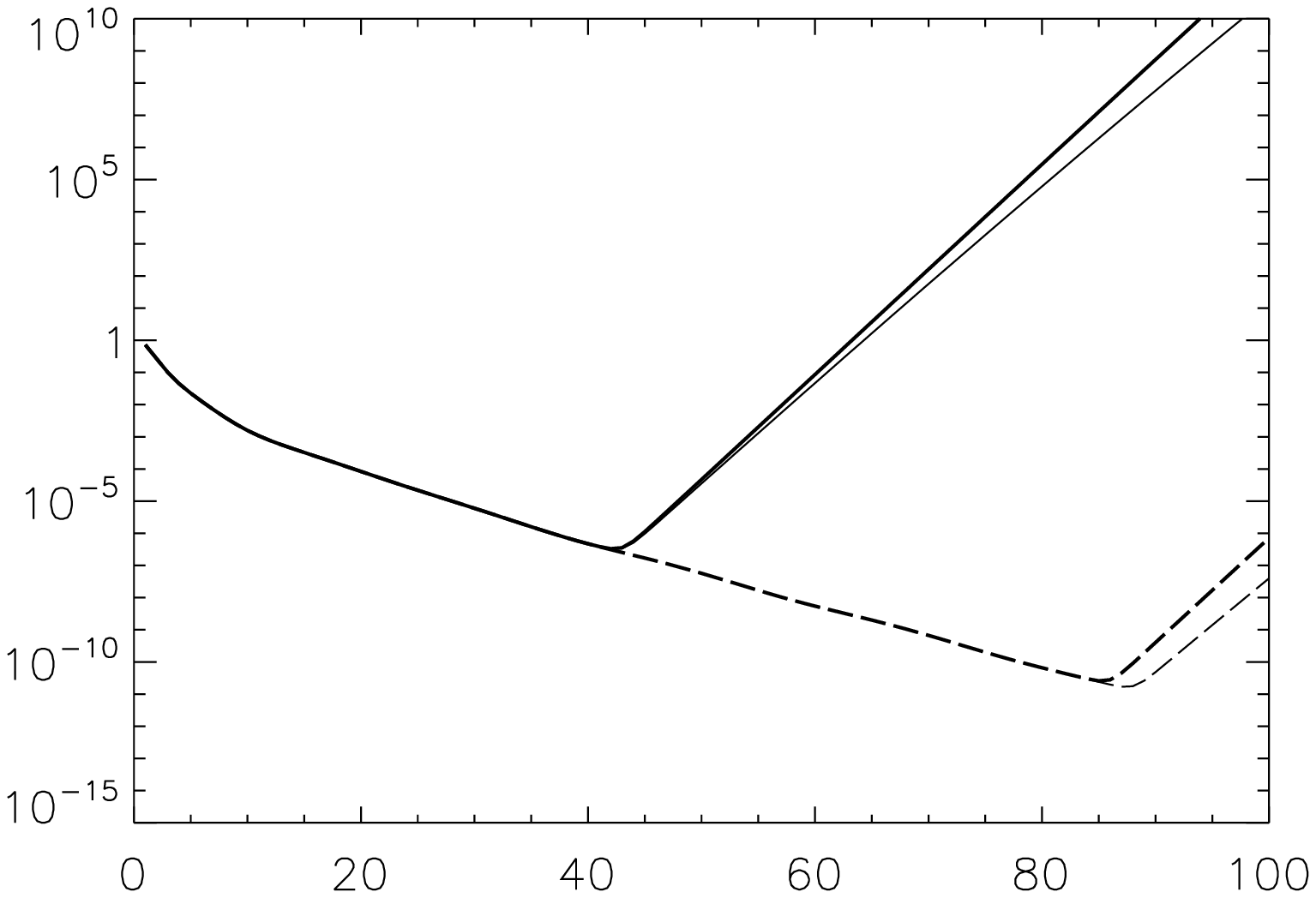,width=67mm}}
\put(4,29){(a)}
\put(76,0){\psfig{file=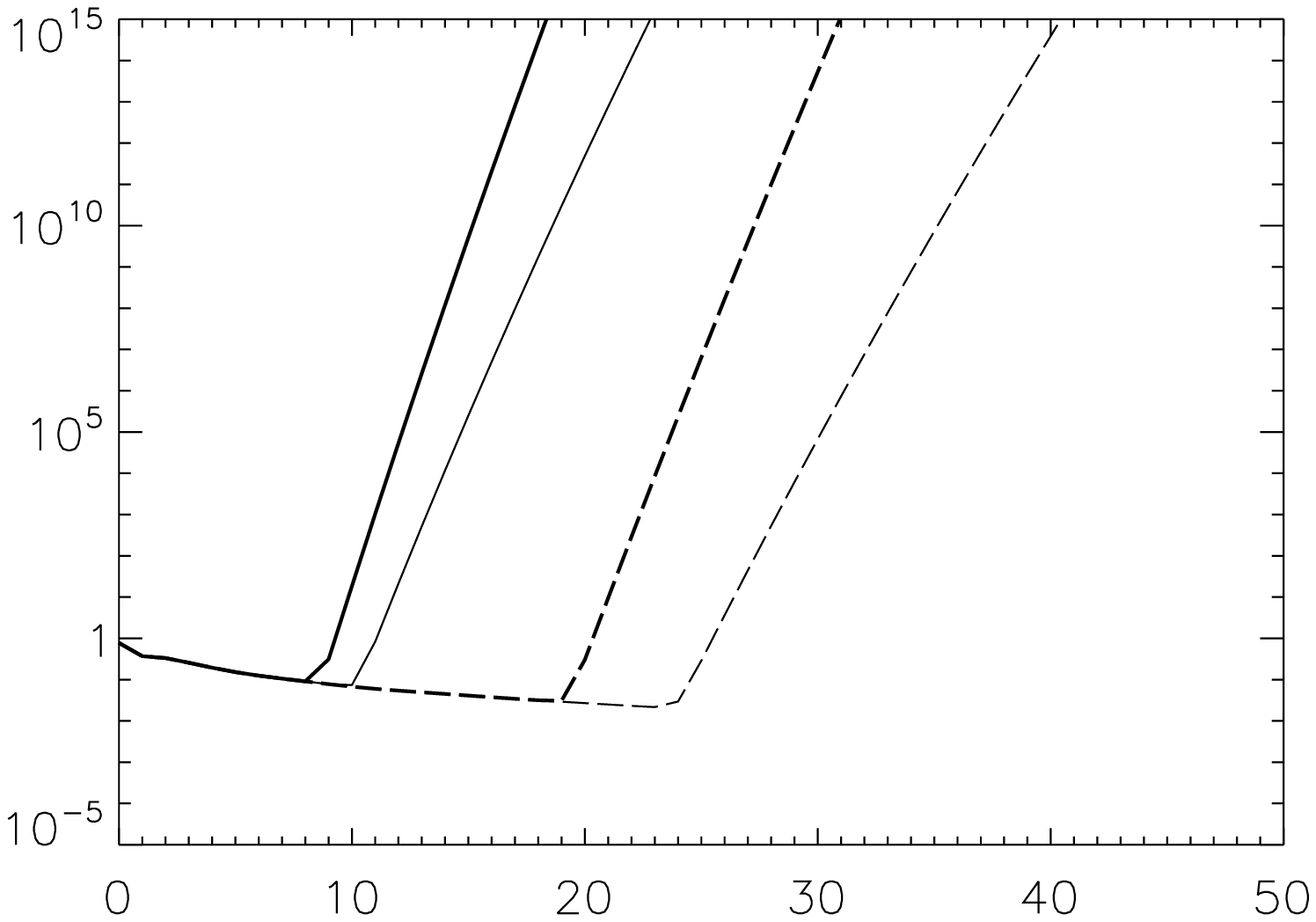,width=67mm}}
\put(70,29){(b)}
\put(16,29){$\|\bxi^{(s)}\|$}
\put(82,29){$\|\omega_s\|$}
\put(53,4){$s$}
\put(119,4){$s$}
\end{picture}
\caption{L$^2$ norms of coefficients of the displacement, $\|\bxi^{(s)}\|$, and
vorticity, $\|\omega_s\|$, in the Lagrangian (a) and Eulerian~(b) time-Taylor
series, respectively, computed at $t=0$ for the 4-mode initial
condition. Solid line: double precision, dashed line: quadruple precision.
Thin lines: resolution $512^2$ harmonics, thick lines: $1024^2$ harmonics.}
\label{coeft0}\end{figure}

The situation is drastically different when working in Lagrangian coordinates.
The time step is not controlled by a CFL condition any more, but solely by
the radius of convergence of the time-Taylor expansion for the
Lagrangian map. This radius is roughly the inverse of the largest velocity
gradient present and does not depend on spatial resolution. Hence, in
principle, there is no reason to limit oneself to low-order time-marching.

We now report results of experimenting with very high-order time-Taylor
expansions in both Lagrangian and Eulerian coordinates. We will
see that the severest limitation stems from the growth of rounding errors.

In Fig.~\ref{coeft0}, the L$^2$ norms of time-Taylor coefficients in Lagrangian
coordinates, obtained by the Cauchy--Lagrangian method, and in Eulerian
coordinates, obtained by the Eulerian Taylor method, are shown for the 4-mode
initial condition, computed in double and quadruple precision up to the orders
beyond which an unnatural behaviour, dominated by rounding errors, sets in.
Taylor coefficients of the displacement, $\bxi^{(s)}$, and of the
vorticity, $\omega_s$, in the expansions of the Lagrangian and
Eulerian solutions, respectively, have been computed around $t=0$
with two resolutions of $512^2$ and $1024^2$ harmonics.

The most striking feature is that the L$^2$ norms show a range
of orders over which the results display very little dependence on precision
and then a rapid transition to a regime in which drastic differences
of many orders of magnitude are observed. This happens around a ``transition
order'' that depends very much on precision and which is considerably larger
(by about a factor 4) in the Lagrangian case than in the Eulerian case. (There
is also some dependence on the resolution.) The exponential behaviour (straight
line in the lin-log coordinate system used) is not observed at the smallest
values of the order~$s$. Indeed, such exponential behaviour $\sim R^{-s}$
of the coefficient of order $s$ of a power series with radius of convergence
$R$ is typically only an asymptotic property at large $s$.

The large difference between the Lagrangian and the Eulerian
computations in the transition order, where rounding noise takes over,
can be seen as a consequence of a phenomenon discovered by Ebin and Marsden
\cite{em}, which is called the {\it loss of derivatives} in the Eulerian
formulation. In this formulation, each time derivative
$\partial_t$ is accompanied by a space derivative $\bv\cdot\nabla$.
As a consequence, if the initial velocity has only a finite number of space
derivatives, it will also have a finite (generally, the same) number of time
derivatives and, thus, cannot be analytic in time. In Lagrangian coordinates,
this is not the case: time-analyticity of the Lagrangian map, for at least
a finite time, holds if
the initial velocity is slightly smoother than once differentiable (see
Section~\ref{s:background}). Ebin and Marsden, although they did not prove
time-analyticity in the Lagrangian formulation, did show that
the variational formulation of the Euler equation, which is
essentially Lagrangian, allows interpreting it as an ordinary
differential equation in function spaces, such as Sobolev spaces,
possessing a finite number of space derivatives. The proof was given
using the language of differential geometry
on infinite-dimensional Riemannian manifolds. Although somewhat simplified by
Bourguignon and Brezis \cite{bobr}, the proof did not suggest any
numerical implementation. This became possible with the
rediscovery of Cauchy's Lagrangian formulation, enabling us to generate
time-Taylor coefficients of arbitrary orders.

Now, we comment on why loss of derivatives in Eulerian coordinates is
responsible for rounding noise that grows very fast with the
order of the time-Taylor coefficient. The 4-mode initial condition
used for the present study has of course space derivatives of
arbitrary order (it is an entire function). Hence time
derivatives of arbitrary orders at $t=0$ also exist. It is easy to see that
spatial Fourier harmonics of the time-Taylor coefficient of order $s$
(not subjected to Galerkin truncation) vanish beyond a wave
number growing linearly with $s$. (This is true both in Lagrangian and
Eulerian coordinates.) Once the problem is discretised in space and
FFTs are used, those wave vectors which should not be excited
will be populated by rounding noise. In the Eulerian case,
this noise is amplified by the presence of the space derivatives in
the recurrence relations, each yielding a factor $\i\bk$ on the Fourier
coefficients. In the Lagrangian case, the recurrence relation is written
in terms of gradients of time-Taylor coefficients, but otherwise we only
perform additions, multiplications and apply bounded Calderon--Zygmund operators.
Thus, amplification of noise at high wave numbers is not so strong
as in the Eulerian recurrence relations.
A further confirmation of this mechanism is obtained by considering
the relative change with the resolution of the order $s$ at which
the transition to noise takes place. When changing from resolution of $512^2$
to $1024^2$ harmonics, the transition orders change by just a few percent
in the Lagrangian calculations and by nearly 20\% in the Eulerian calculations.
(Of course, the higher the resolution, the more operations and thus the more
accumulation of rounding errors takes place.)

The advantage of Lagrangian time-marching over Eulerian time-marching
from the point of view of rounding noise is not limited to initial conditions
comprising a finite number of non-vanishing harmonics. To show
this, we use a multistep method with the same initial data as above, but where
a very broad band of Fourier harmonics become excited beyond the first time
step. We begin by explaining {\it how we choose our time steps in the Eulerian
calculations}. The initial data considered here are analytic in the space
variables and remain analytic in space and time as long as there is no blow-up
(for ever in 2D) \cite{bb}. At any time $t>0$ there is a finite radius of
convergence $R(t)$ of the Eulerian time-Taylor series giving the solution (say,
the vorticity) at time $t+\Delta t$ in terms of the solution at time $t$.
Certainly, we have to take $\Delta t<R(t)$, but this is not enough, because
of a rounding noise problem, closely linked to what we discussed above.

In Eulerian coordinates, the fastest temporal variation of small-scale eddies
comes from them being swept by the large-scale energy-carrying eddies. Denoting
by $\vartheta(t)$ the amplitude (say, the vorticity) associated with such an
eddy, the sweeping can be modelled by the equation
\hbox{$\partial_t\vartheta+({\bm U}\cdot\nabla)\vartheta=0$,}
where $\bm U$ is a large-scale velocity which can be
taken uniform to leading order (the expansion parameter being the ratio of
scales of the swept eddies and the sweeping eddies). Let us denote
by $\vartheta_\bk$ the Fourier coefficients of the small-scale vorticity.
The sweeping equation has the solution
$\vartheta_\bk(t+\Delta t)=\exp(-\i\bk\cdot{\bm U}\Delta t)\vartheta_\bk(t)$.
Expanding the exponential in a Taylor series in $\Delta t$, we find:
\BE\vartheta_\bk(t+\Delta t)=\vartheta_\bk(t)\sum_{s=0}^\infty{(-\i\,
\bk\cdot{\bm U}\,\Delta t)^s\over s!}.\EE{expdev}

\begin{figure}[t!]
\begin{picture}(132,36)(0,0)
\put(10,0){\psfig{file=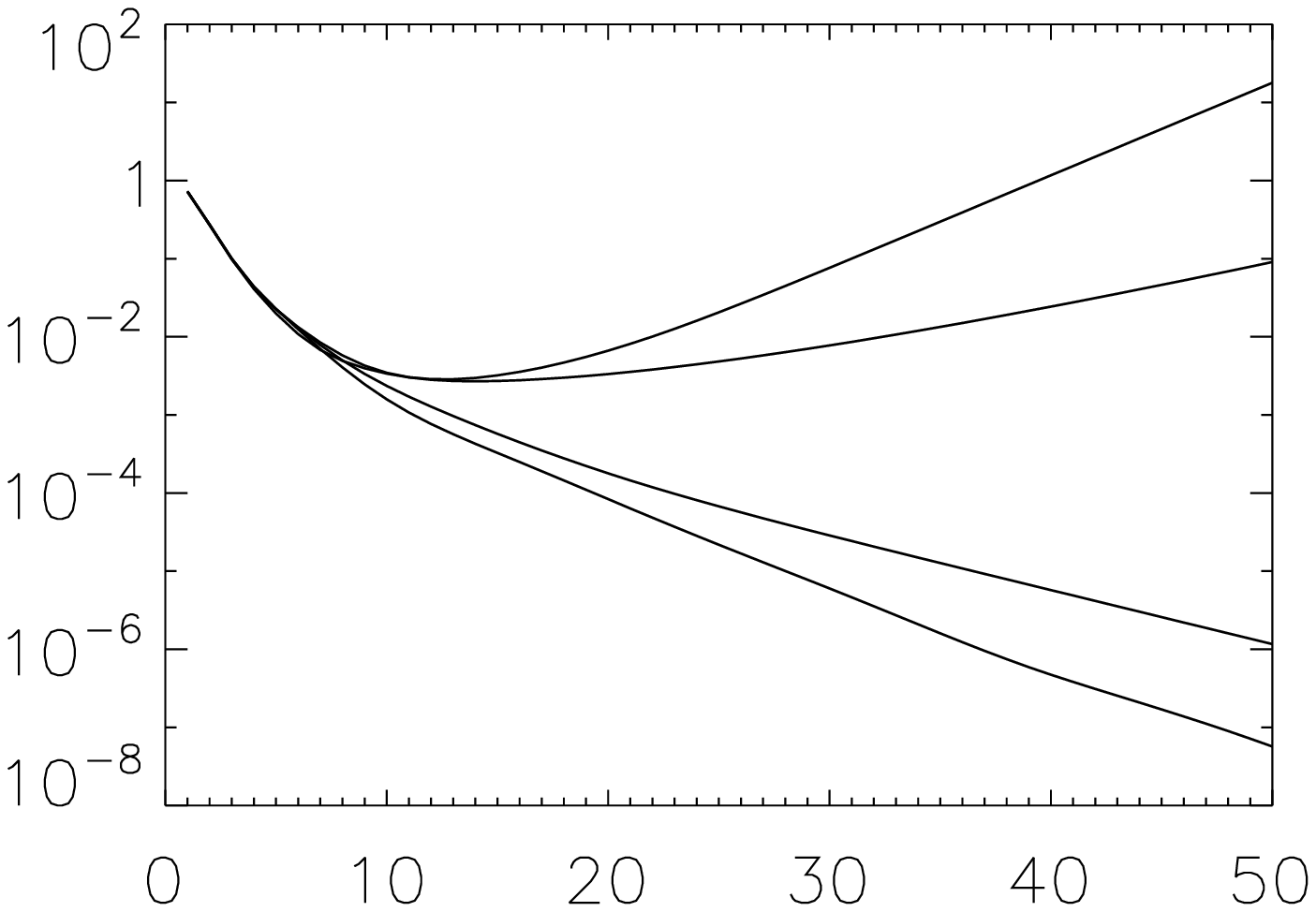,width=67mm}}
\put(4,30){(a)}
\put(76,0){\psfig{file=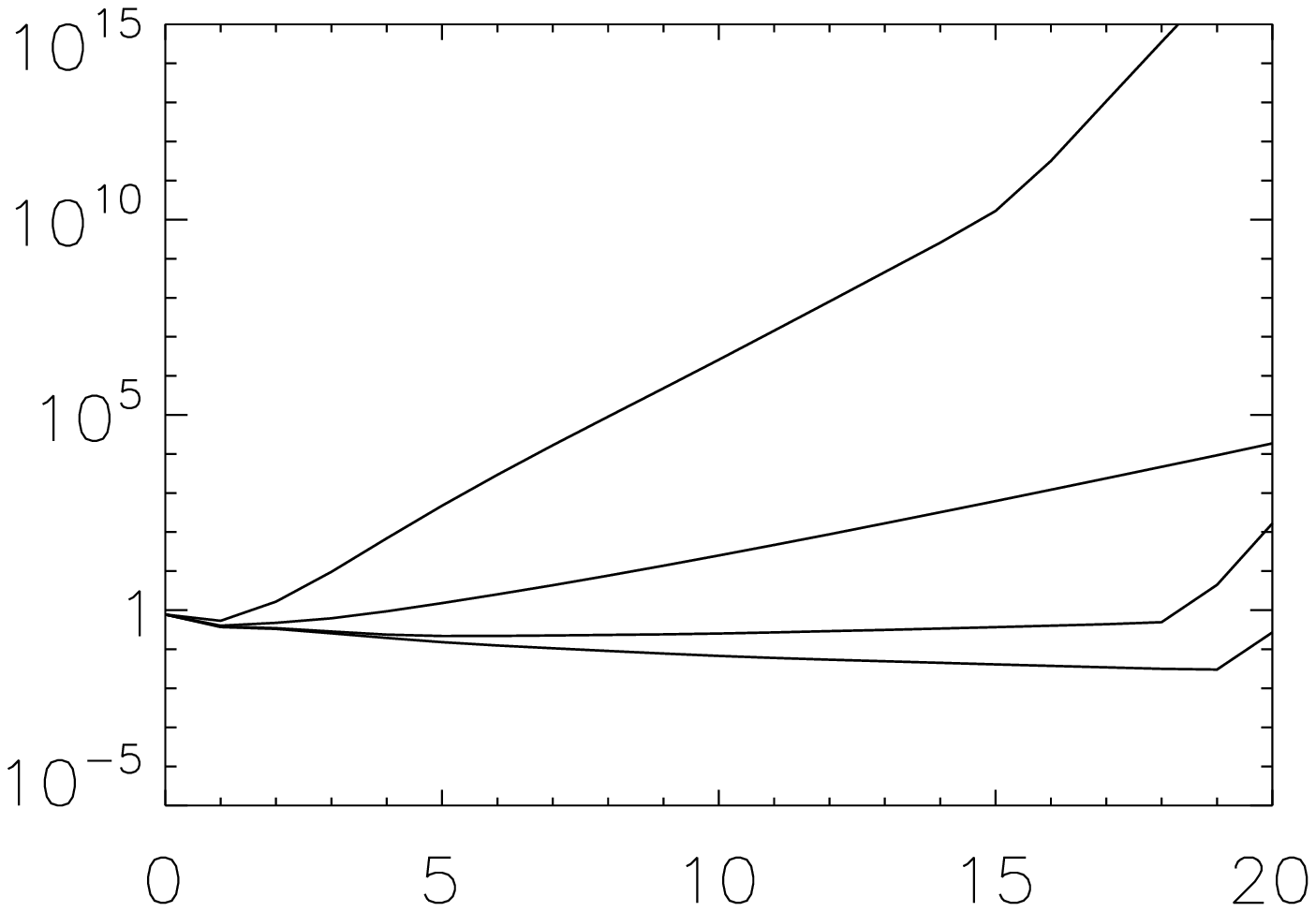,width=67mm}}
\put(70,30){(b)}
\put(17,30){$\|\bxi^{(s)}\|$}
\put(83,30){$\|\omega_s\|$}
\put(53,5){$s$}
\put(119,5){$s$}
\put(43,6){$t=0$}
\put(109,6.5){$t=0$}
\put(49,13){$t=1$}
\put(113.5,11.5){$t=1$}
\put(49,20){$t=3$}
\put(113.5,17.5){$t=3$}
\put(43,28.5){$t=5$}
\put(109,29){$t=5$}
\end{picture}
\caption{Same quantities as in Figs.~\ref{coeft0} (a) and (b), but computed for time-Taylor
expansions around various times (as labelled).}
\label{coeft0135}\end{figure}

We define $C(\bk,{\bm U})\equiv|\bk\cdot{\bm U}\,\Delta t|$. Its maximum
over all wave vectors and velocities present in the flow is the Courant number
${\rm Co}\equiv k_{\max} U_{\max}\Delta t$. Here, $k_{\max}$ is the maximum
wave number, determined by the number of grid points
per spatial period, $N$, and the dealiasing
(typically, $k_{\max}=N/3$), and $U_{\max}$ is the maximum velocity. We observe
that, in the series \rf{expdev}, the moduli of the various terms decrease
with $s$ when ${\rm Co}<1$. Otherwise, for sufficiently small eddies,
swept by sufficiently fast large eddies, i.e., for $C(\bk,{\bm U})>1$,
the moduli increase at first with the order $s$ while $s<C(\bk,{\bm U})$.
Of course, the whole sum, being an imaginary exponential, has a unit modulus,
but this is a sum of a large number of terms, some of which have moduli that
are exponentially large in the Courant number. The accuracy of such a summation
for large Co is low, due to the finite precision of
the computations. This can result in catastrophic loss of
precision -- a well-known phenomenon that can mar numerical summation
of real series with alternating signs (see, e.g., \cite{ac}). To avoid
this difficulty, one can require ${\rm Co}<1$, which is precisely the standard
Courant--Friedrichs--Lewy condition \cite{cfl}. Note that the
CFL criterion for the Eulerian Taylor method does not stem from a
stability condition, but from the finite precision of the computations.

Figure~\ref{coeft0135} shows the L$^2$ norms of time-Taylor coefficients of
various orders, in Lagrangian and Eulerian coordinates, respectively. All
computations are done in quadruple precision with the resolution of $1024^2$
Fourier harmonics. The various graphs are plotted for Taylor expansions
around different output times. Because of the constraints on
time steps that are quite different in Lagrangian and Eulerian\break coordinates,
the number of time steps between successive output times
is quite different in Lagrangian coordinates (typically of the
order of tens) and in Eulerian coordinates (typically of the order of
hundreds). While the Lagrangian and Eulerian time-stepping is performed using
Taylor expansions truncated at the eighth order, at output times they are
pushed to much higher orders. For the Eulerian integrations we observe again
a transition to rounding-noise dominated coefficients at orders roughly
in the range 15--20 (a little beyond order 20 for $t=3$). In the Lagrangian
case this also happens, but at orders about four times larger (not shown in the
figure). Thus, we observe essentially the same phenomenon as for the Taylor
expansions around $t=0$. The explanation is now somewhat different: in the
Eulerian case it is not just the noise that is amplified by the $\i\bk$ factors
stemming from spatial differentiation, but the actual signal, which
can now be much stronger than the direct rounding noise. However,
when we use FFTs,
there is indirect rounding noise, stemming from much lower wave numbers
that have much higher amplitude. Furthermore, in the Eulerian simulations,
this indirect rounding noise gets again amplified by the $\i\bk$ factors
stemming from the space derivatives.

\section{The radius of convergence and depletion}\label{s:depletion}

A well-known advantage of Lagrangian and semi-Lagrangian methods is that
the time step is constrained not by the usual CFL condition
but, roughly, by the inverse of the largest velocity gradient. For our
Cauchy--Lagrangian method, the optimal time step $\Delta t$ for a given
accuracy, when starting from time $t$, is a fixed fraction of the radius
of convergence $\Rei(t)$ of the time-Taylor series for the Lagrangian map
(see Section~\ref{ss:optimal}). It is of interest to know how
$\Rei(t)$ changes in time, since this determines how the time step
should be changed and also the total number of time steps that will be
needed to reach a given termination time $\tf$.
When using initial conditions that have mostly low-wave-number harmonics,
such as the 4-mode flow, the subsequent transfer of excitation
to higher harmonics will be accompanied by some growth in
the velocity gradient: substantially so in 3D and more moderately
in 2D, because of the conservation of vorticity.

\begin{figure}[t!]
\begin{picture}(132,36)(0,0)
\put(10,0){\psfig{file=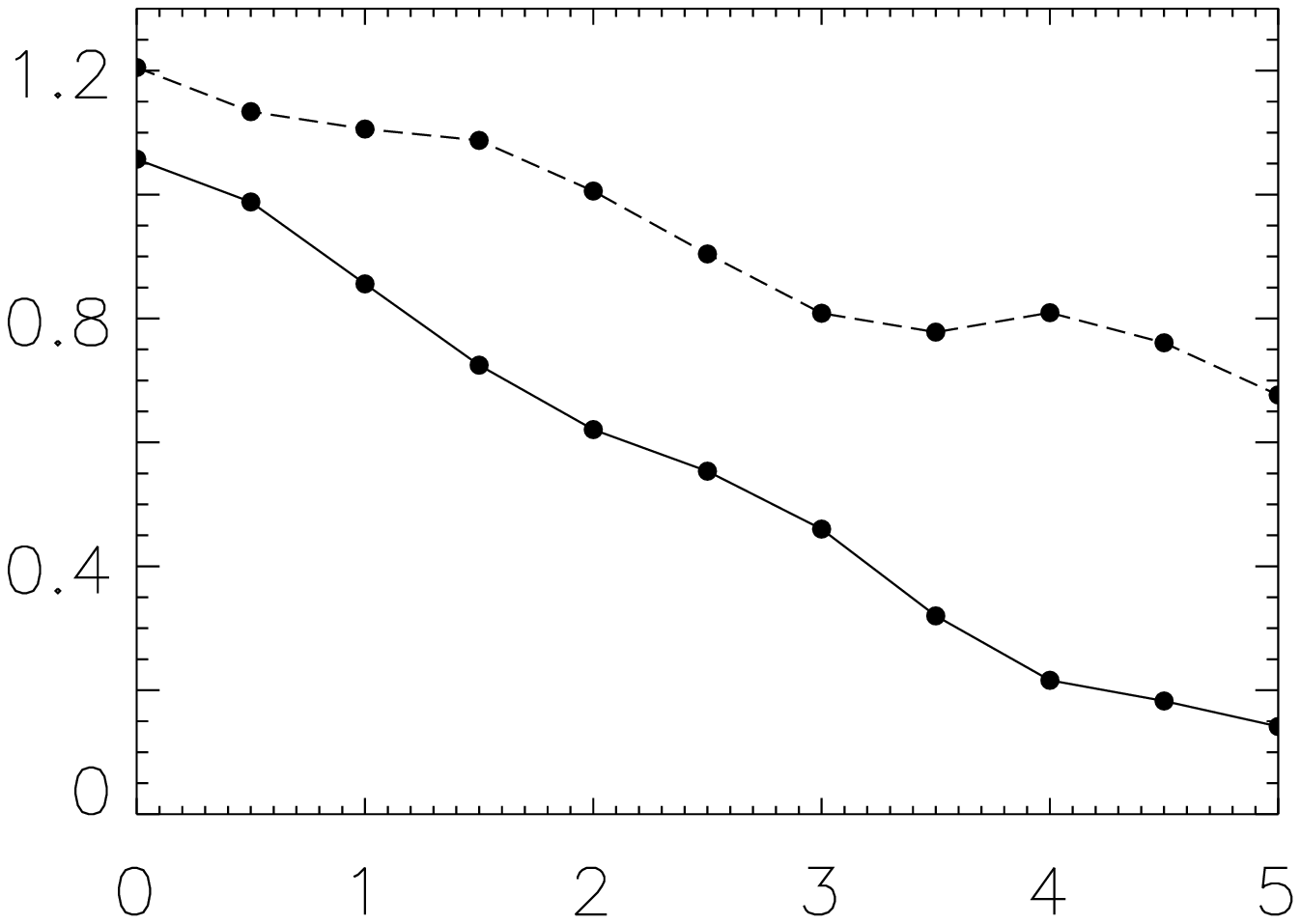,width=67mm}}
\put(4,30){(a)}
\put(76,0){\psfig{file=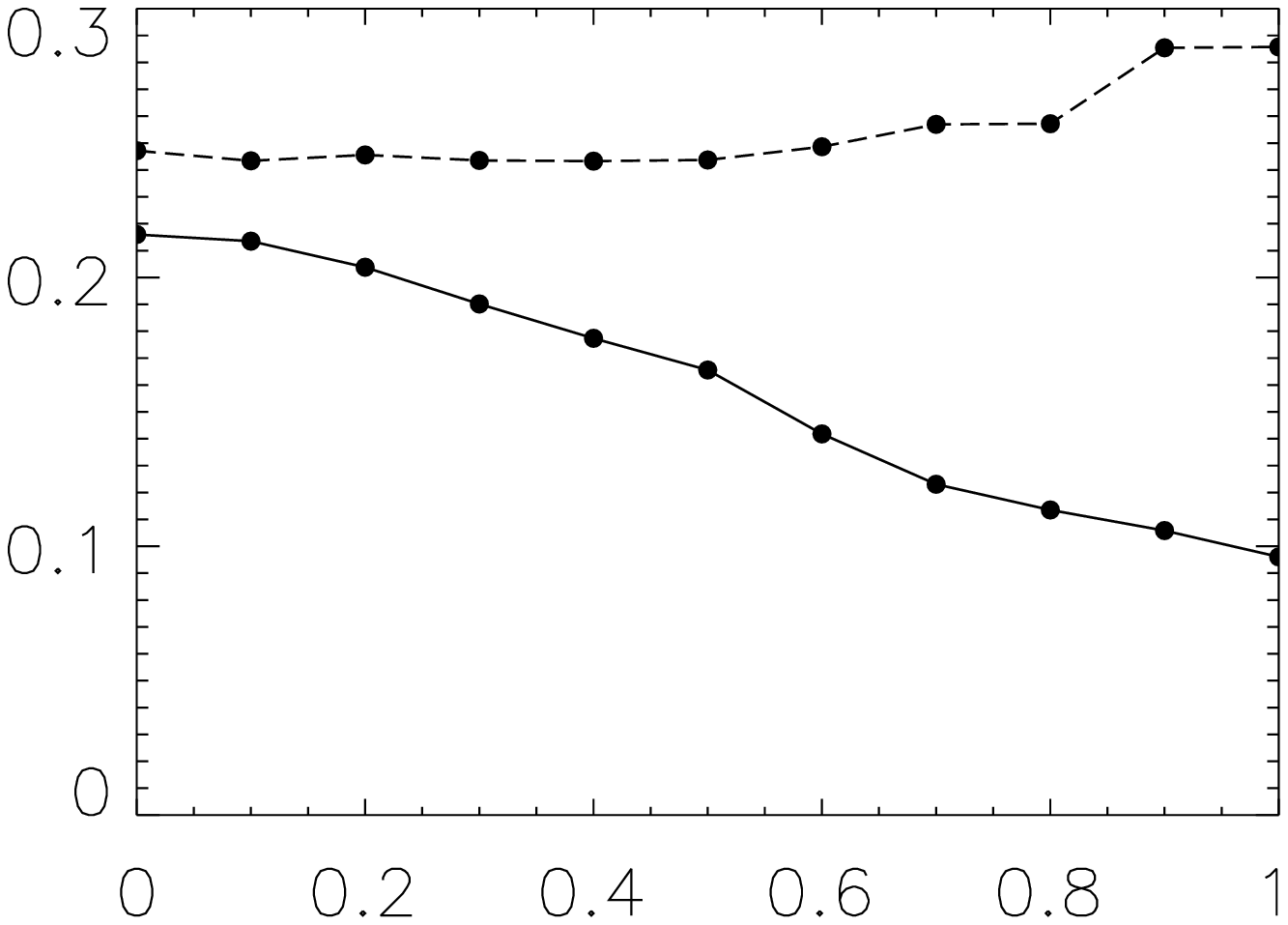,width=67mm}}
\put(70,30){(b)}
\put(16,18){$\Rei(t)$}
\put(82,18){$\Rei(t)$}
\put(49.5,5.5){$t$}
\put(115.5,5.5){$t$}
\end{picture}
\caption{Temporal evolution of the radii of convergence $\Rei(t)$ of the time-Taylor
series of the displacement $\bxi(\ba,t)$ simulated by the Cauchy--Lagrangian method
(dashed lines), and of the vorticity $\omega(\ba,t)$ simulated by the Eulerian
Taylor method (solid lines). Dots show the times at which
the radii were computed. Initial condition: 4-mode (a), random (b).
Resolution: $1024^2$ harmonics, quadruple precision computations.}
\label{rad}\end{figure}

As a consequence, we expect a decrease of the radius $\Rei(t)$.
If $\Rei(t)\to 0$ as $t\to t_\star$ from
below, a blow-up of the solution takes place at the
finite time $t_\star$. Whether the converse holds is so far unknown.
Anyway, in the present numerical studies, we are only concerned with
the 2D case. Before turning to numerical results, let us recall
the proven constraints. Wolibner's theorem \cite{wo} then guarantees
regularity for all $t>0$ when the initial vorticity is
H\"older-continuous of exponent $\alpha_0>0$. More precisely, Wolibner
showed that, at time $t>0$, the vorticity is H\"older-continuous of
exponent $\alpha(t)\ge\alpha_0\,\e^{-t\,\sup|\omega_0|}$, where
$\sup|\omega_0|$ is the supremum of the modulus of the initial vorticity.
The constraint on vorticity conservation does not prevent the velocity gradient
from growing, but it cannot grow faster than $\e^{t\,\sup|\omega_0|}$, a bound
which is actually sharp \cite{KS14}. Consequently, for such 2D
flow, the radius of convergence $\Rei(t)$ has a lower bound proportional
to $\e^{-t\,\sup|\omega_0|}$ (see Section~\ref{s:background}). The actual
numerical results, shown in Fig.~\ref{rad}, indicate that the
decrease in time of the radius of convergence is much slower than
this bound. Consider, for instance, the 4-mode initial
condition, for which the supremum of the initial vorticity is 2.8\,. The
aforementioned lower bound could have the radius of convergence decreasing
by a factor $\e^{-5\times 2.8}\approx 8\times10^{-7}$ over the time
span from 0 to 5. Actually, the radius of convergence decreases only
by about 30\%. As to the random initial condition, its radius does not decrease
at all. For comparison, the radii of convergence for the Eulerian Taylor method
are also shown: they decrease much faster.

This huge discrepancy between the proven bounds and the actual results
is most likely related to the phenomenon of {\it depletion}: in
both 2D and 3D incompressible flow it is frequently found that the
nonlinear effects are growing in time much slower than permitted
by the best proven bounds \cite{fpsm}, probably
because such flow tends to organise itself into structures that have
an almost vanishing nonlinearity, e.g., quasi-one-dimensional flow depending
mostly on a single coordinate. Depletion is definitely why
2D flow is often found numerically to be even more regular than implied
by Wolibner's lower bound, and also why 3D simulated flows show
little reliable evidence for blow-up (one possible exception is
a recently studied swirling axisymmetric flow with a solid boundary \cite{lh}).

\section{Conclusion}\label{s:conclusion}

We have developed a novel Cauchy--Lagrangian (CL) method to solve the initial-value
problem for ideal incompressible flow, governed by the Euler equation.
The method operates in Lagrangian coordinates, relies on temporal analyticity
of particle trajectories in time and employs the recurrence relation \rf{bigXis}
for coefficients of the Lagrangian time-Taylor series, stemming from the Cauchy
invariants. Our experiments in 2D show that this method is highly efficient
(like many semi-Lagrangian methods, it does not suffer from
the Courant--Friedrichs--Lewy constraint) and enables typically an order of magnitude
faster time-marching than the fourth-order Runge--Kutta method and the Eulerian
Taylor method, at the same time being significantly more accurate. As
shown in Section~\ref{ss:efficiency}, the gain in speed can even be higher,
in the range 100--200, at the highest spatial resolutions. Last but
not least, the Cauchy--Lagrangian method allows very-high-order time marching
without the severe rounding error problems that affect high-order
Eulerian time marching. Of course, these impressive
results are so far obtained for just a narrow niche of problems: at present,
incompressible space-periodic two-dimensional flow.

The CL algorithm has been presented
in Section~\ref{s:CLalgorithm} for flow of arbitrary dimension. Since
the CL method allows efficient high-order and high-precision
simulations with only modest growth of rounding errors,
it is a strong candidate for the reliable numerical exploration of the
blow-up problem for the 3D Euler equation (see, e.g., \cite{gi}).

The method, here described for incompressible Euler flow, can be readily
extended for solving certain other problems arising in various areas of physics.
These are described by equations for which known theorems guarantee
the time-analyticity of Lagrangian fluid particle trajectories.
Analyticity in time of the Lagrangian map for at least a finite time, under
the assumption that the initial velocity gradient is, e.g.,
H\"older-continuous, holds for the Euler equation (which applies to inviscid
incompressible fluid flow and to the 2D drift-Poisson equation for an electron
plasma). Furthermore, it also holds for the surface quasi-geostrophic
equation, the incompressible porous medium equation and the Boussinesq
equation. So far, the latter three have been handled not by a recurrence
relation approach, using Cauchy's Lagrangian formulation, but by applying
the Fa\`a di Bruno formula \cite{cs} which gives an explicit but somewhat
cumbersome representation of all the time-Taylor coefficients \cite{cvw}.
Analyticity results, derived from recurrence relations similar to \rf{rrot}--\rf{rgra}, were also proved
for compressible flow of cosmological interest, such as the Euler--Poisson
equations for an Einstein--de Sitter universe \cite{zf}, recently extended
to the more realistic case of a $\Lambda$CDM universe \cite{rvf}.
The Euler--Poisson equations apply to the collisionless flow of dark matter
only as long as multi-streaming (development of caustics with more than
one velocity) is absent or negligible. In such a framework, our methods
are useful in at least two cases: the blow-up problem (called shell-crossing by
cosmologists) and the reconstruction problem, in which the whole dynamical
history of the universe is obtained from just the spatial distribution of
matter at the present epoch \cite{Fetal,Betal}.

But a warning is necessary: not all ideal fluid flow problems are
endowed with time-analytic Lagrangian fluid particle trajectories.
For example, this is apparently not the case of magnetohydrodynamic (MHD)
ideal flow and of barotropic flow.
For non-diffusive MHD flow, one can still write recurrence relations
for the time-Taylor coefficients of the Lagrangian map. However,
the Lagrangian gradient of a time-Taylor coefficient of a higher index
does not involve just the gradients of those of lower index, but also their
space derivatives of order up to three. In other words, in contrast with
the ordinary Euler equation which does not lose derivatives in Lagrangian
coordinates (see Section~\ref{s:rounding}), the MHD equations do, and thus
MHD flow with only limited initial spatial smoothness lacks analyticity
of fluid particle trajectories. Such a loss of derivatives may be viewed
physically as due to the frozenness of the magnetic field. The barotropic
equations govern compressible fluids, in which the pressure is a function
of density only. They play an important role in geophysical fluid dynamics
\cite{mw}. The standard derivation of the Cauchy invariants applies also
to the barotropic equations. However, the second equation, stemming from mass
conservation, is profoundly modified and, in it, loss of derivatives and
hence of analyticity occurs again.

Another class of flows which might be amenable to some variant of the
Cauchy--Lagrangian method is viscous flow (e.g., Navier--Stokes flow) that is
initially analytic in the space variable (if it is not space-analytic
initially, it becomes so at an arbitrary short positive time \cite{ft}).
Such flow is also analytic in time in both Eulerian and Lagrangian
coordinates. There is however a new practical problem: when the equations are
written in Lagrangian coordinates, if we try following Cauchy's approach,
we find that viscous terms involve nonlinearities in
the Lagrangian map of degree four in 2D and of degree 6 in 3D; thus the
recurrence relations become somewhat cumbersome. An alternative and
much simpler strategy for viscous flow is a fractional step method
\cite{yanenko} in which the reversion to Eulerian coordinates by interpolation
is followed by a solution of the viscous heat equation. This is
related to the Trotter-formula approach to Navier--Stokes flow
discussed in \S\,13 of \cite{em}. At high Reynolds numbers, this
method will result in a low-order scheme for the smallest
viscosity-dominated scales but the larger scales will still benefit
from the high accuracy permitted by the CL method.

In this first paper on the Cauchy--Lagrangian method we have
concentrated deliberately on explaining and validating the technique.
We expect that, over the years, this technique will be applied and
improved for a range of situations and problems. It is of course
desirable to handle flow with boundaries: for very smooth boundaries,
this can in principle be done by solving appropriate Helmholtz--Hodge
problems with boundaries. As explained in Section~\ref{s:CLalgorithm}, the basic
ideas for the
Cauchy--Lagrangian method are the same for 2D and 3D flow. At high resolution,
it becomes important to have an efficient parallelisation of both the Lagrangian
time-stepping and the
interpolation phase back to Eulerian coordinates; this is being currently
explored. One of the very challenging applications of such 3D Euler computations
is the possible persistence for all times of an initially assumed smoothness vs.
the loss thereof after a finite time (blow-up problem). In this context, we
finally observe that real-space blow-up cannot take place in
bounded planar 2D flow, but there is already strong evidence that in
2D the Lagrangian time-Taylor series (around $t=0$) has a finite
radius of convergence. Hence, Lagrangian trajectories should have
complex-time singularities, whose nature can be determined by careful
simulations.

\section*{Acknowledgments}
We are grateful to C.~Bardos, A.~Bhatnagar, M.~Blank, J.~van der Hoeven, R.~Pandit,
S.S.~Ray, A.~Shnirelman and B.~Villone for useful discussions.
Research visits of OP and VZ to the Observatoire de la C\^ote d'Azur (France)
were supported by the French Ministry of Higher Education and Research.
Part of this research was supported by the F\'ed\'eration Wolfgang
Doeblin of the Centre National de la Recherche Scienfique under the
grant ``Nouveaux outils math\'ematiques et num\'eriques
pour la structure Eul\'erienne et Lagrangienne des \'ecoulements''.
Some of the calculations were done at the M\'esocentre de
l'Observatoire de la C\^ote d'Azur.

\appendix
\section{Bounds for the space-dependent radii of convergence}\label{a:lpth}

The Lagrangian time-Taylor series \rf{xiSeries} is at the core of the
Cauchy--Lagrangian numerical method. For various reasons explained in the body
of this paper, it is of interest to know its radius of convergence. More
precisely, we want to know how to determine the infimum over the Lagrangian
space of the radii of convergence of the time-Taylor series. A theorem,
announced in Section~\ref{ss:lpth}, allows us to relate this radius to that
of the ordinary Taylor series obtained by replacing each term in \rf{xiSeries}
by its L$^p$ norm. We give now a formal statement and proof of this theorem.

We consider an abstract infinite function series
\BE\xi(\ba,t)=\sum_{s=1}^{\infty}\xi^{(s)}(\ba)t^s,\EE{xiS}
where $\ba$ is in the $d$-dimensional periodicity domain (torus) $\T^d$.
This need not be a Lagrangian solution of the Euler equation; in particular,
it suffices for our purposes to provide arguments for a scalar-valued
series. Let $R(\ba)$ denote the radius of convergence of the time series
at point~$\ba$. We define $\Rei\equiv{\rm ess\,inf}_{\ba}R(\ba)$; we
denote the L$^p$ norm of a function
$\xi(\ba)$ by $\|\xi(\ba)\|_p$ and the radius of convergence of the series
$\sum_{s=1}^\infty\|\xi^{(s)}(\ba)\|_p\,t^s$ by $R_p$.
The Lebesgue space L$^p(\T^d)$ under consideration is that of scalar-valued
functions of the spatial variable $\ba$ of an arbitrary dimension $d$;
we assume that all $\xi^{(s)}(\ba)$ and the sum $\xi(\ba,t)$ are space-periodic.

We prove now the following

\noindent{\it Theorem}. For some $p\ge 1$, let all $\xi^{(s)}(\ba)$ be in L$^p$. Then\\
1. For all $t<R_p$, the sum $\xi(\ba,t)$ belongs to L$^p$, and $\Rei\ge R_p$.\\
2. If $\xi(\ba,t)$ is also in L$^p$ for any $t<R$, then $R_p\ge\Rei$.

\noindent{\it Proof}. The radius of convergence of the series \rf{xiS}, $R(\ba)$,
does not change, when the coefficients $\xi^{(s)}(\ba)$ are replaced by their
absolute values, $|\xi^{(s)}(\ba)|$. Similarly, applying the absolute value
does not modify the Lebesgue space norms. Consequently, without any loss
of generality, we assume $\xi^{(s)}(\ba)\ge0$ for all $s$ and $\ba$, otherwise
replacing the coefficients by their absolute values.

To prove statement 1 of the theorem, we consider the truncated
series $\xi_S(\ba,t)=\sum_{s=1}^S\xi^{(s)}(\ba)t^s$ for a fixed $t$
in the interval $0<t<R_p$. Since the partial sums $\xi_S(\ba,t)$ constitute
a Cauchy sequence in L$^p$, it converges in L$^p$ to a limit function
$\xi(\ba,t)$ also belonging to L$^p$. By a well-known theorem of functional
analysis (see, e.g., \cite{au}), one can choose a subsequence
$\xi_{S_k}(\ba,t)$ that converges almost everywhere to $\xi(\ba,t)$.
(In other words, for our fixed $t$, there exists a set of measure zero
such that, at any point $\ba$ in the complement of this set,
$\lim_{k\to\infty}\xi_{S_k}(\ba,t)=\xi(\ba,t)$, the limit being here understood
as the usual limit for infinite sequences of real numbers.) However,
since $\xi^{(s)}(\ba)\ge0$ for all $s\ge1$, the partial sums
$\xi_S(\ba,t)$ grow monotonically at each point $\ba$ (for the fixed $t$)
on increasing the index $S$, and hence the entire sequence $\xi_S(\ba,t)$
converges for $S\to\infty$ almost everywhere to $\xi(\ba,t)$. Therefore,
any fixed time $t<R_p$ is within the disk of convergence of the series
\rf{xiS}, i.e., $R(\ba)\ge R_p$ for almost all $\ba$, which implies
the first statement of the theorem.

We prove now statement 2. For any $p\ge1$ and any
quantities $\alpha_s\ge0$ such that $\sum_{s=1}^\infty\alpha_s=1$,
clearly, $\sum_{s=1}^\infty\alpha^p_s\le1$. Taking\rule{0ex}{2.4ex}
$\alpha_s=\xi^{(s)}(\ba)t^s/\xi(\ba,t)$ (recall $\xi^{(s)}(\ba)\ge0$ for all
$s$), we thus find
\BE\sum_{s\ge1}(\xi^{(s)}(\ba)\,t^s)^p\le\xi^p(\ba,t)\EE{inp}
for almost all $\ba$ and any $t$ such that $0<t<R(\ba)$, and hence
for any $0<t<\Rei$. Integrating
this inequality over the periodicity domain yields\rule{0ex}{2.4ex}
$\sum_{s\ge1}\|\xi^{(s)}(\ba)\|^p_p\,t^{sp}\le\|\xi(\ba,t)\|^p_p$.
For \hbox{$t>R_p=(\lim\,\sup_{s\to\infty}\|\xi^{(s)}(\ba)\|_p^{1/s})^{-1}$,\rule{0ex}{2.4ex}}
the series in the l.h.s.~of the latter inequality diverges,\rule{0ex}{2.4ex}
but it converges for all $t<\Rei$, for which the inequality holds. Consequently,
$\Rei\le R_p$. This concludes the proof of the theorem.

It remains to comment on the significance of the condition in statement~2
that $\xi(\ba,t)$ belongs to L$^p$ for all $t<R$. One can indeed
construct an abstract counterexample showing that the radii are not necessarily
equal, if this condition is not satisfied.

{\it Counterexample}. For simplicity, we consider functions defined
in the interval $[0,1]$. Suppose $1=a_0>a_1>...>a_s>a_{s+1}>...>0$
are real numbers such that $\lim_{s\to\infty}a_s=0$. For $\alpha>1/p$
and $s\ge1$, consider functions $\xi^{(s)}(a)$, which vanish outside
the interval $a_{s-1}\ge a>a_s$, and $\xi^{(s)}(a)=a^{-\alpha}$ in this interval.
Evidently, the series \rf{xiS} converges for any $t$, i.e., $\Rei=\infty$.
The sum \rf{xiS} belongs to the Lebesgue space L$^p([0,1])$ only for $|t|<R_p$,
since by construction $\|\xi(a,t)\|^p_p=\sum_{s\ge1}\|\xi^{(s)}(a)\|^p_p\,t^{sp}$.
By choosing a suitable sequence of nodes $a_s$, we generate functions
$\xi^{(s)}(a)$ of prescribed norms $\|\xi^{(s)}(a)\|_p=\gamma_s$. Since
$\|\xi^{(s)}(a)\|_p=(p\alpha-1)^{-1/p}(a^{1-p\alpha}_s-a^{1-p\alpha}_{s-1})^{1/p}$,
the condition $\lim_{s\to\infty}a_s=0$ is satisfied whenever $\sum_{s\ge1}\gamma^p_s=\infty$;
otherwise, $\gamma_s\ge0$ are arbitrary. For instance, $\gamma_s=\gamma^s$
is permissible for any $\gamma\ge1$. We have thus constructed an example of series
\rf{xiS}, non-summable in L$^p([0,1])$, for which $\Rei=\infty$, while
the accompanying series \rf{Tp} constructed of L$^p$ norms possesses
an arbitrary prescribed finite radius of convergence $R_p=1/\gamma\le1$;
for such a series, only the first statement of the theorem holds.

{\it Remark}. The Lebesgue space L$^p(\T^d)$ in the statement of the Theorem
can be replaced by a Sobolev function space W$^p_q(\T^d)$ for $p\ge1$ and
$q>0$. Its norm is defined as
$\|\xi(\ba)\|_{q,p}=\|(I-\nabla^2)^{q/2}\xi(\ba)\|_p$,
where $I$ is the identity operator and $\nabla^2$ the Laplacian (the spaces
for $p=2$ are often denoted by H$^q(\T^d)$). The quantity $R_p$ is then
replaced in the modified theorem by $R_{q,p}$, the radius of convergence
of the series $\sum_{s=1}^\infty\|\xi^{(s)}(\ba)\|_{q,p}\,t^s$,\rule{0ex}{2.5ex}
and instead of $\Rei$ we use $R_{q,{\rm inf}}$, the essential infimum over
space (i.e., over the variable $\ba$) and over all $q'$, such that
$0\le q'\le q$, of the radii of convergence of the series
\BE\sum_{s=1}^\infty(-\nabla^2)^{q'/2}\xi^{(s)}(\ba)\,t^s.\EE{sqp}
Convergence of the series
$\sum_{s=1}^\infty\|\xi^{(s)}(\ba)\|_{q,p}\,t^s$\rule{0ex}{2.5ex} implies
convergence of \rf{sqp} in L$^p$ for all $q'\le q$
and, by the arguments in the proof of the Theorem, the inequality
$R_{q,{\rm inf}}\ge R_{q,p}$. Conversely, convergence of \rf{xiS} in W$^p_q$
implies convergence of \rf{sqp} in L$^p$ for $q'=0$ and $q'=q$, and
by the second part of the Theorem the radii of convergence of the two
series, $\sum_{s=1}^\infty\|\xi^{(s)}\|_p\,t^s$,\rule{0ex}{2.5ex} and
$\sum_{s=1}^\infty\|(-\nabla^2)^{q/2}\xi^{(s)}\|_p\,t^s$,\rule{0ex}{2.5ex}
are both bounded from below by $R_{q,{\rm inf}}$, whereby
$R_{q,{\rm inf}}\le R_{q,p}$. Note that $R_{q,p}$ and $R_{q,{\rm inf}}$
are monotonically decreasing functions of $q$ for any $p\ge 1$. The
inequality $R_{q,p}<R_{q',p}$ for some $q>q'$ signals existence of a ring
within the disk of convergence in W$^p_{q'}$, where the series \rf{xiS}
does not converge in W$^p_q$. An intriguing problem is to establish whether
for a solution of the Euler equation $R_{q,{\rm inf}}=R_{q',{\rm inf}}$ and
$R_{q,p}=R_{q',p}$ hold for all sufficiently small $q\ne q'$.

\section{Computation of the radius of convergence of a power series\\
with positive coefficients}\label{a:conver}

The implementation of the Cauchy-Lagrangian method does not require
the evaluation of the radius of convergence of the time-Taylor series
employed. However, to test the method, to analyse its performance and
interpret some of the results (Sections
\ref{s:testing}-\ref{s:depletion}), we need to compute the radii of
convergence of the time-Taylor series for both CL and ET
algorithms. For estimating these, using the result of \ref{a:lpth}, we
can replace the full functional Taylor series
$\sum_{s=1}^\infty\bxi^{(s)}(\ba)\,t^s$ by the ordinary power series
$\sum_{s=1}^\infty\|\bxi^{(s)}(\ba)\|\,t^s$, where $\|\cdot\|$ denotes
the L$^2$ norm. To obtain reasonable precision on the radii, say, a few
percent, this computation actually requires the use of quadruple precision and
the knowledge of enough Taylor coefficients, usually more than is required
for the time-stepping. Due to this added complexity, we do not compute the
radii at every time step, but just frequently enough to be able to
monitor how the radii evolve in time, e.g., every tenth time step.

In this Appendix we describe how we calculate the radii for such power series.

Here is the problem that we want to solve:\\
Given a finite number, $S$, of coefficients of the power series
\BE f(t)=\sum_{s=1}^\infty f_st^s,\EE{pseries}
where $f_s\ge0$, we want to estimate
its radius of convergence, which we denote by
$R$. For our purposes, an accuracy of a few percent suffices.

To begin with, we recall some classical results. By the Cauchy--Hadamard theorem,
\BE\frac{1}{R}={\limsup_{s\rightarrow\infty}\sqrt[s]{|f_s|}}.\EE{hadamard}
A fast way of computing the radius of convergence is the ratio formula
\BE\frac{1}{R}=\lim_{s\rightarrow\infty}\left|\frac{f_{s+1}}{f_s}\right|,\EE{hadrat}
when the limit exists.

With less than a hundred coefficients of the time-Taylor series at hand,
we can obtain the radius of convergence
with a reasonable accuracy using such formulae only, if
we know several terms in the large-$s$ asymptotic expansion of
the coefficients $f_s$. We should derive the asymptotics
of $f_s=\|\bxi^{(s)}(\ba)\|$ from the recurrence relation
\rf{bigXis} for $\bxi^{(s)}$. Unfortunately, no such results
are available so far. To proceed, we look at our problem from an entirely
different angle and ask ourselves the question:
Does the analytic function defined by the series
$\sum_{s=1}^\infty\|\bxi^{(s)}(\ba)\|\,t^s$ have the
``generic'' simple structure near the disk of convergence,
such as many analytic functions \rf{pseries} with real positive
coefficients $f_s$ have?

A standard result is that the radius of
convergence of a Taylor series is the distance in the complex
$t$-plane to the nearest singularity(ies). If there is a single such
singularity and the Taylor coefficients are all positive, the
singularity is on the real positive axis at $t_\star=R$. The simplest
case is when the singularity has a simple structure, e.g., a pole or
a power-law branch point that can be characterised by a local behaviour
$\propto(t_\star-t)^\rho$, where the exponent $\rho$ is a negative
real number. It may then be shown that the Taylor coefficients of large orders
$s$ have the following behaviour:
\BE f_s\approx\gamma\,s^\alpha\e^{\beta s}\quad\hbox{for }s\to\infty,\EE{asimp}
where $\beta=-\ln R$ and $\alpha=-\rho-1$. (The same kind of asymptotics
applies to Fourier series; see, e.g., \cite{ckp}.) In such a case, if we are
able to extract the leading and subleading behaviour of $f_s$,
we may then infer not only the radius of convergence, but also the type
of the nearest singularity on the real axis. One method is the Domb--Sykes
plot of $f_s/f_{s-1}$ versus $1/s$, where the $y$-intercept gives $1/R$ and the
slope near the origin gives the exponent $\alpha$ \cite{domb,hinch}.

We have applied all the standard methods, recalled above. We found
that, when using 50--80 quadruple-precision Taylor coefficients from our
L$^2$ series, such methods allow a determination of the radius of
convergence with an error in the range 10\%--20\%.
We also tried more advanced extrapolation techniques,
such as convergence acceleration methods (see,
e.g., \cite{pa,pa07} and references therein) and the asymptotic
extrapolation technique \cite{jvdh09,pa,pa07}, which have the
potential of capturing both leading and subleading asymptotic behaviour.
Unfortunately, all these techniques failed to reach the asymptotic regimes
and thus gave no improvement over the more standard ones.

We found that for the L$^2$ series with 50-80 terms a rather straightforward
method apparently gives a better accuracy. We again assume that the L$^2$
series has a single singularity and that \rf{asimp} holds asymptotically
for large $s$. It follows that the logarithms of the Taylor
coefficients, $F_s\equiv\ln f_s$, are given by
\BE F_s\approx\ln\gamma+\alpha\ln s+\beta s\qquad\hbox{for~}s\to\infty.\EE{asimf}
We can determine the coefficients $\alpha$, $\beta$ and $\gamma$ by
a least-square fit over the available Taylor coefficients. That is, we minimise
$\sum_{s=1}^Sd_s^2$ for the discrepancies
\BE d_s\equiv c+a\ln s+bs-F_s.\EE{disc}
\vspace*{-\baselineskip}

\begin{figure}[t!]
\begin{picture}(132,30)(0,0)
\put(10,0){\psfig{file=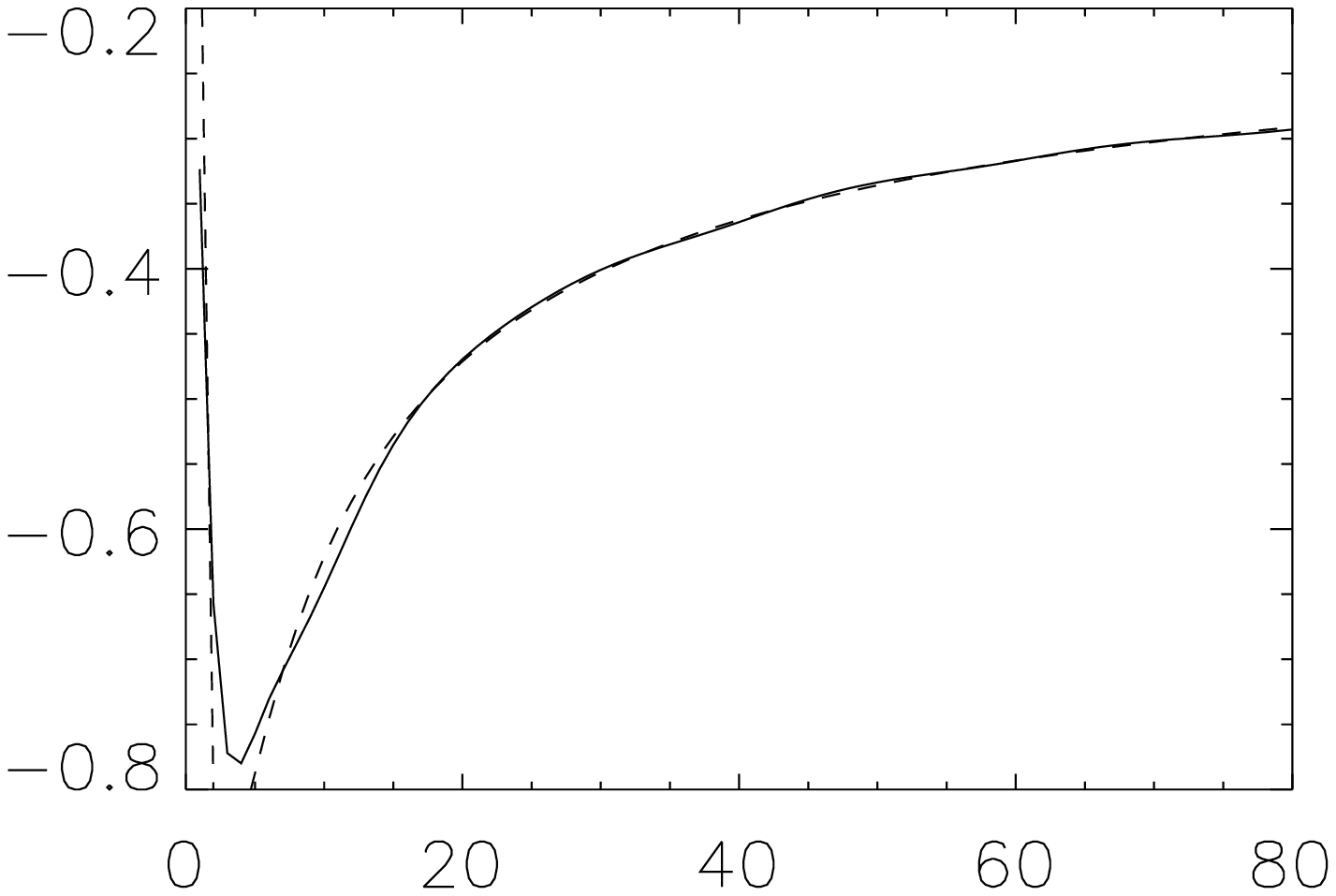,width=67mm}}
\put(4,26){(a)}
\put(76,0){\psfig{file=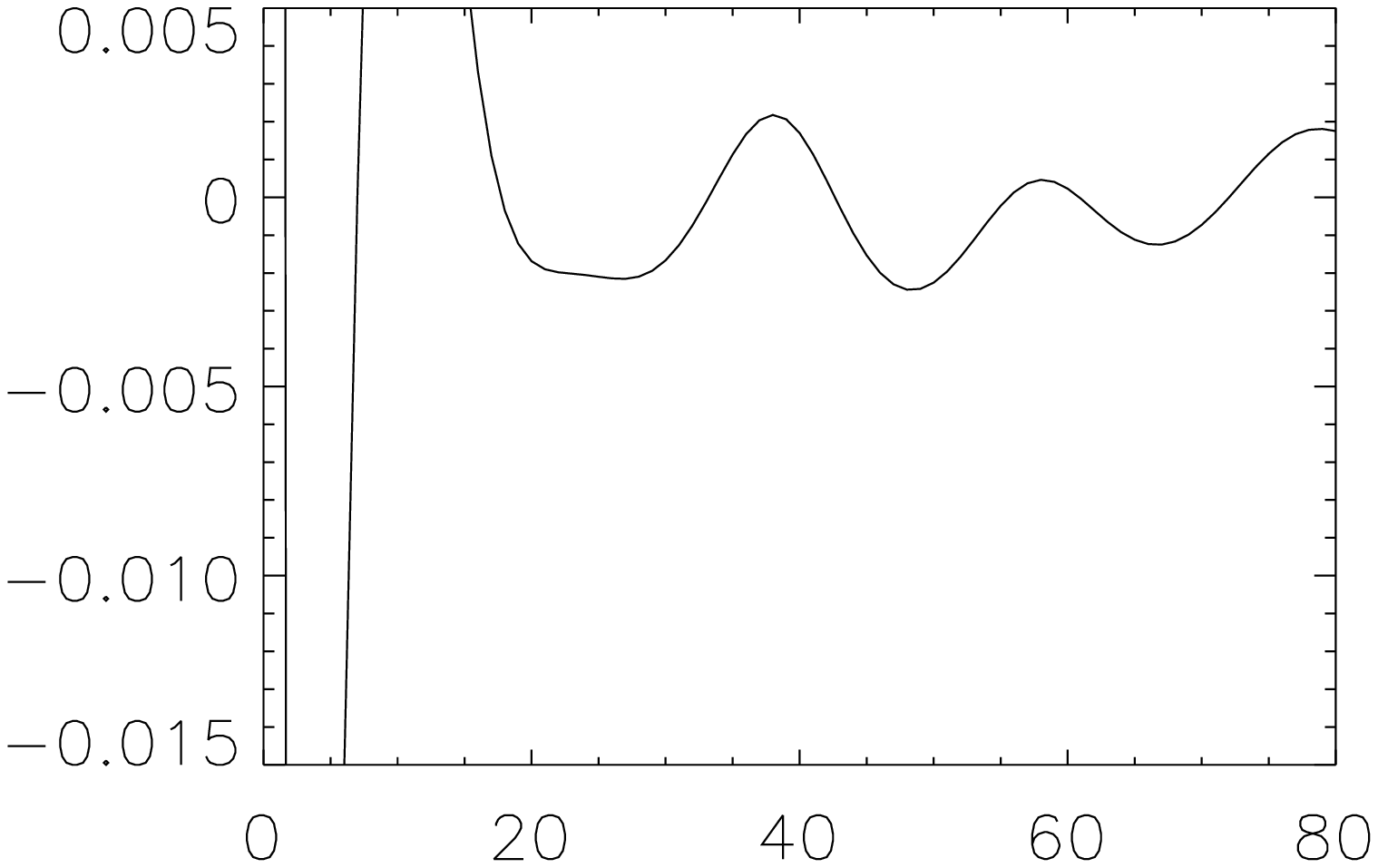,width=67mm}}
\put(70,26){(b)}
\put(51,28){$F_s/s$}
\put(117,28){$d_s/s$}
\put(48,5){$s$}
\put(114,5){$s$}
\end{picture}
\caption{The quantity $F_s/s$ (solid line) and its least-square fit
(dashed line, see \rf{asimf}) (a). Scaled discrepancy $d_s$ \rf{disc} (b).}
\label{figfit}\end{figure}

We now give an example of such determination for the four-mode initial
condition \rf{4mode}. Here $f_s$ are the L$^2$ norms of the first 80
time-Taylor coefficients $\bxi^{(s)}$ for the displacement, computed at $t=0$.
Figure~\ref{figfit}~(a) shows the values $F_s/s$ and the fitted
curve $(c+a\ln s+bs)/s$. (Note that $F_s/s\approx\beta=-\ln R$
for large $s$.) The least square fit gives
$a=-1.94,\ b=-0.187$ and $c=0.116$. Thus, the radius of convergence is
$R=\e^{-b}\approx1.2$\,. To estimate the error on the radius $R$, we show in
Fig.~\ref{figfit}~(b) the scaled discrepancy $d_s/s$ for large $s$. We see
that, beyond $s=20$, the fluctuations do not exceed 1\%. Similar
determinations of the radius $R(t)$ for the same initial four-mode
flow and $0<t<5$ have been made, using between 80 and 50 Taylor
coefficients, yielding errors between 1\% and 3\%.

We have used the same fitting technique to measure the radius $\delta$
of the cylinder of spatial analyticity around the real domain (in both
Lagrangian and Eulerian coordinates). The analogue of \rf{asimp} is then
the following expression for the shell-averaged energy spectrum $E(K)$,
\BE E(K)\approx C K^n\e^{-2\delta K}.\EE{enspecfit}

How can we achieve a much better accuracy on radii of
convergence? We need to know a lot more than 80 coefficients, say, a
few hundred. Unfortunately, as observed in Section~\ref{s:rounding},
in quadruple precision, rounding noise becomes explosively large
just above order 80. We can go a little further by avoiding FFTs and
performing the Fourier-space convolutions explicitly, but this is
computationally very expensive. Otherwise, we have to resort to higher
than quadruple precision. Recent algorithmic developments on FFTs in ``medium''
precision (up to 400 bits) offer interesting possibilities \cite{HL14}.
\end{document}